\documentclass{amsart}
\usepackage{amsmath}%
\usepackage{amsfonts}%
\usepackage{amssymb}%
\usepackage{graphicx}
\theoremstyle{plain}

\numberwithin{equation}{section}
\begin{document}
\title{Expansion and Improvement of Sieve and application in Goldbach's problem}
\author{Cheng Hui Ren}
\email[Chenghui.Ren]{chenghui.ren@gmail.com}%
\date{April, 2009}
\dedicatory{}

\begin{abstract}
This paper expands and improves on the general Sieve method. This expaned and improved Sieve is applied to Goldbach's problem. A new estimate of the exception set in Goldbach's number $E(X)$, an improved lower bound $D_{1,2}(N)$ and upper bound $D(N)$ are proposed.
The proposed values are: $\left|E(X)\right|\leq X^{0.702+\epsilon}$, $D_{1,2}(N)\geq 2.27\frac{C(N)}{\ln^{2}(N)}$, $D(N) \leq 6.916\frac{C(N)}{\ln^{2}(N)}$. \\

keywords: Sieve, Goldbach's number, Goldbach's Conjecture, Twin primes Conjecture. \\

Article class math.NT
\end{abstract}
\maketitle

\section{Introduction}

\noindent \\

\  The modern sieve method began from Brun's fundamental work in 1915. This remained the dominant work until 1941 when Ju. V. Linnik expanded on the sieve method. In 1947, A. Selberg advanced the sieve method by using his upper bound sieve. His sieve weights, are fundamentally different from Brun's thus he brought a structural change into the sieve method. The linear sieve was developed by Rosser in 1950  and Jurkat and Richert in 1965; so far this is the best method of the linear sieve. \\
  
\  This paper expands the sieve function, in this expanded sieve function we can use some weighted sieve results for iteration and improvement of the traditional linear sieve. Applying this sieve method to Goldbach's problem and Twin primes problem, a new estimate of the exception set in Goldbach's number is  obtained; and improved lower bound of $D_{1,2}(N)$ and upper bound of $D(N)$.


\hspace*{\fill}\\
\section{Part I: Expansion and Improvement of Sieve}
\ Set $\mathcal{P}$ is related to parameter $N$ a natural number ,\\

\ Suppose

\begin{equation}
 \mathcal{P}=\mathcal{P}(N) :=\left\{p: p \ \ are \ \ prime \ \ numbers, \ \ (p,N)=1\right\}; \
\end{equation}

\begin{equation}
 \mathcal{P}(z)=\mathcal{P}_{N}(z) =\prod_{p\in \mathcal{P}(N),p<z} p \
\end{equation}
\ usually we omit the parameter $N$.

 Define $\mathcal{A}$ as a number set and
\[ \mathcal{A}_{d}=\sum_{n\in \mathcal{A}, d|n}1 =X\frac{\omega(d)}{d}+r_{d}\leq O(\frac{X\ln^{c_{1}}(X)}{d})\]

\ Here $\omega(d)$ is a multiplicative function, and $\omega(d)$ depends on both $\mathcal{A}$ and $\mathcal{P}$. We suppose $\omega(d)$ fulfils the following expression.
\begin{equation}
\left|\sum_{w \leq p < Z}\frac{\omega(p)\ln(p)}{p}-\ln(\frac{z}{w})\right|\leq L_{2},\   \ 2\leq w\leq z
\end{equation}

\ Define function $\Lambda(n,z)$ and $\Lambda_{\alpha}(n,z,\xi)$, $\alpha > 0$
\\
$$
 \Lambda(n,z)=\begin{cases}
	\begin{array}{ll}
			1 & (n,\mathcal{P}(z))=1 \\
			0 & (n,\mathcal{P}(z))>1 \\
	\end{array}
	\end{cases}
$$

$$
 \Lambda_{k,\alpha}(n,z,\xi)=\begin{cases}
	\begin{array}{ll}
	\xi & (n,\mathcal{P}(z))=1 \\
	k \ln^{\alpha}(q)& (n,\mathcal{P}(z))=q, \ \ is \ \ a \ \ prime.\\
	0 & (n,\mathcal{P}(z)) \ \ have \ \ more \ \ then \ \ 1 \ \ prime \ \ facts.
	\end{array}
	\end{cases}
$$

\ The traditional sieve function is defined as
\begin{equation} 
	S_{0}(\mathcal{A};\mathcal{P},z)=\sum_{n \in \mathcal{A} \atop (n,\mathcal{P}(z))=1,}1=\sum_{n \in \mathcal{A}}\Lambda(n,z)
\end{equation}
\ The expanded sieve function is defined as 
\begin{equation}
 S_{k,\alpha}(\mathcal{A};\mathcal{P},z,\ln^{\alpha}(\xi))=\sum_{n\in \mathcal{A}}\Lambda_{k,\alpha}(n,z,\ln^{\alpha}(\xi))
\end{equation}
\ It is easy to see that we have the following relationship between the traditional and expanded sieves.
\begin{equation}
\ S_{k,\alpha}(\mathcal{A};\mathcal{P},z,\ln^{\alpha}(\xi^{2}))=\ln^{\alpha}(\xi^{2})S_{0}(\mathcal{A};\mathcal{P}(N),z)+k \sum_{2 \leq p<z}\ln^{\alpha}(p)S_{0}(\mathcal{A}_{p};\mathcal{P}(p),z)
\label{L1} 
\end{equation}

\subsection{Lemma 1.1}
\ Suppose $z>z_{1}>2$ \\

\ The well-known identity of Buchstab.
\begin{equation}
\ S_{0}(\mathcal{A};\mathcal{P},z)=S_{0}(\mathcal{A};\mathcal{P},z_{1})-\sum_{z_{1}\leq p<z}S_{0}(\mathcal{A}_{p};\mathcal{P},p)
\end{equation}

\ The opposite identity of the expanded sieve function, 
\begin{equation}
\ S_{k,\alpha}(\mathcal{A};\mathcal{P},z,\ln^{\alpha}(\xi^{2}))=S_{k,\alpha}(\mathcal{A};\mathcal{P},z_{1},\ln^{\alpha}(\xi^{2}))-\sum_{z_{1}\leq p<z}S_{k,\alpha}(\mathcal{A}_{p};\mathcal{P},p,\left(\ln^{\alpha}(\xi^{2})-k \ln^{\alpha}(p)\right))
\end{equation}
\[ -k\sum_{z_{1}\leq p<z}\ln^{\alpha}(p)S_{0}(\mathcal{A}_{p^{2}};\mathcal{P},p) \]
\ And when $\ln^{\alpha}(\xi^{2})>k\ln^{\alpha}(z)$ we have
\[ S_{k,\alpha}(\mathcal{A};\mathcal{P},z,\ln^{\alpha}(\xi^{2}))=\ln^{\alpha}(\xi^{2})S_{0}(\mathcal{A};\mathcal{P},z)-\sum_{2\leq p<z}\frac{\ln^{\alpha}(\xi^{2})-k \ln^{\alpha}(p)}{\ln^{\alpha}(\frac{\xi^{2}}{p})}S_{k(p),\alpha}\left(\mathcal{A}_{p};\mathcal{P},p,\ln^{\alpha}(\frac{\xi^{2}}{p})\right) \]
\[ -k\sum_{z_{1}\leq p<z}\ln^{\alpha}(p)S_{0}(\mathcal{A}_{p^{2}};\mathcal{P},p) \]
\ where 

\[ k(p)=k\frac{\ln^{\alpha}(\frac{\xi^{2}}{p})}{\ln^{\alpha}(\xi)-k \ln^{\alpha}(p)}\geq 0 \]

\begin{proof}
\ Only proof of equation (2.8), by (2.7)
\[ S_{0}(\mathcal{A};\mathcal{P},z)=S_{0}(\mathcal{A};\mathcal{P},2)-\sum_{2\leq p<z}S_{0}(\mathcal{A}_{p};\mathcal{P},p) \]
\[ \ln^{\alpha}(\xi^{2})\sum_{2\leq p<z}S_{0}(\mathcal{A}_{p};\mathcal{P},p)=\sum_{2\leq p<z}(\ln^{\alpha}(\xi^{2})-k \ln^{\alpha}(p))S_{0}(\mathcal{A}_{p};\mathcal{P},p)+k \sum_{2\leq p<z}\ln^{\alpha}(p)S_{0}(\mathcal{A}_{p};\mathcal{P},p) \]
\ The second sum on the right
\[ \sum_{2\leq p<z}\ln^{\alpha}(p)S_{0}(\mathcal{A}_{p};\mathcal{P},p)=\sum_{2\leq p<z}\ln^{\alpha}(p)S_{0}(\mathcal{A}_{p};\mathcal{P},z) \]
\[ +\sum_{2\leq p<q<z}\ln^{\alpha}(p)S_{0}(\mathcal{A}_{pq};\mathcal{P}(p),q) \]
\[ +\sum_{2\leq p<z}\ln^{\alpha}(p)S_{0}(\mathcal{A}_{p^{2}};\mathcal{P}(p),p) \]
\ Thus
\[ S_{k,\alpha}(\mathcal{A};\mathcal{P},z,\ln^{\alpha}(\xi^{2}))= \ln^{\alpha}(\xi)S_{0}(\mathcal{A};\mathcal{P},z)+\sum_{2\leq p<z}k \ln^{\alpha}(p)S_{0}(\mathcal{A}_{p};\mathcal{P}(p),z) \]

\[ =\ln^{\alpha}(\xi)S_{0}(\mathcal{A};\mathcal{P},2)-\sum_{2\leq p<z}(\ln^{\alpha}(\xi)-k \ln^{\alpha}(p))S_{0}(\mathcal{A}_{p};\mathcal{P},p) \]
\[ -\sum_{2\leq q<p<z}k \ln^{\alpha}(q)S_{0}(\mathcal{A}_{qp};\mathcal{P}(q),p) -k\sum_{2\leq p<z}\ln^{\alpha}(p)S_{0}(\mathcal{A}_{p^{2}};\mathcal{P},p) \]

\[ =\ln^{\alpha}(\xi^{2})S_{0}(\mathcal{A};\mathcal{P},2)-\sum_{2\leq p<z}S_{k,\alpha}\left(\mathcal{A}_{p};\mathcal{P},p,\left(\ln^{\alpha}(\xi)-k \ln^{\alpha}(p)\right)\right)-k\sum_{2\leq p<z}\ln^{\alpha}(p)S_{0}(\mathcal{A}_{p^{2}};\mathcal{P},p) \]
\ When $\ln^{\alpha}(\xi^{2})>k\ln^{\alpha}(z)$ we have
\[ S_{k,\alpha}(\mathcal{A};\mathcal{P},z,\ln^{\alpha}(\xi^{2}))=\ln^{\alpha}(\xi^{2})S_{0}(\mathcal{A};\mathcal{P},2)-\sum_{2\leq p<z}\frac{\ln^{\alpha}(\xi)-k \ln^{\alpha}(p)}{\ln^{\alpha}(\frac{\xi^{2}}{p})}S_{k(p),\alpha}\left(\mathcal{A}_{p};\mathcal{P},p,\ln^{\alpha}(\frac{\xi^{2}}{p})\right) \]
\[ -k\sum_{2\leq p<z}\ln^{\alpha}(p)S_{0}(\mathcal{A}_{p^{2}};\mathcal{P},p) \]
\ So we have
\[ S_{k,\alpha}(\mathcal{A};\mathcal{P},z_{1},\ln^{\alpha}(\xi^{2}))-S_{k,\alpha}(\mathcal{A};\mathcal{P},z,\ln^{\alpha}(\xi^{2})) \]
\[ =\sum_{z_{1}\leq p<z}\frac{\ln^{\alpha}(\xi)-k \ln^{\alpha}(p)}{\ln^{\alpha}(\frac{\xi^{2}}{p})}S_{k(p),\alpha}\left(\mathcal{A}_{p};\mathcal{P},p,\ln^{\alpha}(\frac{\xi^{2}}{p})\right) +k\sum_{z_{1}\leq p<z}\ln^{\alpha}(p)S_{0}(\mathcal{A}_{p^{2}};\mathcal{P},p) \]

\ This is the equation (2.8)
\end{proof}

\subsection{Lemma 1.2}
\ For $S_{\alpha}(\mathcal{A};\mathcal{P},z,\xi)$ we have following propertys \
\[
\begin{array}{ll}
\ (1)\ \ S_{k,\alpha}(\mathcal{A};\mathcal{P},z,\ln^{\alpha}(\xi^{2}))\geq 0, & if \ \ \ln^{\alpha}(\xi^{2})\geq 0; \\
\ (2)\ \ S_{k,\alpha}(\mathcal{A};\mathcal{P},z,\ln^{\alpha}(\xi^{2}))\leq S_{k,\alpha}(\mathcal{A};\mathcal{P},z_{1},\ln^{\alpha}(\xi^{2}))), & if \ \ 2\leq z_{1}<z,\ \ and \ \  k\ln^{\alpha}(z)\leq \ln^{\alpha}(\xi^{2}); \\
\ (3)\ \ S_{k,\alpha}(\mathcal{A};\mathcal{P},z,\ln^{\alpha}(\xi_{1}^{2}))\leq S_{k,\alpha}(\mathcal{A};\mathcal{P},z,\ln^{\alpha}(\xi_{2}^{2})), & if \ \ \ln^{\alpha}(\xi_{1}^{2})\leq \ln^{\alpha}(\xi_{2}^{2}) \\

\end{array}
\]
\begin{proof}
\ These equations are easy to see from Lemma 1.1.
\end{proof}

\subsection{Lemma 1.3}

\ Continuum function $h(u)$ define as

$$
 \begin{cases}
\begin{array}{ll}
h(u)=u & 0\leq u \leq 1 \\
(\frac{h(u)}{u})'=-\frac{h(u-1)}{u^{2}} & u>1
\end{array}
\end{cases}
$$

\ We have: \\
\begin{enumerate}
	\item $h(u)$ is an increase function, with the condition $h(\infty)=e^{\gamma}$
	\item $h(u)=h(\infty)+O(e^{-u\ln(u)})$
	\item $h(u)=2u-u\ln(u)-1$,\  \ $1<u\leq 2$
\end{enumerate}
 
\begin{proof}
\ This paper only considers item(3) which can be seen from the function definition.
\end{proof}
\subsection{Lemma 1.4}
\ Selberg in his upper bound sieve used a function\cite{A. Selberg1950} $G_{h}(x,z)$  defined as: 
\[ G_{h}(x,z)=\sum_{l|\mathcal{P}(z) \atop l<x,(l,h)=1}g(l) \]
$$
 g(l)=\begin{cases}
	\begin{array}{ll}
		1 & (l=1) \\
		\frac{w(p)}{p}\left(1-\frac{w(p)}{p}\right)^{-1} & (p,\bar{\mathcal{P}})=1 \\
		\prod_{p|l}g(p)=\frac{w(l)}{l}\prod_{p|l}\left(1-\frac{w(p)}{p}\right)^{-1}, & \mu(l)\neq 0, (p,\bar{\mathcal{P}})=1
	\end{array}
	\end{cases}
$$

\ We have \\

\begin{equation}
 G_{1}(\xi,z)=\frac{1}{C(\omega)}\ln(z) h(\frac{\ln(\xi)}{\ln(z)})(1+O(\frac{1}{\ln(z)}))
\label{GFunction}
\end{equation}

\ Where 

\[ C(\omega)=\prod_{p}(1-\frac{\omega(p)}{p})(1-\frac{1}{p})^{-1} \]

\ and

\begin{equation}
 \sum_{l|\mathcal{P}(\xi) \atop l<\xi}g(l)\ln^{a}(l)=\frac{1}{C(\omega)}\frac{1}{a+1}\ln^{a+1}(\xi)+O(\ln^{a}(\xi))
\end{equation}

\begin{proof}
\ When $a=0$ equation (2.10) is right, when $a\geq 1$
\[ \sum_{l|\mathcal{P}(\xi) \atop l<\xi}g(l)\ln^{a}(l)=\int_{1}^{\xi}\ln^{a}(t)d\left(\sum_{l|\mathcal{P}(\xi) \atop l<t}g(l)\right) \]
\[ =\frac{1}{C(\omega)}\int_{1}^{\xi}\ln^{a}(t)d\ln(t)(1+O(\frac{1}{\ln(t)}))\]
\[ =\frac{1}{C(\omega)}\frac{1}{a+1}\ln^{a+1}(\xi)+O\left(\int_{1}^{\xi}\ln^{a-1}(t)d\ln(t)\right) \]
\[ =\frac{1}{C(\omega)}\frac{1}{a+1}\ln^{a+1}(\xi)+O\left(\ln^{a}(\xi)\right) \]

\ This is the equation (2.10).

\ We know when $\xi\leq z$ $\frac{\ln(\xi)}{\ln(z)}\leq 1$ equation (2.9) is valid\cite{A. Selberg1950}.
\ When $\xi> z$, suppose $z_{0}>z$,$\frac{\ln(\xi)}{\ln(z_{0})}\leq \frac{\ln(\xi)}{\ln(z)}-1$ equation  (2.9) is valid. By definition of the function:
\[ G_{1}(\xi,z_{0})=\sum_{l|\mathcal{P}(z_{0}) \atop l<x}g(l) \]
\[ =\sum_{l|\mathcal{P}(z) \atop l<x}g(l)+\sum_{z\leq p< z_{0}}g(p)\sum_{l|\mathcal{P}(p) \atop l<x/p,(l,p)=1}g(l) \]
\[ =G_{1}(\xi,z)+\sum_{z\leq p< z_{0}}g(p)G_{p}(\xi/p,p) \]
\ Thus
\[ G_{1}(\xi,z)=G_{1}(\xi,z_{0})-\sum_{z\leq p< z_{0}}g(p)G_{p}(\xi/p,p) \]
\ Since $(p,\mathcal{P}(p))=1$, $G_{p}(\xi/p,p)=G_{1}(\xi/p,p)$. So
\[ G_{1}(\xi,z)=\frac{1}{C(\omega)}\ln(z_{0}) h(\frac{\ln(\xi)}{\ln(z_{0})})(1+O(\frac{1}{\ln(z_{0})}))-\sum_{z\leq p< z_{0}}g(p)\frac{1}{C(\omega)}\ln(p)h(\frac{\ln(\xi)}{\ln(p)}-1)(1+O(\frac{1}{\ln(p)})) \]
\[ =\frac{1}{C(\omega)}\ln(z_{0}) h(\frac{\ln(\xi)}{\ln(z_{0})})(1+O(\frac{1}{\ln(z_{0})})) \]
\[ -\frac{1}{C(\omega)}\int_{\frac{\ln(\xi)}{\ln(z_{0})}}^{\frac{\ln(\xi)}{\ln(z)}}\frac{1}{t^{2}}h(t-1)(1+O(\frac{1}{\ln(z)}))dt \]
\[ =\frac{1}{C(\omega)}\ln(z) h(\frac{\ln(\xi)}{\ln(z)})(1+O(\frac{1}{\ln(z)})) \]
\ This is the equation (2.9)
\end{proof}

\subsection{Lemma 1.5 } 
\ This lemma is very similar to the Selberg upper bound Sieve\cite{A. Selberg1950}, only the function $\lambda_{d}$ has a minor differance to that used by Selberg. \\
\ Suppose $\sqrt{\xi}\leq z \leq \xi$

\begin{equation} \lambda_{d}=\mu(d)\prod_{p|d}(1-\frac{w(p)}{p})^{-1}G_{d}(\frac{\xi}{d},\xi)(G_{1}(\xi,\xi))^{-1},\  \ d|P(z), d< \xi \\
\label{lambda1}
\end{equation}

\ We have

\begin{equation}
 \sum_{n \in \mathcal{A}}\left(\sum_{d|(n,P(z))}\lambda_{d}\right)^{2} =X\sum_{l|\mathcal{P}(z)}\frac{y^{2}_{l}}{g(l)}+\sum_{d_{i}|P(z),d_{i}\leq \xi \atop i=1,2}\lambda_{d_{1}}\lambda_{d_{2}}r_{[d_{1},d_{2}]}
\end{equation}

\[ \leq X\frac{C(\omega)}{\ln(z)}\left(\frac{h(u)}{u^{2}}+(2+\ln(u))\left(\frac{1}{u}(1-\frac{1}{u})-\frac{1}{2u}(1-\frac{1}{u})^{2}-\frac{1}{6u}(1-\frac{1}{u})^{3}\right)\right)+O(\frac{1}{\ln^{2}(z)})\]+R

\ Where $u=\frac{\ln(\xi)}{\ln(z)}$, $1\leq u\leq 2$ and\\

\begin{equation}
 y_{l}=\sum_{d|\mathcal{P}(z) \atop l|d}\lambda_{d}\frac{\omega(d)}{d}, \ \ l|\mathcal{P}(z)
\end{equation}

\[ R=\sum_{d_{i}|P(z),d_{i}\leq \xi \atop i=1,2}\lambda_{d_{1}}\lambda_{d_{2}}r_{[d_{1},d_{2}]} \]
\begin{proof}
\ The front half of equation (2.12) is the same as Selberg's upper bound sieve\cite{A. Selberg1950}. According to equation (2.11)
\[ y_{l}=\sum_{d|\mathcal{P}(z) \atop l|d}\frac{\omega(d)}{d}\mu(d)\prod_{p|d}(1-\frac{w(p)}{p})^{-1}G_{d}(\frac{\xi}{d},\xi)(G_{1}(\xi,\xi))^{-1} \]
\[ =\sum_{d|\mathcal{P}(z) \atop l|d}\mu(d)g(d)G_{d}(\frac{\xi}{d},\xi)(G_{1}(\xi,\xi))^{-1} \]
\[ =\mu(l)g(l)(G_{1}(\xi,\xi))^{-1}\sum_{d|\mathcal{P}(z) \atop (d,l)=1}\mu(d)g(d)\sum_{k<\frac{\xi}{dl},k|\mathcal{P}(\xi),(k,dl)=1}g(k) \]
\[ =\mu(l)g(l)(G_{1}(\xi,\xi))^{-1}\sum_{k|\mathcal{P}(\xi),(k,l)=1,k<\xi/l}g(k)\sum_{d|(\mathcal{P}(z)/l,k)}\mu(d) \]
\[ =\mu(l)g(l)(G_{1}(\xi,\xi))^{-1}\sum_{k|\mathcal{P}(\xi),(k,\mathcal{P}(z))=1 \atop(k,l)=1,k<\xi/l}g(k) \]
\ Since $z\geq \sqrt{\xi}$  \\
\[ \sum_{k|\mathcal{P}(\xi),(k,\mathcal{P}(z))=1 \atop(k,l)=1,k<\xi/l}g(k)=(1+\sum_{z\leq p<\xi/l}g(p)) \]
\[ y_{l}=\mu(l)g(l)(G_{1}(\xi,\xi))^{-1}\left(1+\sum_{z\leq p<\xi/l}g(p)\right) \]

\ Thus \\

$$
 y_{l}=\begin{cases}
\begin{array}{ll}
 \mu(l)g(l)(G_{1}(\xi,\xi))^{-1}\left(1+\ln(\frac{\ln(\xi/l)}{\ln(z)})+0(\frac{1}{\ln(z)})\right) & l< \xi/z\\
 \mu(l)g(l)(G_{1}(\xi,\xi))^{-1} & l\geq \xi/z \\
\end{array}
\end{cases}
$$

\[ \sum_{l|\mathcal{P}(z)}\frac{y_{l}^{2}}{g(l)}=(G_{1}(\xi,\xi))^{-2}\sum_{l|\mathcal{P}(z),l<\xi}g(l)\]

\[ +(G_{1}(\xi,\xi))^{-2}\sum_{l|\mathcal{P}(z),l<\xi/z}g(l)\left(2\ln(\frac{\ln(\xi/l)}{\ln(z)})+\ln^{2}(\frac{\ln(\xi/l)}{\ln(z)})+O(\frac{1}{\ln(z)})\right) \]
\[ \leq \frac{G_{1}(\xi,z)}{G_{1}^{2}(\xi,\xi)}+(G_{1}(\xi,\xi))^{-2}\left(2+\ln(\frac{\ln(\xi)}{\ln(z)})\right)\sum_{l|\mathcal{P}(z),l<\xi/z}g(l)\left(\ln(\frac{\ln(\xi/l)}{\ln(z)})+O(\frac{1}{\ln(z)})\right) \]

\ Since

\[ \ln(\frac{\ln(\xi/l)}{\ln(z)})=\ln\left(\frac{\ln(\xi)}{\ln(z)}(1-\frac{\ln(l)}{\ln(\xi)})\right)=\ln(\frac{\ln(\xi)}{\ln(z)})+\ln((1-\frac{\ln(l)}{\ln(\xi)})) \]

\[ \leq \ln(\frac{\ln(\xi)}{\ln(z)})-\frac{\ln(l)}{\ln(\xi)}-\frac{1}{2}\frac{\ln^{2}(l)}{\ln^{2}(\xi)} \]

\ Hence, the sum on the right
\[ \sum_{l|\mathcal{P}(z),l<\xi/z}g(l)\ln(\frac{\ln(\xi/l)}{\ln(z)}) \]
\[ \leq \ln(\frac{\ln(\xi)}{\ln(z)})\sum_{l|\mathcal{P}(z),l<\xi/z}g(l) - \frac{1}{\ln(\xi)}\sum_{l|\mathcal{P}(z),l<\xi/z}g(l)\ln(l)-\frac{1}{2\ln^{2}(\xi)}\sum_{l|\mathcal{P}(z),l<\xi/z}g(l)\ln^{2}(l) \]
\[ =\ln(\frac{\ln(\xi)}{\ln(z)})G_{1}(\xi/z,z)-\frac{1}{2C(\omega)}\frac{\ln^{2}(\xi/z)}{\ln(\xi)}-\frac{1}{6C(\omega)}\frac{\ln^{3}(\xi/z)}{\ln^{2}(\xi)} +O(1)\]

\[ \sum_{l|\mathcal{P}(z)}\frac{y_{l}^{2}}{g(l)}\leq \frac{G_{1}(\xi,z)}{G_{1}^{2}(\xi,\xi)}\]
\[ +(G_{1}(\xi,\xi))^{-2}\left(2+\ln(\frac{\ln(\xi)}{\ln(z)})\right)\left(\ln(\frac{\ln(\xi)}{\ln(z)})G_{1}(\xi/z,z)-\frac{\ln^{2}(\xi/z)}{2C(\omega)\ln(\xi)}-\frac{\ln^{3}(\xi/z)}{6C(\omega)\ln^{2}(\xi)}\right) \]
\[ +O(\frac{1}{\ln^{2}(\xi)}) \]

\[ =\frac{1}{G_{1}(\xi,z)}\left(\frac{G_{1}^{2}(\xi,z)}{G_{1}^{2}(\xi,\xi)}+\frac{G_{1}(\xi,z)G_{1}(\xi/z,z)}{G_{1}^{2}(\xi,\xi)}(2+\ln(\frac{\ln(\xi)}{\ln(z)}))\right) \]
\[ -\frac{1}{G_{1}(\xi,z)}\left(\frac{G_{1}(\xi,z)\ln^{2}(\xi/z)}{2C(\omega)G_{1}^{2}(\xi,\xi)\ln(\xi)}+\frac{G_{1}(\xi,z)\ln^{3}(\xi/z)}{6C(\omega)G_{1}^{2}(\xi,\xi)\ln^{2}(\xi)}\right)(2+\ln(\frac{\ln(\xi)}{\ln(z)})) +O(\frac{1}{\ln^{2}(\xi)})\]

\[ =\frac{1}{G_{1}(\xi,z)}\left(\frac{\ln^{2}(z)h^{2}(\frac{\ln(\xi)}{\ln(z)})}{\ln^{2}(\xi)h^{2}(1)}+\frac{\ln(z)\ln(\xi/z)h(\frac{\ln(\xi)}{\ln(z)})}{\ln^{2}(\xi)h(1)}(2+\ln(\frac{\ln(\xi)}{\ln(z)}))\right)\]
\[ -\frac{1}{G_{1}(\xi,z)}\left(\frac{\ln(z)\ln^{2}(\xi/z)h(\frac{\ln(\xi)}{\ln(z)})}{2\ln^{3}(\xi)h(1)}+\frac{\ln(z)\ln^{3}(\xi/z)h(\frac{\ln(\xi)}{\ln(z)})}{6\ln^{4}(\xi)h(1)}\right) (2+\ln(\frac{\ln(\xi)}{\ln(z)}))\]
\[ +O(\frac{1}{\ln^{2}(\xi)}) \]

\ Setting $u=\frac{\ln(\xi)}{\ln(z)}$, we obtain

\[ \sum_{l|\mathcal{P}(z)}\frac{y_{l}^{2}}{g(l)}\leq \frac{C(\omega)}{\ln(z)h(u)}\left(\frac{h^{2}(u)}{u^{2}h^{2}(1)}+\frac{1}{u}(1-\frac{1}{u})\frac{h(u)}{h(1)}(2+\ln(u))\right)\]
\[ -\frac{C(\omega)}{\ln(z)h(u)}\left(\frac{1}{2u}(1-\frac{1}{u})^{2}\frac{h(u)}{h(1)}+\frac{1}{6u}(1-\frac{1}{u})^{3}\frac{h(u)}{h(1)}\right)(2+\ln(u))+O(\frac{1}{\ln^{2}(\xi)}) \]

\[ =\frac{C(\omega)}{\ln(z)}\left(\frac{h(u)}{u^{2}}+(2+\ln(u))\left(\frac{1}{u}(1-\frac{1}{u})-\frac{1}{2u}(1-\frac{1}{u})^{2}-\frac{1}{6u}(1-\frac{1}{u})^{3}\right)\right)+O(\frac{1}{\ln^{2}(\xi)})\]

\end{proof}

\subsection{Lemma 1.6}

\ Suppose  $1\leq u=\frac{\ln(\xi)}{\ln(z)}\leq 2$, $2\leq p<z$ \\

\[ 1+\lambda_{p}\geq 
 \frac{\ln(p)}{\ln(\xi)}(1+O(\frac{1}{\ln(z)})) \]

\begin{proof}
\ Since 
\[ G_{p}(\frac{\xi}{p},\xi)=G_{1}(\frac{\xi}{p},\xi)-g(p)G_{p}(\frac{\xi}{p^{2}},\xi)\leq G_{1}(\frac{\xi}{p},\xi) \]
\ According to lemma (1.5) 
\[ 1+\lambda_{p}=1+\mu(p)(1-\frac{w(p)}{p})^{-1}G_{p}(\frac{\xi}{p},\xi)G_{1}^{-1}(\xi,\xi) \]
\[ \geq 1-\frac{\ln(\frac{\xi}{p})}{\ln(\xi)}(1+O(\frac{1}{\ln(\xi)}))=\frac{\ln(p)}{\ln(\xi)}(1+O(\frac{1}{\ln(\xi)})) \]
\end{proof}

\subsection{Lemma 1.7}

\ When $2 \leq u=\frac{\ln(\xi^{2})}{\ln(z)}\leq k_{n}\leq 2^{\alpha}$, $\alpha \geq 2$ we have

\begin{equation}
\ S_{k_{n},\alpha}(\mathcal{A};\mathcal{P},z,\ln^{\alpha}(\xi^{2}))=\ln^{\alpha}(\xi^{2})S(\mathcal{A};\mathcal{P},z)+k_{n}\sum_{2\leq p<z}\ln^{\alpha}(p)S(\mathcal{A}_{p};\mathcal{P}(p),z) 
\end{equation} 

\[ \leq C(\omega)\frac{Xe^{-\gamma}\ln^{\alpha}(\xi^{2})\tilde{F}(u)}{\ln(z)}\left(1+O(\frac{1}{\ln(\xi)})\right)+\ln^{\alpha}(\xi^{2})\sum_{d|P(z),d\leq \xi^{2}}3^{v_{1}(d)}\left|r_{d}\right| \]

\[ =C(\omega)Xe^{-\gamma}\ln^{\alpha-1}(\xi^{2})u\tilde{F}(u)\left(1+O(\frac{1}{\ln(\xi)})\right)+\ln^{\alpha}(\xi^{2})\sum_{d|P(z),d\leq \xi^{\alpha}}3^{v_{1}(d)}\left|r_{d}\right| \]

\ Where

\[  \tilde{F}(u)=e^{\gamma}\left(\frac{4h(u/2)}{u^{2}}+(2+\ln(u/2))\frac{2}{u}\left((1-\frac{2}{u})-\frac{1}{2}(1-\frac{2}{u})^{2}-\frac{1}{6}(1-\frac{2}{u})^{3}\right) \right)\]
\begin{proof}
\ According to Lemma 1.6, when $\alpha=2$
\[ 1+\lambda_{p}\geq \frac{\ln(p)}{\ln(\xi)}(1+O(\frac{1}{\ln(\xi)}))=2\frac{\ln(p)}{\ln(\xi^{2})}(1+O(\frac{1}{\ln(\xi)})) \]
\[ (1+\lambda_{p})^{2}\geq 4\frac{\ln^{2}(p)}{\ln^{2}(\xi^{2})}(1+O(\frac{1}{\ln(\xi)})) \]
\ we obtain
\[ S_{4,2}(\mathcal{A};\mathcal{P},z,\ln^{2}(\xi^{2}))\leq \ln^{2}(\xi^{2})\sum_{n \in \mathcal{A}}\left(\sum_{d|(n,\mathcal{P}(z))}\lambda_{d}\right)^{2}(1+O(\frac{1}{\ln(\xi)})) \]
\ According to equation (2.12) get when $\alpha=2$ $Lemma1.7$ is valid. When $\alpha>2$, since
\[ \ln(p)<\ln(z)\leq \frac{1}{2}\ln(\xi^{2}), \ \  2\frac{\ln(p)}{\ln(\xi^{2})}\leq 1 \]
\[ k_{n}\frac{\ln^{\alpha}(p)}{\ln^{\alpha}(\xi^{2})}=2^{\alpha-2}\frac{\ln^{\alpha-2}(p)}{\ln^{\alpha-2}(\xi^{2})}4\frac{\ln^{2}(p)}{\ln^{2}(\xi^{2})}\leq 4\frac{\ln^{2}(p)}{\ln^{2}(\xi^{2})} \]
\[ S_{k_{n},\alpha}(\mathcal{A};\mathcal{P},z,\ln^{\alpha}(\xi^{2}))=\ln^{\alpha}(\xi^{2})S_{0}(\mathcal{A};\mathcal{P},z)+\sum_{2\leq p<z}k_{n}\ln^{\alpha}(p)S_{0}(\mathcal{A}_{p};\mathcal{P}(p),z) \]
\[ \leq \ln^{\alpha-2}(\xi^{2})S_{4,2}(\mathcal{A};\mathcal{P},z,\ln^{2}(\xi^{2})) \]
\ So that $Lemma1.7$ is valid.
\end{proof}

\subsection{Lemma 1.8} Fundamental Lemma of Selberg's sieve\cite{A. Selberg1950}\cite{Halberstam1974}\cite{C.D. Pan and C.B. Pan1992}.\\

\ Suppose
\[  W(z)=\prod_{p<z,(p,N)=1}(1-\frac{w(p)}{p})=\frac{C(\omega)e^{-\gamma}}{\ln(z)}\left(1+O(\frac{1}{\ln(z)})\right) \]

\ Where $C(\omega)$ is defined on $lemma 1.4$

\ For $2\leq z\leq \xi$, we have

\begin{equation} 
\ S_{0}(\mathcal{A};\mathcal{P}(N),z)=XW(z)\left\{1+O\left(exp(-\frac{1}{6}\tau \ln \tau)\right)\right\}
\label{S0}
\end{equation}

\[ +\theta \sum_{d|P(z),d< \xi^{2}}3^{v_{1}(d)}\left|r_{d}\right|, \  \ \left|\theta \right|\leq 1, \]

\ Here $\tau=\frac{\ln(\xi^{2})}{\ln(z)}$, constant of $"O"$ is independent of $\tau$.

\ Similarly, we have
\subsection{Lemma 1.9} \

\ Set $\tau=\frac{\ln(\xi^{2})}{\ln(z)}$, $\alpha>0$ we have

\begin{equation}
\ S_{k,\alpha}(\mathcal{A};\mathcal{P}(N),z,\ln(\xi^{2}))=XW(z)\left(\ln^{\alpha}(\xi^{2})+\frac{k}{\alpha} \ln^{\alpha}(z)\right)\left(1+O(e^{-\frac{1}{6}\tau \ln(\tau)})+O(\frac{1}{\ln(z)})\right)
\label{S1}
\end{equation}

\[ +\theta \ln^{\alpha}(\xi^{2})\sum_{d|P(z),d< \xi^{2}}3^{v_{1}(d)}\left|r_{d}\right|, \  \ \left|\theta \right|\leq 1, \]

\begin{proof}By equation (2.6)
\[ S_{k,\alpha}(\mathcal{A};\mathcal{P},z,\ln(\xi^{2}))=\ln^{\alpha}(\xi^{2})S(\mathcal{A};P,z)+\sum_{p<z}\ln^{\alpha}(p)S_{0}(\mathcal{A}_{p};\mathcal{P}(p),z) \]
\ According to equation (2.15) on $Lemma 1.8$

\[ S_{k,\alpha}(\mathcal{A};\mathcal{P},z,\ln^{\alpha}(\xi^{2}))=XW(z)\ln^{\alpha}(\xi^{2})\left\{1+O\left(exp(-\frac{1}{6}\tau \ln \tau)\right)\right\} \]
\[ +kXW(z)\sum_{p<z}\frac{w(p)\ln^{\alpha}(p)}{p}\left\{1+O\left(exp(-\frac{1}{6}\tau_{p} \ln \tau_{p})\right)\right\} \]
\[ +\theta \ln^{\alpha}(\xi^{2})\sum_{d|P(z),d< \xi^{2}}3^{v_{1}(d)}\left|r_{d}\right|
 +\theta \sum_{p<Z}\sum_{d|\frac{P(z)}{p},d< \frac{\xi^{2}}{p}}3^{v_{1}(d)}\ln^{\alpha}(p)\left|r_{pd}\right|, \  \ \left|\theta \right|\leq 1, \]

\ Where $\tau_{p}=\frac{\ln(\xi^{2}/p)}{\ln(z)}$. Since $p<z$, so 
\[ exp(-\frac{1}{6}\tau_{p} \ln \tau_{p})\leq O(exp(-\frac{1}{6}\tau \ln \tau)) \]
\ and

\[ \sum_{p<Z}\sum_{d|\frac{\mathcal{P}(z)}{p},d< \frac{\xi^{2}}{p}}3^{v_{1}(d)}\ln^{\alpha}(p)\left|r_{pd}\right| \]
\[ \leq \sum_{d|P(z),d< \xi^{\alpha+1}}3^{v_{1}(d)}\left|r_{d}\right|\sum_{p|d}\ln(p) \leq \sum_{d|P(z),d< \xi^{2}}3^{v_{1}(d)}\ln(d)\left|r_{d}\right| \]

\ So we obtain
\[ S_{k,\alpha}(\mathcal{A};\mathcal{P},z,\ln^{\alpha}(\xi^{2}))=XW(z)\left(\ln^{\alpha}(\xi^{2})+k\sum_{p<z}\frac{w(p)\ln^{\alpha}(p)}{p})\right)\left\{1+O\left(exp(-\frac{1}{6}\tau \ln \tau)\right)\right\} \]
\[ +\theta \ln^{\alpha}(\xi^{2})\sum_{d|P(z),d< \xi^{2}}3^{v_{1}(d)}\left|r_{d}\right|, \ \ \left|\theta \right|\leq 1, \]

\[ =XW(z)\left(\ln^{\alpha}(\xi^{2})+\frac{k}{\alpha}\ln^{\alpha}(z)+O(\ln^{\alpha-1}(z))\right)\left\{1+O\left(exp(-\frac{1}{6}\tau \ln \tau)\right)\right\} \]
\[ +\theta \ln^{\alpha}(\xi^{2})\sum_{d|P(z),d< \xi^{2}}3^{v_{1}(d)}\left|r_{d}\right|, \ \ \left|\theta \right|\leq 1, \]

\end{proof}

\subsection{Lemma 1.10} (The theorem of Jurkat-Richert\cite{W.B. Jurkat and H.-E. Richert1965}):

\ Suppose $2\leq z \leq \xi$, we have

\begin{equation}
\ S_{0}(\mathcal{A};\mathcal{P},z)\leq 
XW(z)\left\{F(\frac{\ln(\xi^{2})}{\ln(z)})+O\left(\frac{1}{\ln^{\frac{1}{14}}(\xi)}\right)\right\}\end{equation}

\[ +\sum_{d|P(z),d<\xi^{2}}3^{\nu_{1}(d)}\left|r_{d}\right| \] 

\ And

\begin{equation}
\ S_{0}(\mathcal{A};\mathcal{P},z)\geq XW(z)\left\{f(\frac{\ln(\xi^{2})}{\ln(z)})+O\left(\frac{1}{\ln^{\frac{1}{14}}(\xi)}\right)\right\}\end{equation}

\[ -\sum_{d|P(z),d<\xi^{2}}3^{\nu_{1}(d)}\left|r_{d}\right| \]
\ Where continuum functions $F(u)$, and $f(u)$ are defined as:

$$
 \begin{cases}
\begin{array}{lll}
	F(u)=\frac{2e^{\gamma}}{u}, & f(u)=0, & 1\leq u\leq 2 \\
	(uF(u))'=f(u-1), & (uf(u))'=F(u-1), & u>2
\end{array}
\end{cases}
$$

\subsection{Corollary 1.10.1}
\ When $\frac{\ln(\xi^{2})}{\ln(z)}=u \geq 2$, $\alpha\geq 1$ 

\[ S_{k,\alpha}(\mathcal{A};\mathcal{P},z,\ln^{\alpha}(\xi^{2}))\leq XW(z)\left\{\ln^{\alpha}(\xi^{2})F(\frac{\ln(\xi^{2})}{\ln(z)})+\frac{k\ln^{\alpha}(z)}{\alpha}F(\frac{\ln(\xi^{2})}{\ln(z)}-1)+O(\frac{1}{\ln^{\frac{1}{14}}(\xi)})\right\} \]

\[ +\sum_{d|\mathcal{P}(z),d<\xi^{2}}3^{v_{1}(d)}|r_{d}| \]

\[ =\ln^{\alpha-1}(\xi^{2})Xe^{-\gamma}C(\omega)\left\{uF(u)+\frac{k}{\alpha u^{\alpha-1}}F(u-1)+O(\frac{1}{\ln^{\frac{1}{14}}(\xi)})\right\} \]

\[ +\ln^{\alpha}(\xi^{2})\sum_{d|\mathcal{P}(z),d<\xi^{2}}3^{v_{1}(d)}|r_{d}| \]

\ And
\[ S_{k,\alpha}(\mathcal{A};\mathcal{P},z,\ln^{\alpha}(\xi^{2})) \geq \ln^{\alpha-1}(\xi^{2})Xe^{-\gamma}C(\omega)\left\{uf(u)+\frac{k}{\alpha u^{\alpha-1}}f(u-1)+O(\frac{1}{\ln^{\frac{1}{14}}(\xi)})\right\} \]

\[ -\ln^{\alpha}(\xi^{2})\sum_{d|\mathcal{P}(z),d<\xi^{2}}3^{v_{1}(d)}|r_{d}| \]

\begin{proof}By equation (2.6), $Lemma1.10$
\[ S_{k,\alpha}(\mathcal{A};\mathcal{P},z,\ln^{\alpha}(\xi^{2}))=\ln^{\alpha}(\xi^{2})S_{0}(\mathcal{A};\mathcal{P},z) \]
\[ +\sum_{p<z}k\ln^{\alpha}(p)S_{0}(\mathcal{A}_{p};\mathcal{P}(p),z) \]
\[ \leq Xe^{-\gamma}C(\omega)\ln^{\alpha-1}(\xi^{2})uF(u)(1+\frac{1}{\ln^{\frac{1}{14}}(z)})+\ln^{\alpha}(\xi^{2})R \]
\[ +\sum_{p<z}\frac{\ln^{\alpha}(p)}{\ln(\xi^{2})}\frac{X}{p}e^{-\gamma}C(\omega)uF(u_{p})(1+\frac{1}{\ln^{\frac{1}{14}}(z)})+\sum_{p<z}\ln^{\alpha}(p)R_{p} \]
\ Where $u_{p}=\frac{\ln(\xi^{2}/p)}{\ln(z)}\geq \frac{\ln(\xi^{2})}{\ln(z)}-1=u-1$. \\
\[ R=\sum_{d|\mathcal{P}(z),d<\xi^{2}}3^{v_{1}(d)}|r_{d}| \]
\[ R_{p}=\sum_{d|\mathcal{P}(z),d<\xi^{2}/p,(d,p)=1}3^{v_{1}(d)}|r_{pd}| \]
\[ \sum_{p<z}\ln^{\alpha}(p)R_{p}\leq \ln^{\alpha}(\xi^{2})R \]
\ So that when $\alpha \geq 1$
\[ \sum_{p<z}\frac{\ln^{\alpha}(p)}{\ln(\xi^{2})}\frac{X}{p}e^{-\gamma}C(\omega)uF(u_{p})\leq Xe^{-\gamma}C(\omega)uF(u-1)\ln^{-1}(\xi^{2})\sum_{p<z}\frac{\ln^{\alpha}(p)}{p} \]
\[ =Xe^{-\gamma}C(\omega)uF(u-1)\ln^{-1}(\xi^{2})\frac{\ln^{\alpha}(z)}{\alpha}(1+O(\frac{1}{\ln(z)})) =\ln^{\alpha-1}(\xi^{2})Xe^{-\gamma}C(\omega)uF(u-1)\frac{1}{\alpha}\frac{\ln^{\alpha}(z)}{\ln^{\alpha}(\xi^{2})}(1+O(\frac{1}{\ln(z)})) \]
\[ =\ln^{\alpha-1}(\xi^{2})Xe^{-\gamma}C(\omega)\frac{F(u-1)}{\alpha u^{\alpha-1}}(1+O(\frac{1}{\ln(z)})) \]
\ Combining these to get the quation of upper bound. Proof the equation of lower bound is same as this.

\end{proof}
\subsection{Lemma 1.11 }

\ Suppose , $2\leq w<z\leq \xi^{2}$, and $k\geq 0$, we have

\begin{equation}
\	\ln^{\alpha}(\xi^{2})S_{0}(\mathcal{A};\mathcal{P},z)+k\sum_{2\leq p<w}\ln^{\alpha}(p)S_{0}(\mathcal{A}_{p};\mathcal{P}(p),z)+\frac{1}{2}\ln^{\alpha}(\xi^{2})\sum_{w\leq p<z}S_{0}(\mathcal{A}_{p};\mathcal{P}(p),z)
\end{equation}

\[ =S_{k,\alpha}(\mathcal{A};\mathcal{P},w,\ln^{\alpha}(\xi^{2}))\]
\[ -\sum_{w\leq p<z}\left(\frac{1}{2}\ln^{\alpha}(\xi^{2})S_{0}(\mathcal{A}_{p};\mathcal{P},p)+k\sum_{q<w}\ln^{\alpha}(q)S_{0}(\mathcal{A}_{pq};\mathcal{P}(q),p)+\frac{1}{2}\sum_{w\leq q<p}\ln^{\alpha}(\xi^{2})S_{0}(\mathcal{A}_{pq};\mathcal{P}(q),p)\right) \]
\[ -\frac{1}{2}\sum_{w\leq p<z}\ln^{\alpha}(\xi^{2})S_{0}(\mathcal{A}_{p^{2}};\mathcal{P},p)\]

\[ \geq S_{k,\alpha}(\mathcal{A};\mathcal{P},w,\ln^{\alpha}(\xi^{2}))-\sum_{w\leq p<z}S_{k,\alpha}(\mathcal{A}_{p};\mathcal{P},w,\frac{1}{2}\ln^{\alpha}(\xi^{2}))-\frac{1}{2}\sum_{w\leq p<z}\ln^{\alpha}(\xi^{2})S_{0}(\mathcal{A}_{p^{2}};\mathcal{P},p) \]

\begin{proof} similar of proof equation (2.8) on $Lemma 1.1$\
\[ S_{0}(\mathcal{A};\mathcal{P},z)=S_{0}(\mathcal{A};\mathcal{P},2)-\sum_{2\leq p<z}S_{0}(\mathcal{A}_{p};\mathcal{P},p) \]
\ The sum on the right
\[ \ln^{\alpha}(\xi^{2})\sum_{2\leq p<z}S_{0}(\mathcal{A}_{p};\mathcal{P},p)=\sum_{2\leq p<w}\left(\ln^{\alpha}(\xi^{2})-k\ln^{\alpha}(p)\right)S_{0}(\mathcal{A}_{p};\mathcal{P}(p),p)\]
\[ +\sum_{2\leq p<w}k\ln^{\alpha}(p)S_{0}(\mathcal{A}_{p};\mathcal{P}(p),p)+\frac{1}{2}\sum_{w\leq p<z}\ln^{\alpha}(\xi^{2})S_{0}(\mathcal{A}_{p};\mathcal{P}(p),p) +\frac{1}{2}\sum_{w\leq p<z}\ln^{\alpha}(\xi^{2})S_{0}(\mathcal{A}_{p};\mathcal{P}(p),p)\]

\[ =\sum_{2\leq p<w}\left(\ln^{\alpha}(\xi^{2})-k\ln^{\alpha}(p)\right)S_{0}(\mathcal{A}_{p};\mathcal{P}(p),p)+\frac{1}{2}\sum_{w\leq p<z}\ln^{\alpha}(\xi^{2})S_{0}(\mathcal{A}_{p};\mathcal{P}(p),p)\]
\[ +\sum_{2\leq p<w}k\ln^{\alpha}(p)S_{0}(\mathcal{A}_{p};\mathcal{P}(p),z)+\sum_{2\leq p<w,p\leq  q<z}k\ln^{\alpha}(p)S_{0}(\mathcal{A}_{pq};\mathcal{P}(p),q) \]
\[ +\sum_{w\leq p<z}\frac{1}{2}\ln^{\alpha}(\xi^{2})S_{0}(\mathcal{A}_{p};\mathcal{P}(p),z)+\sum_{w\leq p<z,p\leq q<z}\frac{1}{2}\ln^{\alpha}(\xi^{2})S_{0}(\mathcal{A}_{pq};\mathcal{P}(p),q) \]
\ So that
\[ \ln^{\alpha}(\xi^{2})S_{0}(\mathcal{A};\mathcal{P},z)+\sum_{2\leq p<w}k\ln^{\alpha}(p)S_{0}(\mathcal{A}_{p};\mathcal{P}(p),z)+\sum_{w\leq p<z}\frac{1}{2}\ln^{\alpha}(\xi^{2})S_{0}(\mathcal{A}_{p};\mathcal{P}(p),z) \]
\[ =\ln^{\alpha}(\xi^{2})S_{0}(\mathcal{A};\mathcal{P},2)-\sum_{2\leq p<w}\left(\left(\ln^{\alpha}(\xi^{2})-k\ln^{\alpha}(p)\right)S_{0}(\mathcal{A}_{p};\mathcal{P}(p),p)+\sum_{2\leq q\leq p}k\ln^{\alpha}(q)S_{0}(\mathcal{A}_{pq};\mathcal{P}(q),p)\right) \]
\[ -\sum_{w\leq p<z}\left(\frac{1}{2}\ln^{\alpha}(\xi^{2})S_{0}(\mathcal{A};\mathcal{P},p)+\sum_{q<w}k\ln^{\alpha}(q)S_{0}(\mathcal{A}_{pq};\mathcal{P}(q),p)+\sum_{w\leq q \leq p}\frac{1}{2}\ln^{\alpha}(\xi^{2})S_{0}(\mathcal{A}_{pq};\mathcal{P}(q),p)\right) \]
\ On the other hand
\[ S_{k,\alpha}(\mathcal{A};\mathcal{P},w,\ln^{\alpha}(\xi^{2}))=\ln^{\alpha}(\xi^{2})S_{0}(\mathcal{A};\mathcal{P},w)+\sum_{2\leq p<w}k\ln^{\alpha}(p)S_{0}(\mathcal{A}_{p};\mathcal{P}(p),w) \]
\[ =\ln^{\alpha}(\xi^{2})S_{0}(\mathcal{A};\mathcal{P},2)-\sum_{2\leq p<w}\left(\left(\ln^{\alpha}(\xi^{2})-k\ln^{\alpha}(p)\right)S_{0}(\mathcal{A}_{p};\mathcal{P},p)+\sum_{2\leq q\leq p}k\ln^{\alpha}(p)S_{0}(\mathcal{A}_{pq};\mathcal{P}(q),p)\right) \]
\ We obtain
\[ \ln^{\alpha}(\xi^{2})S_{0}(\mathcal{A};\mathcal{P},z)+\sum_{2\leq p<w}k\ln^{\alpha}(p)S_{0}(\mathcal{A}_{p};\mathcal{P}(p),z)+\sum_{w\leq p<z}\frac{1}{2}\ln^{\alpha}(\xi^{2})S_{0}(\mathcal{A}_{p};\mathcal{P}(p),z)\]
 \[ =\ln^{\alpha}(\xi)S_{k,\alpha}(\mathcal{A};\mathcal{P},w,\ln^{\alpha}(\xi^{2}))\]
\[ -\sum_{w\leq p<z}\left(\frac{1}{2}\ln^{\alpha}(\xi^{2})S_{0}(\mathcal{A}_{p};\mathcal{P},p)+\sum_{q<w}k\ln^{\alpha}(q)S_{0}(\mathcal{A}_{pq};\mathcal{P}(q),p)+\sum_{w\leq q \leq p}\frac{1}{2}\ln^{\alpha}(\xi^{2})S_{0}(\mathcal{A}_{pq};\mathcal{P}(q),p)\right) \]
\[ \sum_{w\leq p<z} \left(\frac{1}{2}\ln^{\alpha}(\xi^{2})S_{0}(\mathcal{A}_{p};\mathcal{P},p)+\sum_{w\leq q\leq p}\frac{1}{2}\ln^{\alpha}(\xi^{2})S_{0}(\mathcal{A}_{pq};\mathcal{P}(q),p)\right) \]
\[ =\frac{1}{2}\sum_{w\leq p<z}\ln^{\alpha}(\xi^{2})S_{0}(\mathcal{A}_{p};\mathcal{P},w) \]
\[ +\frac{1}{2}\ln^{\alpha}(\xi^{2})\sum_{w\leq p<z}\sum_{w\leq q<p}\left(S_{0}(\mathcal{A}_{pq};\mathcal{P}(q),p)-S_{0}(\mathcal{A}_{pq};\mathcal{P}(q),q)\right) \]
\[ +\frac{1}{2}\sum_{w\leq p<z}\ln^{\alpha}(\xi^{2})S_{0}(\mathcal{A}_{p^{2}};\mathcal{P},p) \]
\[ \leq \frac{1}{2}\sum_{w\leq p<z}\ln^{\alpha}(\xi^{2})S_{0}(\mathcal{A}_{p};\mathcal{P},w)+\frac{1}{2}\sum_{w\leq p<z}\ln^{\alpha}(\xi^{2})S_{0}(\mathcal{A}_{p^{2}};\mathcal{P},p) \]
\ and
\[ \sum_{w\leq p<z}\sum_{q<w}k\ln^{\alpha}(q)S_{0}(\mathcal{A}_{pq};\mathcal{P}(q),p)\leq \sum_{w\leq p<z}\sum_{q<w}k\ln^{\alpha}(q)S_{0}(\mathcal{A}_{pq};\mathcal{P}(q),w) \]

\ Hence

\[ \ln^{\alpha}(\xi^{2})S_{0}(\mathcal{A};\mathcal{P},z)+\sum_{2\leq p<w}k\ln^{\alpha}(p)S_{0}(\mathcal{A}_{p};\mathcal{P}(p),z)+\sum_{w\leq p<z}\frac{1}{2}\ln^{\alpha}(\xi^{2})S_{0}(\mathcal{A}_{p};\mathcal{P}(p),z)\]
 \[ \geq  S_{k,\alpha}(\mathcal{A};\mathcal{P},w,\ln^{\alpha}(\xi^{2}))\]
 
 \[ -\frac{1}{2}\sum_{w\leq p<z}\ln^{\alpha}(\xi^{2})S_{0}(\mathcal{A}_{p};\mathcal{P},w)-\sum_{w\leq p<z}\sum_{q<w}k\ln^{\alpha}(q)S_{0}(\mathcal{A}_{pq};\mathcal{P}(q),w) \]
\[ -\frac{1}{2}\sum_{w\leq p<z}\ln^{\alpha}(\xi^{2})S_{0}(\mathcal{A}_{p^{2}};\mathcal{P},p) \]

\[ =S_{k,\alpha}(\mathcal{A};\mathcal{P},w,\ln^{\alpha}(\xi^{2}))-\sum_{w\leq p<z}S_{k,\alpha}(\mathcal{A}_{p};\mathcal{P},w,\frac{1}{2}\ln^{\alpha}(\xi^{2})) -\frac{1}{2}\sum_{w\leq p<z}\ln^{\alpha}(\xi^{2})S_{0}(\mathcal{A}_{p^{2}};\mathcal{P},p) \]

\end{proof}	

\subsection{Lemma 1.12}
\ Suppose  $0=k_{0}\leq k_{1}<\cdots < k_{n}= 2^{\alpha}$,  $k_{n+2}=\frac{4^{\alpha}}{2}$, $k_{n+3}=\frac{4.5^{\alpha}}{2}$, $k_{n+4}=\frac{5^{\alpha}}{2}$, $k_{n+5}=\frac{5.5^{\alpha}}{2}$. \\
\[ k_{n+1}=\min\left(\frac{3^{\alpha}}{2},\frac{2^{\alpha}}{\frac{2^{\alpha}}{10^{\alpha}}+(1-\frac{1}{10})^{\alpha}}\right) \]

\ For example, when $\alpha=2$, and $n=16$, we have \\

\ $k_{0}=0$, $k_{h}=0.25 h$, $(0<h \leq n)$ \\

\ $k_{n+1}=4.5$, $k_{n+2}=8$, $k_{n+3}=10.125$, $k_{n+4}=12.5$, $k_{n+5}=15.125$ \\

\ When $\alpha=3$, and $n=16$, we have \\

\ $k_{0}=0$, $k_{h}=0.5 h$, $(0<h \leq n)$ \\

\ $k_{n+1}=10.85482$, $k_{n+2}=32$, $k_{n+3}=45.5625$, $k_{n+4}=62.5$, $k_{n+5}=83.1875$ \\

\ Two function serials $F_{\alpha}^{(i)}(k_{l},u)$ and $f_{\alpha}^{(i)}(k_{l},u)$, $(i=0,1,\cdots)$ are defined as

$$
\begin{array}{ll}
	F_{\alpha}^{(0)}(0,u)=\frac{2F(2)}{u}, & 0<u \leq 1 \\
	F_{\alpha}^{(0)}(0,u)=F(u), &  u\geq 1 \\
	F_{\alpha}^{(0)}(k_{n},u)=\frac{u\tilde{F}(u)}{u+\frac{k_{n}}{\alpha u^{\alpha-1}}} & 2 \leq u\leq 4 \\
	F_{\alpha}^{(0)}(k_{n},u)=\frac{\left(uF(u)+\frac{k_{n}}{\alpha u^{\alpha-1}}F(u-1)\right)}{u+\frac{k_{n}}{\alpha u^{\alpha-1}}} & u\geq 4 \\
	F_{\alpha}^{(0)}(k_{l},u)=\frac{(1-\frac{k_{l}}{k_{n}})uF_{\alpha}^{(0)}(0,u)+\frac{k_{l}}{k_{n}}(u+\frac{k_{n}}{\alpha u^{\alpha-1}})F_{\alpha}^{(0)}(k_{n},u)}{u+\frac{k_{l}}{\alpha u^{\alpha-1}}} & u\geq 2, \ \ 0 <k_{l}<k_{n} \\
	F_{\alpha}^{(0)}(k_{l},u)=\frac{(2+\frac{k_{l}}{\alpha 2^{\alpha-1}})F_{\alpha}^{(0)}(k_{l},2)}{u+\frac{k_{l}}{\alpha u^{\alpha-1}}} & k_{l}^{\frac{1}{\alpha}}< u \leq 2, u\geq 1; 0\leq k_{l} \leq k_{n} \\
\\
	f_{\alpha}^{(0)}(0,u)=0 & 0< u < 1 \\
	f_{\alpha}^{(0)}(0,u)=f(u) & u\geq 1 \\
	f_{\alpha}^{(0)}(k_{n+1},u)=\frac{uf(u)+\frac{k_{n+1}}{\alpha u^{\alpha-1}}f(u-1)}{u+\frac{k_{n+1}}{\alpha u^{\alpha-1}}} & u\geq 2 \\
	f_{\alpha}^{(0)}(k_{n+1},u)=0 & 0 < u < 2 \\
	f_{\alpha}^{(0)}(k_{l},u)=\frac{(1-\frac{k_{l}}{k_{n+1}})uf_{\alpha}^{(0)}(0,u)+\frac{k_{l}}{k_{n+1}}(u+\frac{k_{n+1}}{\alpha u^{\alpha-1}})f_{\alpha}^{(0)}(k_{n+1},u)}{u+\frac{k_{l}}{\alpha u^{\alpha-1}}} & u\geq 1, \ \ 0<k_{l}\leq k_{n} \\
\end{array}
$$
\\

\subsection{Lemma 1.13}
\ Suppose $\alpha \geq 2$, $0\leq k_{l}\leq k_{n}$, and $k_{l}^{\frac{1}{\alpha}}<u=\frac{\ln(\xi^{2})}{\ln(z)}$, $i-1=0$ we have:
 
\begin{equation}
S_{k_{l},\alpha}\left(\mathcal{A};\mathcal{P},z,\ln^{\alpha}(\xi^{2})\right) \leq Xe^{-\gamma}C(\omega)\ln^{\alpha-1}(\xi^{2})(u+\frac{k_{l}}{\alpha u^{\alpha-1}})F_{\alpha}^{(i-1)}(k_{l},u)(1+O(\frac{1}{\ln^{\frac{1}{14}}(\xi)})) 
\label{StartUp}
\end{equation}

\[ +\ln^{\alpha}(\xi^{2})\sum_{d|\mathcal{P}(z),d<\xi^{2}}3^{v_{1}(d)}|r_{d}| \]

\begin{equation}
S_{k_{l},\alpha}\left(\mathcal{A};\mathcal{P},z,\ln^{\alpha}(\xi^{2})\right) \geq Xe^{-\gamma}C(\omega)\ln^{\alpha-1}(\xi^{2})(u+\frac{k_{l}}{\alpha u^{\alpha-1}})f_{\alpha}^{(i-1)}(k_{l},u)(1+O(\frac{1}{\ln^{\frac{1}{14}}(\xi)}))
\label{StartDown}
\end{equation}

\[ -\ln^{\alpha}(\xi^{2})\sum_{d|\mathcal{P}(z),d<\xi^{2}}3^{v_{1}(d)}|r_{d}| \]
\begin{proof}
\ According to $Lemma 10$, $Corollary 1.10.1$ and $Lemma 1.7$ we know when $k_{l}=0$, or $k_{l}=k_{n}$ equation (2.20) is valid, and when $k_{l}=0$, or $k_{l}=k_{n+1}$  equation (2.21) is valid. When $0<k_{l}<k_{n}$, Suppose $u=\frac{\ln(\xi^{2})}{\ln(z)}$
\[ S_{k_{l},\alpha}(\mathcal{A};\mathcal{P},z,\ln^{\alpha}(\xi^{2}))=(1-\frac{k_{l}}{k_{n}})S_{0,\alpha}(\mathcal{A};\mathcal{P},z,\ln^{\alpha}(\xi^{2}))+\frac{k_{l}}{k_{n}}S_{k_{n},\alpha}(\mathcal{A};\mathcal{P},z,\ln^{\alpha}(\xi^{2}))\]
\[ \leq (1-\frac{k_{l}}{k_{n}})\ln^{\alpha-1}(\xi^{2})uF_{\alpha}^{(0)}(0,u)(1+o(\frac{1}{\ln^{\frac{1}{14}}(\xi)})) \]
\[ +\frac{k_{l}}{k_{n}}\ln^{\alpha-1}(\xi^{2})(u+\frac{k_{n}}{\alpha u^{\alpha-1}})F_{\alpha}^{(0)}(k_{n},u)(1+o(\frac{1}{\ln^{\frac{1}{14}}(\xi)})) \]
\[ +\ln^{\alpha}(\xi^{2})\sum_{d|\mathcal{P}(z),d<\xi^{2}}3^{v_{1}(d)}|r_{d}| \]
\[ =(u+\frac{k_{l}}{\alpha u^{\alpha-1}})F_{\alpha}^{(0)}(k_{l},u)(1+o(\frac{1}{\ln^{\frac{1}{14}}(\xi)}))+\ln^{\alpha}(\xi^{2})\sum_{d|\mathcal{P}(z),d<\xi^{2}}3^{v_{1}(d)}|r_{d}| \]

\ Proof of equation (2.21) is the same as equation (2.20)

\end{proof}
\subsection{Lemma 1.14}
\ Suppose $\ln^{\alpha+c_{1}+2}(\xi^{2})<w<z$, $\frac{\ln(\xi^{2})}{\ln(w)}=v$,$\frac{\ln(\xi^{2})}{\ln(z)}=u$ \\

\ If for any $i-1 \geq 0$ equation (2.20) and (2.21) are correct, when $l\leq n+1$, suppose

\[ u_{0}=\max\left(\min(3,u_{l}),2\right)\]
\ Where $u_{l}$ is the solution of equation
\[ k_{l}(t)=k_{l}\frac{(t-1)^{\alpha}}{t^{\alpha}-k_{l}}=k_{n} \]

\ Define the function:
\[ \hat{f}_{\alpha,1}^{(i)}(k_{l},u,v)=0, \ \ 1\leq u \leq u_{0} \]
\ And when $u_{0}< u<v$ \\
\[ \hat{f}_{\alpha,1}^{(i)}(k_{l},u,v)=\frac{(v+\frac{k_{l}}{\alpha v^{\alpha-1}})f_{\alpha}^{(i-1)}(k_{l},v)-\int_{u}^{v}\frac{t^{\alpha}-k_{l}}{t^{\alpha}(t-1)}\left(t-1+\frac{k_{l}(t)}{\alpha(t-1)^{\alpha-1}}\right)\breve{F}_{\alpha}^{(i-1)}(k_{l}(t),t-1)dt}{u+\frac{k_{l}}{\alpha u^{\alpha-1}}} \]

\ Where 
\[  0\leq k_{l}(t)=k_{l}\frac{(t-1)^{\alpha}}{t^{\alpha}-k_{l}}\leq k_{n} \]
\ and the function
\[ \breve{F}_{\alpha}^{(i-1)}(k_{l}(t),t-1)=\frac{\beta(t-1+\frac{k_{l1}}{2(t-1)})F_{\alpha}^{(i-1)}(k_{l_{1}},t-1)+(1-\beta)(t-1+\frac{k_{l_{2}}}{2(t-1)})F_{\alpha}^{(i-1)}(k_{l_{2}},t-1)}{t-1+\frac{k_{l}(t)}{2(t-1)}} \]
\ Where
\[ k_{l_{1}}=\min_{0<m\leq n \atop k_{m}\geq k_{l}(t)}(k_{m}),\ \ k_{l_{2}}=\max_{0\leq m<n \atop k_{m}\leq k_{l}(t)}(k_{m}) \]
\[ \beta= \beta (t)=\frac{k_{l}(t)-k_{l_{2}}}{k_{l_{1}}-k_{l_{2}}}\]
\ Is the solution of equation

\[ \beta k_{l_{1}}+(1-\beta)k_{l_{2}}=k_{l}(t) \]
\ we have: 

\begin{equation}
 S_{k_{l},\alpha}(\mathcal{A};\mathcal{P},z,\ln^{\alpha}(\xi^{2}))\geq \max_{v\geq  u}\left(Xe^{-\gamma}C(\omega)\ln^{\alpha-1}(\xi^{2})(u+\frac{k_{l}}{\alpha u^{\alpha-1}})\hat{f}_{\alpha,1}^{(i)}(k_{l},u,v)(1+O(\frac{1}{\ln^{\frac{1}{14}}(\xi^{2})}))\right)
\end{equation}

\[ -\ln^{\alpha}(\xi^{2})\sum_{d|\mathcal{P}(z),d<\xi^{2}}3^{v_{1}(d)}|r_{d}| \]
\begin{proof}
\ setting: 
\[ k(p)=k_{l}\frac{\ln^{\alpha}(\frac{\xi^{2}}{p})}{\ln^{\alpha}(\xi^{2})-k_{l}\ln^{\alpha}(p)} \]
\[ u_{p}=\frac{\ln(\xi^{2}/p)}{\ln(p)}=\frac{\ln(\xi^{2})}{\ln(p)}-1 \]
\ Since 
\[ \sum_{w\leq p<z}\ln^{\alpha}(p)S_{0}(\mathcal{A}_{p^{2}};\mathcal{P},p)\leq O(\sum_{w\leq p<z}\frac{X\ln^{\alpha+c_{1}}(p)}{p^{2}})\leq O(\frac{X\ln^{\alpha+c_{1}}(\xi^{2})}{w}) \leq \frac{X}{\ln^{2}(\xi^{2})}\]
\ According to equation (2.8) on $Lemma 1.1$
\[ S_{k_{l},\alpha}(\mathcal{A};\mathcal{P},z,\ln^{\alpha}(\xi^{2}))=\ln^{\alpha}(\xi^{2})S_{0}(\mathcal{A};\mathcal{P},w)-\sum_{w\leq p<z}\frac{\ln^{\alpha}(\xi^{2})-k_{l} \ln^{\alpha}(p)}{\ln^{\alpha}(\frac{\xi^{2}}{p})}S_{k(p),\alpha}\left(\mathcal{A}_{p};\mathcal{P},p,\ln^{\alpha}(\frac{\xi^{2}}{p})\right)\]
\[ -k\sum_{z_{1}\leq p<z}\ln^{\alpha}(p)S_{0}(\mathcal{A}_{p^{2}};\mathcal{P},p) \]
\[ =\ln^{\alpha}(\xi^{2})S_{0}(\mathcal{A};\mathcal{P},w)-\sum_{w\leq p<z}\frac{\ln^{\alpha}(\xi^{2})-k_{l} \ln^{\alpha}(p)}{\ln^{\alpha}(\frac{\xi^{2}}{p})}S_{k(p),\alpha}\left(\mathcal{A}_{p};\mathcal{P},p,\ln^{\alpha}(\frac{\xi^{2}}{p})\right)+O(\frac{X}{\ln^{2}(\xi^{2})}) \]

\[ S_{k(p),\alpha}\left(\mathcal{A}_{p};\mathcal{P},p,\ln^{\alpha}(\frac{\xi^{2}}{p})\right)= \beta S_{k_{l_{1}},\alpha}\left(\mathcal{A}_{p};\mathcal{P},p,\ln^{\alpha}(\frac{\xi^{2}}{p})\right)+(1-\beta)S_{k_{l_{2}},\alpha}\left(\mathcal{A}_{p};\mathcal{P},p,\ln^{\alpha}(\frac{\xi^{2}}{p})\right) \]
\[ \leq \frac{X}{p}e^{-\gamma}C(\omega)\ln^{\alpha-1}(\frac{\xi^{2}}{p})\left(\beta (u_{p}+\frac{k_{l_{1}}}{\alpha u_{p}^{\alpha-1}})F_{\alpha}^{i-1}(k_{l_{1}},u_{p})+(1-\beta)(u_{p}+\frac{k_{l_{2}}}{\alpha u_{p}^{\alpha-1}})F_{\alpha}^{i-1}(k_{l_{2}},u_{p})\right)\times \]
\[ \left(1+O(\frac{1}{\ln^{\frac{1}{14}}(w)})\right)+\ln^{\alpha}(p)R_{p} \]
\[ =\frac{X}{p}e^{-\gamma}C(\omega)\ln^{\alpha-1}(\frac{\xi^{2}}{p})(u_{p}+\frac{k_{p}}{\alpha u_{p}^{\alpha-1}})\breve{F}_{\alpha}^{i-1}(k_{p},u_{p})\left(1+O(\frac{1}{\ln^{\frac{1}{14}}(p)})\right)+\ln^{\alpha}(p)R_{p} \]

\ Since 
\[ \ln^{\alpha}(\xi)-k_{l} \ln^{\alpha}(p)=\ln^{\alpha}(\xi)(1-\frac{k_{l}}{u_{p}^{\alpha}})> \ln^{\alpha}(\xi)\left(1-\frac{k_{l}}{(u_{0}-1)^{\alpha}}\right) \geq 0 \]
\ we have
\[ \sum_{w\leq p<z}\frac{\ln^{\alpha}(\xi)-k_{l} \ln^{\alpha}(p)}{\ln^{\alpha}(\frac{\xi^{2}}{p})}S_{k(p),\alpha}\left(\mathcal{A}_{p};\mathcal{P},p,\ln^{\alpha}(\frac{\xi^{2}}{p})\right)\]
\[ \leq Xe^{-\gamma}C(\omega)\sum_{w\leq p <z}\frac{\omega(p)(\ln^{\alpha}(\xi)-k_{l} \ln^{\alpha}(p))}{p\ln^{\alpha}(\frac{\xi^{2}}{p})}\ln^{\alpha-1}(\frac{\xi^{2}}{p})(u_{p}+\frac{k_{p}}{\alpha u_{p}^{\alpha-1}})\breve{F}_{\alpha}^{i-1}(k_{p},u_{p})\left(1+O(\frac{1}{\ln^{\frac{1}{14}}(w)})\right)\]
\[ +\sum_{w\leq p<z}\ln^{\alpha}(p)R_{p} \]
\[ =Xe^{-\gamma}C(\omega)\int_{w}^{z}\frac{\ln^{\alpha}(\xi)-k_{l} \ln^{\alpha}(t)}{t\ln(t)\ln^{\alpha}(\frac{\xi^{2}}{t})}\ln^{\alpha-1}(\frac{\xi^{2}}{t})(u_{t}+\frac{k_{l}(t)}{\alpha u_{t}^{\alpha-1}})\breve{F}_{\alpha}^{i-1}(k_{l}(t),u_{t})dt\left(1+O(\frac{1}{\ln^{\frac{1}{14}}(w)})\right)\]
\[ +\sum_{w\leq p<z}\ln^{\alpha}(p)R_{p} \]
\[ =Xe^{-\gamma}C(\omega)\int_{u}^{v}\ln^{\alpha-1}(\xi^{2})\frac{1-\frac{k_{l}}{t^{\alpha}}}{t(1-\frac{1}{t})}(t-1+\frac{k_{l}(t)}{\alpha (t-1)^{\alpha-1}})\breve{F}_{\alpha}^{i-1}(k_{l}(t),t-1)dt\left(1+O(\frac{1}{\ln^{\frac{1}{14}}(w)})\right)\]
\[ +\sum_{w\leq p<z}\ln^{\alpha}(p)R_{p} \]
\ Where
\[ \sum_{w\leq p<z}\ln^{\alpha}(p)R_{p}\leq \ln^{\alpha}(\xi^{2})R \]
\ we obtain
\[ S_{k_{l},\alpha}(\mathcal{A};\mathcal{P},z,\ln^{\alpha}(\xi^{2}))\geq Xe^{-\gamma}C(\omega)\ln^{\alpha-1}(\xi^{2})(v+\frac{k_{l}}{\alpha v^{\alpha-1}})f^{i-1}_{\alpha}(k_{l},v)+\ln^{\alpha}(\xi^{2})R\]
\[ -Xe^{-\gamma}C(\omega)\int_{u}^{v}\ln^{\alpha-1}(\xi^{2})\frac{1-\frac{k_{l}}{t^{\alpha}}}{t(1-\frac{1}{t})}(t-1+\frac{k_{l}(t)}{\alpha (t-1)^{\alpha-1}})\breve{F}_{\alpha}^{i-1}(k_{l}(t),t-1)dt\left(1+O(\frac{1}{\ln^{\frac{1}{14}}(w)})\right)\]
\[ +\sum_{w\leq p<z}\ln^{\alpha}(p)R_{p} \]
\[ =\ln^{\alpha-1}(\xi^{2})(u+\frac{k_{l}}{\alpha u^{\alpha-1}})f^{(i)}_{\alpha,1}(k_{l},u,v)\left(1+O(\frac{1}{\ln^{\frac{1}{14}}(w)})\right)+\ln^{\alpha}(\xi^{2})R \]

\end{proof}
\subsection{Lemma 1.15}
\ Same condition of Lemma (1.14), Suppose $v>u$, define function
\
\[ \hat{f}_{\alpha,2}^{(i)}(k_{l},u,v)=f_{\alpha}^{(i-1)}(k_{l},u), \ \ u \geq v \]

\[ \hat{f}_{\alpha,2}^{(i)}(k_{l},u,v)=\frac{1}{u+\frac{k_{l}}{\alpha u^{\alpha-1}}}\left((v+\frac{k_{l}}{\alpha v^{\alpha-1}})f_{\alpha}^{(i-1)}(k_{l},v)-\frac{1}{2}d_{1}(u,v)-d_{2}(u,v)\right) \ \ 1< u\leq v \]

\ Where
\[ d_{1}(u,v)=\int_{u}^{v}\frac{1}{t-1}\left(v-v/t+\frac{k_{l}(t)}{\alpha (v-v/t)^{\alpha-1}}\right)\breve{F}_{\alpha}(k_{l}(t),v-v/t)dt \]
\ Function $\breve{F}_{\alpha}(k_{l}(t),v-v/t)$ define same as $lemma 1.14$

\[ d_{2}(u,v)=uF_{\alpha}^{(i-1)}(0,u)\int_{u<t<v \atop t^{\alpha}\geq 2k_{l}}\frac{1}{t}(\frac{1}{2}-\frac{k_{l}}{t^{\alpha}})\frac{1}{1-1/t}dt \]
\[ =uF_{\alpha}^{(i-1)}(0,u)\int_{u_{1}}^{v}(\frac{1}{2}-\frac{k_{l}}{t^{\alpha}})\frac{1}{t-1}dt \]

\ $u_{1}=\max((2k_{l})^{\frac{1}{\alpha}},u)$, $(2k_{l})^{\frac{1}{\alpha}}$ is a solution of equation
\[ \frac{t^{\alpha}}{2}-k_{l}=0 \]

\ We have:

\begin{equation}
 S_{k_{l},\alpha}(\mathcal{A};\mathcal{P},z,\ln^{\alpha}(\xi^{2}))\geq \max_{v\geq  u}\left(X\ln^{\alpha-1}(\xi^{2})(u+\frac{k_{l}}{\alpha u^{\alpha-1}})\hat{f}_{\alpha,2}^{(i)}(k_{l},u,v)(1+O(\frac{1}{\ln^{\frac{1}{14}}(\xi^{2})}))\right)
\end{equation} 

\[ -\ln^{\alpha}(\xi^{2})\sum_{d|\mathcal{P}(z),d<\xi^{2}}3^{v_{1}(d)}|r_{d}| \]
\begin{proof}
\ Since $\ln^{\alpha+c_{1}+2}(\xi^{2})<w<z$, $v=\frac{\ln(\xi^{2})}{\ln(w)}$,$u=\frac{\ln(\xi^{2})}{\ln(z)}$,$v\geq 3$, $0<u<v$,  when $u\geq u_{1}$, 
\[ \frac{\ln^{\alpha}(\xi^{2})}{2}\geq k_{l}\ln^{\alpha}(p) \] 
\ according to $lemma 1.11$ we have
\[ S_{\alpha,k_{l}}(\mathcal{A};\mathcal{P},z,\ln^{\alpha}(\xi^{2}))\geq \ln^{\alpha}(\xi^{2})S_{0}(\mathcal{A};\mathcal{P},z)+k_{l}\sum_{2\leq p<w}\ln^{\alpha}(p)S_{0}(\mathcal{A}_{p},\mathcal{P}(p),z) \]
\[ +\frac{1}{2}\ln^{\alpha}(\xi^{2})\sum_{w\leq p<z}S_{0}(\mathcal{A}_{p};\mathcal{P},z)-\sum_{w\leq p<z, \atop \ln^{\alpha}(\xi^{2})>2k_{l}\ln^{\alpha}(p)}\left(\frac{1}{2}\ln^{\alpha}(\xi^{2})-k_{l}\ln^{\alpha}(p)\right)S_{0}(\mathcal{A}_{p};\mathcal{P},z) \]

\[ \geq S_{k_{l},\alpha}(\mathcal{A};\mathcal{P},w\ln^{\alpha}(\xi^{2})-\sum_{w\leq p<z}S_{k_{l},\alpha}(\mathcal{A}_{p};\mathcal{P},w,\frac{1}{2}\ln^{\alpha}(\xi^{2})) \]
\[ -\sum_{w\leq p<z, \atop \ln^{\alpha}(\xi^{2})>2k_{l}\ln^{\alpha}(p)}\left(\frac{1}{2}\ln^{\alpha}(\xi^{2})-k_{l}\ln^{\alpha}(p)\right)S_{0}(\mathcal{A}_{p};\mathcal{P},z) +O\left(\frac{X}{\ln^{2}(\xi^{2})}\right)\]
\ Where used
\[ \sum_{w\leq p<z}\ln^{\alpha}(p)S_{0}(\mathcal{A}_{p^{2}};\mathcal{P},p)\leq O(\sum_{w\leq p<z}\frac{X\ln^{\alpha+c_{1}}(p)}{p^{2}})\leq O(\frac{X\ln^{\alpha+c_{1}}(\xi^{2})}{w}) \leq \frac{X}{\ln^{2}(\xi^{2})}\]

\ It is easy to see that
\[ S_{k_{l},\alpha}(\mathcal{A};\mathcal{P},w\ln^{\alpha}(\xi^{2})\geq \ln^{\alpha-1}(\xi^{2})Xe^{-\gamma}C(\omega)vf_{\alpha}^{(i-1)}(k_{l},v)(1+O(\frac{1}{\ln^{\frac{1}{14}}(\xi^{2})})) \]

\ Suppose 
\[ 0 \leq k_{l}(p)=2k_{l}\frac{\ln^{\alpha}(\xi^{2}/p)}{\ln^{\alpha}(\xi^{2})}\leq k_{n} \]

\[ \sum_{w\leq p<z}S_{k_{l},\alpha}(\mathcal{A}_{p};\mathcal{P},w,\frac{1}{2}\ln^{\alpha}(\xi^{2}))=\sum_{w\leq p<z}\frac{1}{2}\frac{\ln^{\alpha}(\xi^{2})}{\ln^{\alpha}(\xi^{2}/p)}S_{k_{l}(p),\alpha}(\mathcal{A}_{p};\mathcal{P},w,\ln^{\alpha}(\xi^{2}/p)) \]
\[ \leq \frac{1}{2} Xe^{-\gamma}C(\omega)\sum_{w\leq p<z}\frac{w(p)}{p}\frac{\ln^{\alpha}(\xi^{2})}{\ln^{\alpha}(\xi^{2}/p)}\ln^{\alpha-1}(\xi^{2}/p)\left(\frac{\ln(\xi^{2}/p)}{\ln(w)}+k_{l}(p)(\frac{\ln(w)}{\ln(\xi^{2}/p)})^{\alpha-1}\right)\breve{F}_{\alpha}(k_{l}(p),\frac{\ln(\xi^{2}/p)}{\ln(w)})\times \]
\[ \left(1+O(\frac{1}{\ln^{\frac{1}{14}}(w)})\right) +\sum_{w\leq p <z}\ln^{\alpha}(p)R_{p}\]

\[ =\frac{1}{2} Xe^{-\gamma}C(\omega)\ln^{\alpha-1}(\xi^{2})\int_{u}^{v}\frac{1}{t(1-\frac{1}{t})}(v-v/t+\frac{k_{l}(t)}{\alpha (v-v/t)^{\alpha-1}})\breve{F}_{\alpha}(k_{l}(t),v-v/t)dt\left(1+O(\frac{1}{\ln^{\frac{1}{14}}(w)})\right)\]
\[ +\sum_{w\leq p <z}\ln^{\alpha}(p)R_{p} \]
\[ =\frac{1}{2} Xe^{-\gamma}C(\omega)\ln^{\alpha-1}(\xi^{2})d_{1}(u,v)(1+O(\frac{1}{\ln^{\frac{1}{14}}(w)})+\sum_{w\leq p <z}\ln^{\alpha}(p)R_{p} \]
\ and 
\[ \sum_{w\leq p<z, \atop \ln^{\alpha}(\xi^{2})>2k_{l}\ln^{\alpha}(p)}\left(\frac{1}{2}\ln^{\alpha}(\xi^{2})-k_{l}\ln^{\alpha}(p)\right)S_{0}(\mathcal{A}_{p};\mathcal{P},z) \]
\[ \leq \sum_{w\leq p<z, \atop \ln^{\alpha}(\xi^{2})>2k_{l}\ln^{\alpha}(p)}\frac{\omega(p)}{p}\left(\frac{1}{2}\ln^{\alpha}(\xi^{2})-k_{l}\ln^{\alpha}(p)\right)Xe^{-\gamma}C(\omega)\ln^{-1}(z)F_{\alpha}^{(i-1)}(0,\frac{\ln(\xi^{2}/p)}{\ln(z)})(1+O(\frac{1}{\ln^{\frac{1}{14}}(w)}))\]
\[ +\sum_{w\leq p<z}\ln^{\alpha}(p)R_{p} \]
\[ \leq Xe^{-\gamma}C(\omega)\ln^{\alpha-1}(\xi^{2})uF_{\alpha}^{(i-1)}(0,u)\int_{u<t<v \atop t^{\alpha}\geq 2k_{l}}\frac{1}{t}(\frac{1}{2}-\frac{k_{l}}{t^{\alpha}})\frac{1}{1-1/t}dt(1+O(\frac{1}{\ln^{\frac{1}{14}}(w)}))+\sum_{w\leq p<z}\ln^{\alpha}(p)R_{p} \]
\[ =Xe^{-\gamma}C(\omega)\ln^{\alpha-1}(\xi^{2})d_{2}(u,v)(1+O(\frac{1}{\ln^{\frac{1}{14}}(w)}))+\sum_{w\leq p<z}\ln^{\alpha}(p)R_{p} \]
\ Where used equation:
\[ (u-u/t)F_{\alpha}(0,u-u/t)\leq uF_{\alpha}(0,u), \ \ t \geq u \geq 1, \ \ t>1 \]
 
\begin{equation}
uF_{\alpha}(0,u-u/t) \leq uF_{\alpha}(0,u)/(1-1/t), \ \ t \geq u \geq 1, \ \ t>1
\end{equation}

\ Combining these tree equations we obtain $Lemma1.15$.
\end{proof}

\subsection{Lemma 1.16}
\ Same condition as Lemma (1.14). Suppose $0\leq k_{m}<k_{l}<k_{h}\leq k_{n}$ define function

\[ \hat{f}_{\alpha,3}^{(i)}(k_{l},u)=\frac{1}{u+\frac{k_{l}}{\alpha u^{\alpha-1}}}\max_{0\leq h<l \atop l<m\leq n+1}\left(\beta(u+\frac{k_{h}}{\alpha u^{\alpha-1}})f_{\alpha}^{(i-1)}(k_{h},u)+(1-\beta)(u+\frac{k_{m}}{\alpha u^{\alpha-1}})f_{\alpha}^{(i-1)}(k_{m},u)\right) \]
 
\[ \beta=\frac{k_{m}-k_{l}}{k_{m}-k_{h}} \]
\ Is the solution of the function \\
\[ \beta(k_{h})+(1-\beta)k_{m}=k_{l} \]
\ We have:

\begin{equation}
 S_{k_{l},\alpha}(\mathcal{A};\mathcal{P},z,\ln^{\alpha}(\xi^{2}))\geq \max_{0\leq k_{m}<k_{l}<k_{h}\leq k_{n}}\left(Xe^{-\gamma}C(\omega)\ln^{\alpha-1}(\xi^{2})(u+\frac{k_{l}}{\alpha u^{\alpha-1}})\hat{f}_{\alpha,3}^{(i)}(k_{l},u)(1+O(\frac{1}{\ln^{\frac{1}{14}}(z)}))\right)
\end{equation} 

\[ -\ln^{\alpha}(\xi^{2})\sum_{d|\mathcal{P}(z),d<\xi^{2}}3^{v_{1}(d)}|r_{d}| \]
\begin{proof}
\ For any $0\leq k_{m}<k_{l}<k_{h}\leq k_{n}$, and
\[ \beta=\frac{k_{m}-k_{l}}{k_{m}-k_{h}} \]
\ we have
\[ S_{k_{l},\alpha}(\mathcal{A};\mathcal{P},z,\ln^{\alpha}(\xi^{2}))=\beta S_{k_{h},\alpha}(\mathcal{A};\mathcal{P},z,\ln^{\alpha}(\xi^{2}))+(1-\beta)S_{k_{m},\alpha}(\mathcal{A};\mathcal{P},z,\ln^{\alpha}(\xi^{2})) \]
\[ \geq Xe^{-\gamma}C(\omega)\ln^{\alpha-1}(\xi^{2})\left(\beta(u+\frac{k_{h}}{\alpha u^{\alpha-1}})f_{\alpha}^{(i-1)}(k_{h},u)+(1-\beta)(u+\frac{k_{m}}{\alpha u^{\alpha-1}})f_{\alpha}^{(i-1)}(k_{m},u)\right)\times \]
\[ \left(1+O(\frac{1}{\ln^{\frac{1}{14}}(z)})\right) -\ln^{\alpha}(\xi^{2})R\]
\[ =Xe^{-\gamma}C(\omega)\ln^{\alpha-1}(\xi^{2})(u+\frac{k_{l}}{\alpha u^{\alpha-1}})\hat{f}_{\alpha,3}^{(i)}(k_{l},u)\left(1+O(\frac{1}{\ln^{\frac{1}{14}}(z)})\right) -\ln^{\alpha}(\xi^{2})R\]

\end{proof}
\subsection{Lemma 1.17}
\ Same condition of Lemma (1.14), suppose $0\leq k_{l}<k_{n}$, define function

\[ \hat{f}_{\alpha,4}^{(i)}(k_{l},u)=\frac{1}{u+\frac{k_{l}}{\alpha u^{\alpha-1}}}\max_{h>l}\left((u+\frac{k_{h}}{\alpha u^{\alpha-1}})f_{\alpha}^{i-1}(k_{h},u)-(k_{h}-k_{l})uF_{\alpha}^{i-1}(0,u)\int_{0}^{\frac{1}{u}}\frac{t^{\alpha-1}}{1-t}dt\right) \]

\ We have:

\begin{equation}
 S_{k_{l},\alpha}(\mathcal{A};\mathcal{P},z,\ln^{\alpha}(\xi^{2}))\geq X\ln^{\alpha-1}(\xi^{2})(u+\frac{k_{l}}{\alpha u^{\alpha-1}})\hat{f}_{\alpha,4}^{(i)}(k_{l},u)\left(1+O(\frac{1}{\ln^{\frac{1}{14}}(\xi^{2})})\right)
\end{equation} 

\[ -\ln^{\alpha}(\xi^{2})\sum_{d|\mathcal{P}(z),d<\xi^{2}}3^{v_{1}(d)}|r_{d}| \]
\begin{proof}
\ Suppose $k_{h}>k_{l}$
\[ S_{k_{l},\alpha}(\mathcal{A};\mathcal{P},z,\ln^{\alpha}(\xi^{2}))=S_{k_{h},\alpha}(\mathcal{A};\mathcal{P},z,\ln^{\alpha}(\xi^{2})) -\sum_{2\leq p<z}(k_{h}-k_{l})\ln^{\alpha}(p)S_{0}(\mathcal{A};\mathcal{P}(p),z) \]
\[ \geq Xe^{-\gamma}C(\omega)\ln^{\alpha-1}(\xi^{2})(u+\frac{k_{h}}{\alpha u^{\alpha-1}})f_{\alpha}^{(i-1)}(k_{h},u)\left(1+O(\frac{1}{\ln^{\frac{1}{14}}(z)})\right) -\ln^{\alpha}(\xi^{2})R\]
\[ -Xe^{-\gamma}C(\omega)(k_{h}-k_{l})\ln^{\alpha-1}(\xi^{2})\sum_{2\leq p<z}\frac{\omega(p)\ln^{\alpha}(p)}{p\ln^{\alpha}(\xi^{2})}uF_{\alpha}^{(i-1)}(0,\frac{\ln(\xi^{2}/p)}{\ln(z)})\left(1+O\frac{1}{\ln^{\frac{1}{14}}(p)}\right) -\sum_{2\leq p<z}\ln^{\alpha}(p)R_{p}\]
\[ =Xe^{-\gamma}C(\omega)\ln^{\alpha-1}(\xi^{2})(u+\frac{k_{h}}{\alpha u^{\alpha-1}})f_{\alpha}^{(i-1)}(k_{h},u)\left(1+O(\frac{1}{\ln^{\frac{1}{14}}(z)})\right) \]
\[ -Xe^{-\gamma}C(\omega)\ln^{\alpha-1}(\xi^{2})(k_{h}-k_{l})uF_{\alpha}^{(i-1)}(0,u)\int_{0}^{\frac{1}{u}}\frac{t^{\alpha-1}}{1-t}dt\left(1+O(\frac{1}{\ln^{\frac{1}{14}}(z)})\right) -\ln^{\alpha}(\xi^{2})R\]
\[ =Xe^{-\gamma}C(\omega)\ln^{\alpha-1}(\xi^{2})(u+\frac{k_{l}}{\alpha u^{\alpha-1}})\hat{f}_{\alpha,4}^{(i)}(k_{l},u)\left(1+O(\frac{1}{\ln^{\frac{1}{14}}(z)})\right) -\ln^{\alpha}(\xi^{2})R \]
\ The last step is performed using equation (2.24)
\end{proof}

\subsection{Lemma 1.18}

\ Same condition as Lemma (1.14), Suppose  $n+1< l\leq n+4$\\
\ Set $4\leq v=(2k_{l})^{\frac{1}{\alpha}}\leq 5$
\[ \hat{f}_{\alpha,1}^{(i)}(k_{l},u)=\frac{1}{u+\frac{k_{l}}{\alpha u^{\alpha-1}}}\left(vf_{0}^{(i-1)}(v)-\frac{v}{2}\int_{u}^{v}\frac{F_{\alpha}^{i-1}(0,v-v/t)}{t}dt\right) \]

\[ \hat{f}_{\alpha,2}^{(i)}(k_{l},u)=\frac{1}{u+\frac{k_{l}}{\alpha u^{\alpha-1}}}\left((u+\frac{k_{l-1}}{\alpha u^{\alpha-1}})f_{\alpha}^{i-1}(k_{l-1},u)\right) \]

\[ f_{\alpha}^{(i)}(k_{l},u)=\max \left(f_{\alpha}^{(i-1)}(k_{l},u),\hat{f}_{\alpha,1}^{(i)}(k_{l},u),\hat{f}_{\alpha,2}^{(i)}(k_{l},u)\right) \]

\ We have:

\begin{equation}
 S_{k_{l},\alpha}(\mathcal{A};\mathcal{P},z,\ln^{\alpha}(\xi^{2}))\geq X\ln^{\alpha-1}(\xi^{2})(u+\frac{k_{l}}{\alpha u^{\alpha-1}})f_{\alpha}^{(i)}(k_{l},u)(1+O(\frac{1}{\ln^{\frac{1}{14}}(\xi^{2})}))
\end{equation} 

\[ -\ln^{\alpha}(\xi^{2})\sum_{d|\mathcal{P}(z),d<\xi^{2}}3^{v_{1}(d)}|r_{d}| \]

\begin{proof}
\ Suppose $\ln^{\alpha+c_{1}+2}(\xi^{2})\leq w<z$, $v=\frac{\ln(\xi^{2})}{\ln(w)}\leq (2k_{l})^{\frac{1}{\alpha}}\leq 5$, when $w\leq p<z$ we have that
\[ \frac{1}{2}\ln^{\alpha}(\xi^{2})\leq k_{l}\ln^{\alpha}(p)\]
\ So that
\[ S_{k_{l},\alpha}(\mathcal{A};\mathcal{P},z,\ln^{\alpha}(\xi^{2}))\geq \ln^{\alpha}(\xi^{2})S_{0}(\mathcal{A};\mathcal{P},z)+\frac{1}{2}\ln^{\alpha}(\xi^{2})\sum_{w\leq p<z}S_{0}(\mathcal{A}_{p};\mathcal{P}(p),z) \]
\[ \geq \ln^{\alpha}(\xi^{2})S_{0}(\mathcal{A};\mathcal{P},w)-\frac{1}{2}\ln^{\alpha}(\xi^{2})\sum_{w\leq p<z}S_{0}(\mathcal{A}_{p};\mathcal{P},w) +O(\frac{X}{\ln^{2}(X)})\]

\[ \geq Xe^{-\gamma}C(\omega)\ln^{\alpha-1}(\xi^{2})\left(vf_{\alpha}^{(i-1)}(0,v)-\frac{1}{2}v\int_{u}^{v}\frac{F_{\alpha}^{(i-1)}(0,v-v/t)}{t} dt\right)\left(1+O(\frac{1}{\ln^{\frac{1}{14}}(\xi^{2})})\right) -\ln^{\alpha}(\xi^{2})R\]
\[ =Xe^{-\gamma}C(\omega)\ln^{\alpha-1}(\xi^{2})(u+\frac{k_{l}}{\alpha u^{\alpha-1}})\hat{f}_{\alpha,1}^{(i)}(k_{l},u)\left(1+O(\frac{1}{\ln^{\frac{1}{14}}(\xi^{2})})\right) -\ln^{\alpha}(\xi^{2})R\]
\ On the other hand we have: when $k_{l}>0$
\[ S_{k_{l},\alpha}(\mathcal{A};\mathcal{P},z,\ln^{\alpha}(\xi^{2}))\geq S_{k_{l-1},\alpha}(\mathcal{A};\mathcal{P},z,\ln^{\alpha}(\xi^{2})) \]
\[ \geq Xe^{-\gamma}C(\omega)\ln^{\alpha-1}(\xi^{2})(u+\frac{k_{l-1}}{\alpha u^{\alpha-1}})f_{\alpha}^{(i-1)}(k_{l-1},u)\left(1+O(\frac{1}{\ln^{\frac{1}{14}}(\xi^{2})})\right) -\ln^{\alpha}(\xi^{2})R\]
\[ =Xe^{-\gamma}C(\omega)\ln^{\alpha-1}(\xi^{2})(u+\frac{k_{l}}{\alpha u^{\alpha-1}})\frac{(u+\frac{k_{l-1}}{\alpha u^{\alpha-1}})}{u+\frac{k_{l}}{\alpha u^{\alpha-1}}}f_{\alpha}^{(i-1)}(k_{l-1},u)\left(1+O(\frac{1}{\ln^{\frac{1}{14}}(\xi^{2})})\right) -\ln^{\alpha}(\xi^{2})R\]
\[ =Xe^{-\gamma}C(\omega)\ln^{\alpha-1}(\xi^{2})(u+\frac{k_{l}}{\alpha u^{\alpha-1}})\hat{f}_{\alpha,2}^{(i)}(k_{l},u)\left(1+O(\frac{1}{\ln^{\frac{1}{14}}(\xi^{2})})\right) -\ln^{\alpha}(\xi^{2})R\]

\end{proof}

\subsection{Lemma 1.19}
\ Same condition as Lemma (1.14), 

\ When $i>0$, Suppose  $0\leq l\leq n$\\

\[ u_{0}=\max\left(\min(3,u_{l}),2\right), \ \ u_{0}\leq u<v \]
\ Where $u_{l}$ is the solution of the equation
\[ k_{l}\frac{(t-1)^{\alpha}}{t^{\alpha}-k_{l}}= k_{n} \]

\ define the function \\

\[ \hat{F}_{\alpha,1}^{(i)}(k_{l},u,v)=\frac{(v+\frac{k_{l}}{\alpha v^{\alpha-1}})F_{\alpha}^{(i-1)}(k_{l},v)-\int_{u}^{v}\frac{t^{\alpha}-k_{l}}{t^{\alpha}(t-1)}\left(t-1+\frac{k_{l}(t)}{\alpha (t-1)^{\alpha-1}}\right)\breve{f}_{\alpha}^{(i-1)}(k_{l}(t),t-1)dt}{u+\frac{k_{l}}{\alpha u^{\alpha-1}}} \]

\ Where 
\[  0\leq k_{l}(t)=k_{l}\frac{(t-1)^{\alpha}}{t^{\alpha}-k_{l}}\leq k_{n} \]
\ and 
\[ \breve{f}_{\alpha}^{(i-1)}(k_{l}(t),t-1)=\frac{\beta(t-1+\frac{k_{l1}}{\alpha(t-1)^{\alpha-1}})f_{\alpha}^{(i-1)}(k_{l_{1}},t-1)+(1-\beta)(t-1+\frac{k_{l2}}{\alpha(t-1)^{\alpha-1}})f_{\alpha}^{(i-1)}(k_{l_{2}},t-1)}{t-1+\frac{k_{l}(t)}{\alpha(t-1)^{\alpha-1}}} \]

\[ k_{l_{1}}=\min_{0<m\leq n \atop k_{m}\geq k_{l}(t)}(k_{m}),\ \ k_{l_{2}}=\max_{0\leq m<n \atop k_{m}\leq k_{l}(t)}(k_{m}) \]
\[ \beta= \beta(t)=\frac{k_{l}(t)-k_{l_{2}}}{k_{l_{1}}-k_{l_{2}}}\]
\ Is the solution of the equation

\begin{equation}
\beta k_{l_{1}}+(1-\beta)k_{l_{2}}=k_{l}(t) 
\end{equation}

\[ \hat{F}_{\alpha,1}^{(i)}(k_{l},u,v)=(u_{0}+\frac{k_{l}}{2u_{0}})F_{\alpha}^{(i-1)}(k_{l},u_{0})/(u+\frac{k_{l}}{2u}), \ \ k_{l}^{\frac{1}{\alpha}}\leq u <u_{0} \]
\ We have:

\begin{equation}
 S_{k_{l},\alpha}(\mathcal{A};\mathcal{P},z,\ln^{\alpha}(\xi^{2}))\leq Xe^{-\gamma}C(\omega)\ln^{\alpha-1}(\xi^{2})(u+\frac{k_{l}}{\alpha u^{\alpha-1}})\hat{F}_{\alpha,1}^{(i)}(k_{l},u,v)(1+O(\frac{1}{\ln^{\frac{1}{14}}(\xi^{2})}))
\end{equation} 

\[ +\ln^{\alpha}(\xi^{2})\sum_{d|\mathcal{P}(z),d<\xi^{2}}3^{v_{1}(d)}|r_{d}| \]
\begin{proof}
\ When $u>u_{0}$, $k_{l}(t)\leq k_{n}$, in equation (2.28) $\beta$ will have a positive solution. Proof of equation (2.29) is the same as the proof of $lemma 1.14$. \\
\ When $k_{l}^{\frac{1}{\alpha}}\leq u\leq u_{0}=\frac{\ln(\xi^{2})}{\ln(w)}$ according to equation (2.8), we have
\[ S_{k_{l},\alpha}(\mathcal{A};\mathcal{P},z,\ln^{\alpha}(\xi^{2})=S_{k_{l},\alpha}(\mathcal{A};\mathcal{P},w,\ln^{\alpha}(\xi^{2}) \]
\[ -\sum_{w\leq p<z}\frac{\ln^{\alpha}(\xi^{2})-k_{l}\ln^{\alpha}(p)}{\ln^{\alpha}(\frac{\xi^{2}}{p})}S_{k_{p},\alpha}(\mathcal{A}_{p};\mathcal{P},z,\ln^{\alpha}(\frac{\xi^{2}}{p})) +O\left(\frac{X}{\ln^{2}(X)}\right)\]
\ Since 
\[ \frac{\ln(\xi^{2})}{\ln(p)}\geq u \geq k_{l}^{\frac{1}{\alpha}} \]
\[ \ln^{\alpha}(\xi^{2})-k_{l}\ln^{\alpha}(p)\geq 0 \]
\ and
\[ \ln^{\alpha}(\frac{\xi^{2}}{p})>0 \]
\ The sum on the right
\[ \sum_{w\leq p<z}\frac{\ln^{\alpha}(\xi^{2})-k_{l}\ln^{\alpha}(p)}{\ln^{\alpha}(\frac{\xi^{2}}{p})}S_{k_{p},\alpha}(\mathcal{A}_{p};\mathcal{P},z,\ln^{\alpha}(\frac{\xi^{2}}{p}))\geq 0 \]

\ We obtain
\[ S_{k_{l},\alpha}(\mathcal{A};\mathcal{P},z,\ln^{\xi^{2}}) \leq S_{k_{l},\alpha}(\mathcal{A};\mathcal{P},w,\ln^{\xi^{2}}) \]
\[ \leq Xe^{-\gamma}C(\omega)\ln^{\alpha-1}(\xi^{2})(u_{0}+\frac{u_{0}}{\alpha u_{0}^{\alpha-1}})F_{\alpha}^{(i-1)}(k_{l},u_{0})\left(1+O(\frac{1}{\ln^{\frac{1}{14}}(\xi^{2})})\right) +\ln^{\alpha}(\xi^{2})R\]
\[ =Xe^{-\gamma}C(\omega)\ln^{\alpha-1}(\xi^{2})(u+\frac{k_{l}}{\alpha u^{\alpha-1}})\hat{F}_{\alpha,1}^{(i)}(k_{l},u,v)\left(1+O(\frac{1}{\ln^{\frac{1}{14}}(\xi^{2})})\right) +\ln^{\alpha}(\xi^{2})R\]

\end{proof}
\subsection{Lemma 1.20}
\ Same condition as Lemma (1.14), define the function

\[ \hat{F}_{\alpha,2}^{(i)}(k_{l},u)=\frac{\min_{h>l}\left((u+\frac{k_{h}}{\alpha u^{\alpha-1}})F_{\alpha}^{i-1}(k_{h},u)-(k_{h}-k_{l})u\int_{0}^{\frac{1}{u}}t^{\alpha-1}f_{\alpha}^{i-1}(0,u-ut)dt\right)}{u+\frac{k_{l}}{\alpha u^{\alpha-1}}}
\]
\ We have

\begin{equation}
 S_{k_{l},\alpha}(\mathcal{A};\mathcal{P},z,\ln^{\alpha}(\xi^{2}))\leq Xe^{-\gamma}C(\omega)\ln^{\alpha-1}(\xi^{2})(u+\frac{k_{l}}{\alpha u^{\alpha-1}})\hat{F}_{\alpha,2}^{(i)}(k_{l},u)(1+O(\frac{1}{\ln^{\frac{1}{14}}(\xi^{2})}))
\end{equation} 

\[ +\ln^{\alpha}(\xi^{2})\sum_{d|\mathcal{P}(z),d<\xi^{2}}3^{v_{1}(d)}|r_{d}| \]
\begin{proof}
\ Proof of this lemma is the same as the proof of $lemma 1.17$.
\end{proof}

\subsection{Lemma 1.21}
\ Suppose $\ln^{\alpha+c_{1}+2}(\xi^{2})<w<z$, $v=\frac{\ln(\xi^{2})}{\ln(w)}$, $u=\frac{\ln(\xi^{2})}{\ln(z)}$, $1 \leq u \leq t \leq v$, $1 \leq a \leq k_{n}$. Continuum function $k(a,t)$ satisfy:
\[ 0<k(a,t)\leq k_{n+1} \]

\ Define functions $k(t)$, $R(t)$, $H(\alpha,v,t)$ as
\[ k(t) = \min_{t \leq t_{1} \leq v}k(a,t_{1}) \]
\[ R(t)=\frac{1-\frac{k(a,t)}{t^{\alpha}}}{(1-\frac{1}{t})^{\alpha}} \]
\[ H(\alpha,v,t)=(1-1/t)^{\alpha-1}\int_{0}^{\frac{t}{v(t-1)}}\frac{t_{1}^{\alpha-1}}{1-t_{1}}dt_{1} \]
\ When $k_{h}<=\min_{u\leq t\leq v}k(t)$ define the function $\hat{F}_{\alpha,3}^{i}(k_{h},u,v)$ as
\[ \hat{F}_{\alpha,3}^{i}(k_{h},u,v)=\frac{(v+\frac{k_{h}}{\alpha v^{\alpha-1}})F_{\alpha}^{i-1}(k_{h},v)-\int_{u}^{v}\frac{1-\frac{k(a,t)}{t^{\alpha}}}{(t-1)}(t-1+\frac{k(t)}{R(t)\alpha (t-1)^{\alpha-1}})\breve{f}_{\alpha}^{i-1}(k(t)/R(t),t-1)dt }{u+\frac{k_{h}}{\alpha u^{\alpha-1}}} \]
\[ +\frac{\int_{u}^{v}(k(t)-k_{h})\frac{t-1}{t}F_{\alpha}^{(i-1)}(0,t-1)H(\alpha,v,t)dt}{u+\frac{k_{h}}{\alpha u^{\alpha-1}}} \]

\ We have

\begin{equation}
 S_{k_{l},\alpha}(\mathcal{A};\mathcal{P},z,\ln^{\alpha}(\xi^{2}))\leq Xe^{-\gamma}C(\omega)\ln^{\alpha-1}(\xi^{2})(u+\frac{k_{h}}{\alpha u^{\alpha-1}})\hat{F}_{\alpha,3}^{i}(k_{h},u,v)(1+O(\frac{1}{\ln^{\frac{1}{14}}(\xi^{2})}))
\end{equation} 

\[ +\ln^{\alpha}(\xi^{2})\sum_{d|\mathcal{P}(z),d<\xi^{2}}3^{v_{1}(d)}|r_{d}| \]

\begin{proof} Define the function
\[ S(v,k(u))= \sum_{\xi^{2/v} \leq p < \xi^{2/u}}(\ln^{\alpha}(\xi^{2})-k(a,\frac{\ln(\xi^{2})}{\ln(p)})\ln^{\alpha}(p))S_{\alpha,0}(\mathcal{A}_{p};\mathcal{P},p,) \]
\[ +\sum_{\xi^{2/v} \leq p < \xi^{2/u}}\sum_{q<p}k_{h}(a,\frac{\ln(\xi^{2})}{\ln(q)})\ln^{\alpha}(q)S_{\alpha,0}(\mathcal{A}_{pq};\mathcal{P},p) \]
\ Where
$$
 k_{h}(a,\frac{\ln(\xi^{2})}{\ln(q)})=\begin{cases}
\begin{array}{ll}
k(a,\frac{\ln(\xi^{2})}{\ln(q)}) & \frac{\ln(\xi^{2})}{\ln(q)}<v \\
k_{h} & \frac{\ln(\xi^{2})}{\ln(q)} \geq v \\
\end{array}
\end{cases}
$$
$$
 k_{h}(a,t)=\begin{cases}
\begin{array}{ll}
k(a,t) & t<v \\
k_{h} & t \geq v \\
\end{array}
\end{cases}
$$

\ For any $k_{h} \leq \min_{u\leq t \leq v}(k(t))$ we have 
\[ S(v,k_{h})= \sum_{\xi^{2/v} \leq p < \xi^{2/u}}(\ln^{\alpha}(\xi^{2})-k(a,\frac{\ln(\xi^{2})}{\ln(p)})\ln^{\alpha}(p))S_{\alpha,0}(\mathcal{A}_{p};\mathcal{P},p,) \]
\[ +\sum_{\xi^{2/v} \leq p < \xi^{2/u}}\sum_{q<p}k_{h}(a,\frac{\ln(\xi^{2})}{\ln(q)})\ln^{\alpha}(q)S_{\alpha,0}(\mathcal{A}_{pq};\mathcal{P},p) \]

\[ \geq \sum_{\xi^{2/v} \leq p < \xi^{2/u}}R(\frac{\ln(\xi^{2})}{\ln(p)})S_{\alpha,k(\frac{\ln(\xi^{2})}{\ln(p)})/R(\frac{\ln(\xi^{2})}{\ln(p)})}(\mathcal{A}_{p};\mathcal{P},p,\ln^{\alpha}(\xi^{2}/p)) \]
\[ -\sum_{\xi^{2/v} \leq p < \xi^{2/u}}(k(\frac{\ln(\xi^{2})}{\ln(p)})-k_{h})\sum_{q<\xi^{2/v}}\ln^{\alpha}(q)S_{\alpha,0}(\mathcal{A}_{pq};\mathcal{P},p) \]

\ Where
\[ R(\frac{\ln(\xi^{2})}{\ln(p)})=\frac{\ln^{\alpha}(\xi^{2})-k(a,\frac{\ln(\xi^{2})}{\ln(p)})\ln^{\alpha}(p)}{(\ln(\xi^{2})-\ln(p))^{\alpha}} =\frac{1-k(a,t)/t^{\alpha}}{(1-1/t)^{\alpha}}=R(t), \  \ t=\frac{\ln(\xi^{2})}{\ln(p)}\]

\ Infer the sum of the second sum on the right
\[ \sum_{q<\xi^{2/v}}\ln^{\alpha}(q)S_{\alpha,0}(\mathcal{A}_{pq};\mathcal{P},p) \leq \sum_{q<\xi^{2/v}}\ln^{\alpha}(q)\frac{XC(w)e^{-\gamma}\omega(pq)}{pq\ln(p)}F_{\alpha}^{(i-1)}(0,\frac{\ln(\xi^{2}/pq)}{\ln(p)})(1+o(1))+\sum_{q<\xi^{2/v}}\ln^{\alpha}(q)R_{pq} \]
\[ \leq \frac{XC(w)e^{-\gamma}\omega(p)}{p}\frac{\ln(\xi^{2}/p)}{\ln(p)}F_{\alpha}^{(i-1)}(0,\frac{\ln(\xi^{2}/p)}{\ln(p)})\sum_{1<\xi^{2/v}}\frac{\omega(q)\ln^{\alpha}(q)}{q\ln(\xi^{2}/p/q)}(1+o(1))+\ln^{\alpha}(p)R_{p} \]
\[ =\frac{XC(w)e^{-\gamma}\omega(p)}{p}\frac{\ln(\xi^{2}/p)}{\ln(p)}F_{\alpha}^{(i-1)}(0,\frac{\ln(\xi^{2}/p)}{\ln(p)})\ln^{\alpha-1}(\xi^{2})\int_{0}^{1/v}\frac{t_{1}^{\alpha-1}}{1-\frac{\ln(p)}{\ln(\xi^{2})}-t_{1}}dt_{1}(1+o(1))+\ln^{\alpha}(p)R_{p} \]

\ Where used
\[ \frac{1}{\ln(p)}F_{\alpha}^{(i-1)}(0,\frac{\ln(\xi^{2}/pq)}{\ln(p)})=\frac{\ln(\xi^{2}/pq)}{\ln(p)\ln(\xi^{2}/pq)}F_{\alpha}^{(i-1)}(0,\frac{\ln(\xi^{2}/pq)}{\ln(p)}) \leq \frac{\ln(\xi^{2}/p)F_{\alpha}^{(i-1)}(0,\frac{\ln(\xi^{2}/p)}{\ln(p)})}{\ln(p)\ln(\xi^{2}/pq)} \]
\ Thus
\[ \sum_{\xi^{2/v} \leq p < \xi^{2/u}}(k(\frac{\ln(\xi^{2})}{\ln(p)})-k_{h})\sum_{q<\xi^{2/v}}\ln^{\alpha}(q)S_{\alpha,0}(\mathcal{A}_{pq};\mathcal{P},p) \]
\[ \leq XC(w)e^{-\gamma}\ln^{\alpha-1}(\xi^{2})\int_{u}^{v}(k(t)-k_{h})\frac{t-1}{t}F_{\alpha}^{(i-1)}(0,t-1)dt\int_{0}^{1/v}\frac{t_{1}^{\alpha-1}}{1-1/t-t_{1}}dt_{1}(1+o(1))+\sum_{w\leq p<z}\ln^{\alpha}(p)R_{p} \]
\[ =XC(w)e^{-\gamma}\ln^{\alpha-1}(\xi^{2})\int_{u}^{v}(k(t)-k_{h})\frac{t-1}{t}F_{\alpha}^{(i-1)}(0,t-1)H(\alpha,v,t)dt(1+o(1))+\sum_{w\leq p<z}\ln^{\alpha}(p)R_{p} \]
\ Where used
\[ \int_{0}^{1/v}\frac{t_{1}^{\alpha-1}}{1-1/t-t_{1}}dt_{1}=(1-1/t)^{\alpha-1}\int_{0}^{\frac{t}{v(t-1)}}\frac{t_{1}^{\alpha-1}}{1-t_{1}}dt_{1}=H(\alpha,v,t) \]
\ Further
\[ S(v,k_{h}) \geq Xe^{-\gamma}C(\omega)\ln^{\alpha-1}(\xi^{2}) \int_{u}^{v}\frac{R(t)}{t}(1-\frac{1}{t})^{\alpha-1}(t-1+\frac{k(t)}{R(t)\alpha (t-1)^{\alpha-1}})\breve{f}_{\alpha}^{i-1}(k(t)/R(t),t-1)dt\times \]
\[ \left(1+O(\frac{1}{\ln^{\frac{1}{14}}(\xi^{2})})\right)-\ln^{\alpha}(\xi^{2})R \]
\[ -Xe^{-\gamma}C(\omega)\ln^{\alpha-1}(\xi^{2})\int_{u}^{v}(k(t)-k_{h})\frac{t-1}{t}F_{\alpha}^{(i-1)}(0,t-1)H(\alpha,v,t)dt \left(1+O(\frac{1}{\ln^{\frac{1}{14}}(\xi^{2})})\right)-\ln^{\alpha}(\xi^{2})R \]

\[ =Xe^{-\gamma}C(\omega)\ln^{\alpha-1}(\xi^{2})\int_{u}^{v}\frac{1-\frac{k(a,t)}{t^{\alpha}}}{t(1-\frac{1}{t})}(t-1+\frac{k(t)}{R(t)\alpha (t-1)^{\alpha-1}})\breve{f}_{\alpha}^{i-1}(k(t)/R(t),t-1)dt\times \]
\[ \left(1+O(\frac{1}{\ln^{\frac{1}{14}}(\xi^{2})})\right)-\ln^{\alpha}(\xi^{2})R \]
\[ -Xe^{-\gamma}C(\omega)\ln^{\alpha-1}(\xi^{2})\int_{u}^{v}(k(t)-k_{h})\frac{t-1}{t}F_{\alpha}^{(i-1)}(0,t-1)H(\alpha,v,t)dt \left(1+O(\frac{1}{\ln^{\frac{1}{14}}(\xi^{2})})\right)-\ln^{\alpha}(\xi^{2})R\]
\ Finally we obtain
\[ S_{\alpha,k_{h}}(\mathcal{A};\mathcal{P},z,\ln^{\alpha}(\xi^{2}))\leq \ln^{\alpha}(\xi^{2})S_{\alpha,0}(\mathcal{A};\mathcal{P},z)+k_{h}\sum_{p<\xi^{\frac{2}{v}}}\ln^{\alpha}(P)S_{\alpha,0}(\mathcal{A}_{p};\mathcal{P}(p),z)\]
\[ +\sum_{\xi^{\frac{2}{v}} \leq p <z}k(a,\frac{\ln(\xi^{2})}{\ln(p)})\ln^{\alpha}(p)S_{0}(\mathcal{A}_{p};\mathcal{P}(p),z) \]
\[ = S_{\alpha,k_{h}}(\mathcal{A};\mathcal{P},\xi^{\frac{2}{v}},\ln^{\alpha}(\xi^{2}))-S(v,k_{h}) \]

\[ \leq Xe^{-\gamma}C(\omega)\ln^{\alpha-1}(\xi^{2})(v+\frac{k_{h}}{\alpha v^{\alpha-1}})F_{\alpha}^{i-1}(k_{h},v)\left(1+O(\frac{1}{\ln^{\frac{1}{14}}(\xi^{2})})\right)+\ln^{\alpha}(\xi^{2})R\]
\[ -Xe^{-\gamma}C(\omega)\ln^{\alpha-1}(\xi^{2})\int_{u}^{v}\frac{1-\frac{k(a,t)}{t^{\alpha}}}{t(1-\frac{1}{t})}(t-1+\frac{k(t)}{R(t)\alpha (t-1)^{\alpha-1}})\breve{f}_{\alpha}^{i-1}(k(t)/R(t),t-1)dt \]
\[ \left(1+O(\frac{1}{\ln^{\frac{1}{14}}(\xi^{2})})\right)+\ln^{\alpha}(\xi^{2})R\]
\[ +Xe^{-\gamma}C(\omega)\ln^{\alpha-1}(\xi^{2})\int_{u}^{v}(k(t)-k_{h})\frac{t-1}{t}F_{\alpha}^{(i-1)}(0,t-1)H(\alpha,v,t)dt \left(1+O(\frac{1}{\ln^{\frac{1}{14}}(\xi^{2})})\right)+\ln^{\alpha}(\xi^{2})R\]
\[ =Xe^{-\gamma}C(\omega)\ln^{\alpha-1}(\xi^{2})(u+\frac{k_{h}}{\alpha u^{\alpha-1}})\hat{F}_{\alpha,3}^{i}(k_{h},u,v)\left(1+O(\frac{1}{\ln^{\frac{1}{14}}(\xi^{2})})\right)+\ln^{\alpha}(\xi^{2})R \]

\end{proof}
\ In this study will use two forms of the function $k(a,t)$ to do the interation. \\
\[
\begin{array}{ll}
 k(a,t)=\min_{t\leq x<v}(k_{n+1},\frac{(x-1)^{\alpha}}{a}) & 1\leq a\leq k_{n} \\
 k(a,t)=\frac{a}{(1-\frac{1}{t})^{\alpha}+\frac{a}{v^{\alpha}}} & 1\leq a\leq k_{n} \\
\end{array}
\]

\subsection{Lemma 1.22}
\ Suppose $0<k_{l}<k_{n}$
\[ \hat{F}_{\alpha,4}^{(i)}(k_{l},u)=\frac{\beta (u+\frac{k_{l_{1}}}{\alpha u^{\alpha-1}})F_{\alpha}^{(i-1)}(k_{l_{1}},u)+(1-\beta)(u+\frac{k_{l_{2}}}{\alpha u^{\alpha-1}})F_{\alpha}^{(i-1)}(k_{l_{2}},u)}{u+\frac{k_{l}}{\alpha u^{\alpha-1}}} \]

\ Where 
\[ k_{l_{1}}=\max_{0\leq h\leq l}k_{h}; \  \ k_{l_{1}}=\min_{k_{l} \leq h\leq l}k_{n} \]
\[ \beta=\frac{k_{l}-k_{l_{2}}}{k_{l_{1}}-k_{l_{2}}} \]
\ Is the solution of the equation
\[ \beta (u+\frac{k_{l_{1}}}{\alpha u^{\alpha-1}})+(1-\beta)(u+\frac{k_{l_{2}}}{\alpha u^{\alpha-1}})=u+\frac{k_{l}}{\alpha u^{\alpha-1}} \]
\ We have

\begin{equation}
 S_{k_{l},\alpha}(\mathcal{A};\mathcal{P},z,\ln^{\alpha}(\xi^{2}))\leq X\ln^{\alpha-1}(\xi^{2})(u+\frac{k_{l}}{\alpha u^{\alpha-1}})\hat{F}_{\alpha,4}^{(i)}(k_{l},u)(1+O(\frac{1}{\ln^{\frac{1}{14}}(\xi^{2})}))
\end{equation} 

\[ +\ln^{\alpha}(\xi^{2})\sum_{d|\mathcal{P}(z),d<\xi^{2}}3^{v_{1}(d)}|r_{d}| \]
\begin{proof}
\ Proof of this lemma is the same as the proof of $lemma 1.16$.
\end{proof}

\subsection{Lemma 1.23}
\ Suppose $0<k_{l}\leq k_{n}$, $0<u=\frac{\ln^(\xi^{2})}{\ln(z)}<k_{l}^{\frac{1}{\alpha}}$,
 $u_{1}=k_{l}^{\frac{1}{\alpha}}$\\
 
\[ \hat{F}_{\alpha,5}^{(i)}(k_{l},u)=\frac{1}{u+\frac{k_{l}}{u^{\alpha-1}}}\frac{k_{l}}{u^{\alpha}}(u_{1}+\frac{k_{l}}{u_{1}^{\alpha-1}})F_{\alpha}^{(i-1)}(k_{l},u_{1}) \]
\ We have
\[ S_{k_{l},\alpha}(\mathcal{A};\mathcal{P},z,\ln^{\alpha}(\xi^{2}))\leq XC(\omega)e^{-\gamma}(u+\frac{k_{l}}{u^{\alpha-1}})\hat{F}_{\alpha,5}^{i}(k_{l},u)(1+o(\frac{1}{\ln^{\frac{1}{14}}(\xi)})) +\ln^{\alpha}(\xi^{2})R\]

\begin{proof}
\ Suppose $z_{1}$ is the solution of the equation
\[ \frac{\ln(\xi^{2})}{\ln(z_{1})}= k_{l}^{\frac{1}{\alpha}} =u_{1}\]
\ When $u\leq u_{1}$
\[ S_{k_{l},\alpha}(\mathcal{A};\mathcal{P},z,\ln^{\alpha}(\xi^{2}))\leq S_{k_{l},\alpha}(\mathcal{A};\mathcal{P},z,\frac{k_{l}}{u^{\alpha}}\ln^{\alpha}(\xi^{2})) \]
\[ \leq S_{k_{l},\alpha}(\mathcal{A};\mathcal{P},z_{1},\frac{k_{l}}{u^{\alpha}}\ln^{\alpha}(\xi^{2})) \]
\[ \leq \frac{k_{l}}{u^{\alpha}}S_{k_{l},u}(\mathcal{A};\mathcal{P},z_{1},\ln^{\alpha}(\xi^{2}))\leq XC(\omega)e^{-\gamma}\frac{k_{l}}{u^{\alpha}}(u_{1}+\frac{k_{l}}{u_{1}^{\alpha-1}})F_{\alpha}^{(i-1)}(k_{l},u_{1})(1+o(\frac{1}{\ln^{\frac{1}{14}}(\xi)})) +\ln^{\alpha}(\xi^{2})R\]
\[ =XC(\omega)e^{-\gamma}(u+\frac{k_{l}}{u^{\alpha-1}})\hat{F}_{\alpha,5}^{(i)}(k_{l},u)(1+o(\frac{1}{\ln^{\frac{1}{14}}(\xi)})) +\ln^{\alpha}(\xi^{2})R\]

\end{proof}
\ In this study, the functin $\hat{F}_{\alpha,5}^{(i)}(k_{l},u)$ in region $u<k_{l}^{\frac{1}{\alpha}}$ is not always used to perform the iteration.

\subsection{Theorem 1}
\ Combining these lemmas, we obtain:
\ for any $i-1 \geq 0$ the equations \eqref{StartUp} and \eqref{StartDown} are correct, define the functions 
\[ F_{\alpha}^{(i)}(k_{l},u)=\min_{k_{l}^{\frac{1}{\alpha}}\leq u<v}\left(F_{\alpha}^{(i-1)}(k_{l},u),F_{\alpha,1}^{(i)}(k_{l},u,v),\hat{F}_{\alpha,2}^{(i)}(k_{l},u),\hat{F}_{\alpha,3}^{(i)}(k_{l},u.v),\hat{F}_{\alpha,4}^{(i)}(k_{l},u)\right), \ \ u\geq k_{l}^{\frac{1}{\alpha}} \]
\[ F_{\alpha}^{(i)}(k_{l},u)=\hat{F}_{\alpha,5}^{(i)}(k_{l},u), \ \ 0<u<k_{l}^{\frac{1}{\alpha}} \]
\ And

\[ f_{\alpha}^{(i)}(k_{l},u)=\max_{0<u<v} \left(f_{\alpha}^{(i-1)}(k_{l},u),\hat{f}_{\alpha,1}^{(i)}(k_{l},u,v),\hat{f}_{\alpha,2}^{(i)}(k_{l},u,v),\hat{f}_{\alpha,3}^{(i)}(k_{l},u),\hat{f}_{\alpha,4}^{(i)}(k_{l},u)\right),  \ \ u>0 \]

\ We have:
\ When , $0\leq k_{l}\leq k_{n}$

\begin{equation}
 S_{k_{l},\alpha}(\mathcal{A};\mathcal{P},z,\ln^{\alpha}(\xi^{2}))\leq Xe^{-\gamma}C(\omega)\ln^{\alpha-1}(\xi^{2})(u+\frac{k_{l}}{\alpha u^{\alpha-1}})F_{\alpha}^{(i)}(k_{l},u)(1+O(\frac{1}{\ln^{\frac{1}{14}}(\xi^{2})}))
\end{equation} 

\[ +\ln^{\alpha}(\xi^{2})\sum_{d|\mathcal{P}(z),d<\xi^{2}}3^{v_{1}(d)}|r_{d}| \]
\ And
\ When , $0\leq k_{l}\leq k_{n+4}$

\begin{equation}
 S_{k_{l},\alpha}(\mathcal{A};\mathcal{P},z,\ln^{\alpha}(\xi^{2}))\geq Xe^{-\gamma}C(\omega)\ln^{\alpha-1}(\xi^{2})(u+\frac{k_{l}}{\alpha u^{\alpha-1}})f_{\alpha}^{(i)}(k_{l},u)(1+O(\frac{1}{\ln^{\frac{1}{14}}(\xi^{2})}))
\end{equation} 

\[ -\ln^{\alpha}(\xi^{2})\sum_{d|\mathcal{P}(z),d<\xi^{2}}3^{v_{1}(d)}|r_{d}| \]

\ These two functions  are used in series to create an iterative program, using flowing calculate order. \\


\ in the beginning set $\alpha=2$, and
\begin{itemize}
\item each $k_{l}$ is peformed for $4$ cycles \\
\item $k_{l}$ from $f_{n}$ to $0$, is peformed for $8$ cycles. \\
\end{itemize}
\ Parameter $v$ for each iteration are: \\
\\
$$
\begin{array}{ll}
\ Lemma 1.14 & v=10 \\
\ Lemma 1.15 & v=3,3.5,4,4.5,5 \\
\ Lemma 1.18 & v=4,4.5,5 \\
\ Lemma 1.19 & v=10 \\
\ Lemma 1.21 & v=3,2.75,2.5,2.25 \\
\end{array}
$$
\\

\ Table $1$ is the values of $e^{-\gamma}(u+\frac{k_{l}}{2u})F_{2}(k_{l},u)$ and $e^{-\gamma}(u+\frac{k_{l}}{2u})f_{2}(k_{l},u)$ at some points:
\begin{table}[h]
	\centering
	\caption{$e^{-\gamma}(u+\frac{k_{l}}{2u})F_{2}(k_{l},u)$ and $e^{-\gamma}(u+\frac{k_{l}}{2u})f_{2}(k_{l},u)$}
		\begin{tabular}{|c|c|c|c|c|c|c|c|c|c|c|} \hline
			 u & \multicolumn{2}{|c|}{\textbf{3}} & \multicolumn{2}{|c|}{\textbf{2.5}} & \multicolumn{2}{|c|}{\textbf{2}} & \multicolumn{2}{|c|}{\textbf{1.5}} & \multicolumn{2}{|c|}{\textbf{1}} \\ \hline
			 $k_{l}$ & F & f & F & f & F & f & F & f & F & f \\ \hline
			4.5 &   & 1.96635 &  & 1.67752 &  & 1.28813 &  & 0.70587 &  & 0 \\ \hline
			4 & 2.15054 & 1.95519 & 2.02792 & 1.66557 & 1.99927 & 1.27638 & 3.48815 & 0.63944 & 7.74209 & 0 \\ \hline
			3.75 & 2.12684 & 1.94961 & 2.00342 & 1.65960 & 1.97867 & 1.27052 & 3.25499 & 0.60623 & 7.24307 & 0 \\ \hline
			3.5 & 2.10295 & 1.91658 & 1.97861 & 1.61773 & 1.95809 & 1.22884 & 3.02183 & 0.56581 & 6.67440 & 0 \\ \hline
			3.25 & 2.08865 & 1.88156 & 1.96646 & 1.56967 & 1.94637 & 1.18084 & 2.79741 & 0.52540 & 6.25374 & 0 \\ \hline
			3 & 2.07414 & 1.84654 & 1.95397 & 1.51728 & 1.93465 & 1.12785 & 2.57543 & 0.48498 & 5.76573 & 0 \\ \hline
			2.75 & 2.06162 & 1.81152 & 1.94537 & 1.46425 & 1.92659 & 1.06885 & 2.35346 & 0.44457 & 5.27772 & 0 \\ \hline
			2.5 & 2.04885 & 1.77650 & 1.93640 & 1.41122 & 1.91853 & 1.00295 & 2.13149 & 0.40415 & 4.78970 & 0 \\ \hline
			2.25 & 2.03817 & 1.74147 & 1.93415 & 1.35820 & 1.91807 & 0.92819 & 1.91837 & 0.36374 & 4.31075 & 0 \\ \hline
			2 & 2.02748 & 1.70645 & 1.93190 & 1.30516 & 1.91761 & 0.84215 & 1.91761 & 0.32332 & 3.83178 & 0 \\ \hline
			1.75 & 2.01680 & 1.67143 & 1.92965 & 1.25213 & 1.91714 & 0.75020 & 1.91714 & 0.28290 & 3.35282 & 0 \\ \hline
			1.5 & 2.00611 & 1.63641 & 1.92740 & 1.19910 & 1.91668 & 0.65825 & 1.91668 & 0.24249 & 2.87389 & 0 \\ \hline
			1.25 & 1.99543 & 1.60138 & 1.92515 & 1.14607 & 1.91622 & 0.56629 & 1.91622 & 0.20207 & 2.39495 & 0 \\ \hline
			1 & 1.98475 & 1.56637 & 1.92290 & 1.09305 & 1.91576 & 0.47434 & 1.91576 & 0.16166 & 1.91602 & 0 \\ \hline
			0.75 & 1.97406 & 1.53134 & 1.92065 & 1.04002 & 1.91529 & 0.38198 & 1.91529 & 0.12125 & 1.91529 & 0 \\ \hline
			0.5 & 1.96338 & 1.49632 & 1.91840 & 0.98699 & 1.91483 & 0.28962 & 1.91483 & 0.08083 & 1.91483 & 0 \\ \hline
			0.25 & 1.95269 & 1.46130 & 1.91615 & 0.93396 & 1.91437 & 0.19726 & 1.91437 & 0.04041 & 1.91437 & 0 \\ \hline
			0 & 1.94201 & 1.42628 & 1.91390 & 0.88809 & 1.91390 & 0.10490 & 1.91390 & 0 & 1.91390 & 0 \\ \hline

		\end{tabular}
\end{table} \\

\ In order to improve these results, we set $\alpha=3.5$, and consider  \\

\begin{equation}
S_{0,\alpha}(\mathcal{A};\mathcal{P},z,\ln^{\alpha}(\xi^{2}))=\ln^{\alpha-2}(\xi^{2})S_{0,2}(\mathcal{A};\mathcal{P},z,\ln^{2}(\xi^{2}))=\ln^{\alpha}(\xi^{2})S_{0}(\mathcal{A};\mathcal{P},z)
\end{equation}

\begin{equation}
S_{k_{n},\alpha}(\mathcal{A};\mathcal{P},z,\ln^{\alpha}(\xi^{2}))= \ln^{\alpha}(\xi^{2})S_{0}(\mathcal{A};\mathcal{P},z)+\sum_{p<z}k_{n}\ln^{\alpha}(p)S_{0}(\mathcal{A}_{p};\mathcal{P},z)
\end{equation}

\[ =\ln^{\alpha}(\xi^{2})S_{0}(\mathcal{A};\mathcal{P},z) +\ln^{\alpha-2}(\xi^{2})\sum_{p<z}k_{n}\frac{\ln^{\alpha-2}(p)}{\ln^{\alpha-2}(\xi^{2})}\ln^{2}(p)S_{0}(\mathcal{A}_{p};\mathcal{P},z) \]
\[ \leq \ln^{\alpha-2}(\xi^{2})\left(S_{k_{n}\frac{\ln^{\alpha-2}(z)}{\ln^{\alpha-2}(\xi^{2})},\alpha}(\mathcal{A};\mathcal{P},z,\ln^{2}(\xi^{2}))\right) \]
\ So that we can set\\

\begin{itemize}
\item $uF_{3.5}^{(0)}(0,u)=uF_{2}^{(i)}(0,u)$ \\
\item $uf_{3.5}^{(0)}(0,u)=uf_{2}^{(i)}(0,u)$ \\
\end{itemize}

\ and \\

\[ (u+\frac{k_{n}}{3.5 u^{3.5-1}})F_{3.5}^{(0)}(k_{n},u)=(u+\frac{k_{n}}{u^{3.5-2}}\frac{1}{2u})F_{2}^{(i)}(\frac{k_{n}}{u^{3.5-2}},u), \  \ u\geq 2 \]

\ The iteration is repeated for 4 cycles for each $k_{l}$, and 4 cycles from $k_{l}=0$ to $k_{n}$. The results for $k_{l}=0$ are kept and reset $\alpha=2$, set\\

\begin{itemize}
\item $F_{2}^{(0)}(0,u)=F_{3.5}^{(i)}(0,u)$ \\
\item $f_{2}^{(0)}(0,u)=f_{3.5}^{(i)}(0,u)$ \\
\end{itemize}

\ Then set $\alpha=4$ and repeat the iteration, finally again set $\alpha=2$, repeat the iteration. \\
\ Table $2$ is the resulting values of $e^{-\gamma}(u+\frac{k_{l}}{2u})F_{2}(k_{l},u)$ and $e^{-\gamma}(u+\frac{k_{l}}{2u})f_{2}(k_{l},u)$ at some points:
 \\
\ Table $3$ is the resulting values of $e^{-\gamma}uF_{2}(0,u)$ and $e^{-\gamma}uf_{2}(0,u)$, $(k_{l}=0)$ listed from 1.8 to 5 in $0.1$ increments \\

\begin{table}[h]
	\centering
	\caption{$e^{-\gamma}(u+\frac{k_{l}}{2u})F_{2}(k_{l},u)$ and $e^{-\gamma}(u+\frac{k_{l}}{2u})f_{2}(k_{l},u)$}
		\begin{tabular}{|c|c|c|c|c|c|c|c|c|c|c|} \hline
			 u & \multicolumn{2}{|c|}{\textbf{3}} & \multicolumn{2}{|c|}{\textbf{2.5}} & \multicolumn{2}{|c|}{\textbf{2}} & \multicolumn{2}{|c|}{\textbf{1.5}} & \multicolumn{2}{|c|}{\textbf{1}} \\ \hline
			 $k_{l}$ & F & f & F & f & F & f & F & f & F & f \\ \hline
			4.5 &   & 1.97453 &  & 1.68857 &  & 1.30463 &  & 0.72792 &  & 0 \\ \hline
			4 & 2.14267 & 1.96384 & 2.01741 & 1.67727 & 1.99419 & 1.29362 & 3.45385 & 0.66592 & 7.62429 & 0 \\ \hline
			3.75 & 2.11819 & 1.95849 & 1.99171 & 1.67161 & 1.96864 & 1.28812 & 3.22027 & 0.63491 & 7.13005 & 0 \\ \hline
			3.5 & 2.093361 & 1.92817 & 1.96583 & 1.63301 & 1.943359 & 1.24989 & 2.98668 & 0.59687 & 6.63582 & 0 \\ \hline
			3.25 & 2.07708 & 1.89786 & 1.95059 & 1.59028 & 1.92889 & 1.20735 & 2.76343 & 0.55452 & 6.15192 & 0 \\ \hline
			3 & 2.06046 & 1.86754 & 1.93519 & 1.54486 & 1.91420 & 1.15984 & 2.54245 & 0.51186 & 5.67028 & 0 \\ \hline
			2.75 & 2.04534 & 1.83723 & 1.92323 & 1.49944 & 1.90294 & 1.10631 & 2.32147 & 0.46921 & 5.18865 & 0 \\ \hline
			2.5 & 2.03013 & 1.80691 & 1.91113 & 1.45402 & 1.89168 & 1.04593 & 2.10049 & 0.42655 & 4.70702 & 0 \\ \hline
			2.25 & 2.01723 & 1.77660 & 1.90779 & 1.40859 & 1.89018 & 0.97690 & 1.89044 & 0.38390 & 4.23630 & 0 \\ \hline
			2 & 2.00432 & 1.74628 & 1.90445 & 1.36317 & 1.88868 & 0.89681 & 1.88868 & 0.34124 & 3.76557 & 0 \\ \hline
			1.75 & 1.99142 & 1.71596 & 1.90111 & 1.31775 & 1.88718 & 0.81618 & 1.8718 & 0.29859 & 3.29484 & 0 \\ \hline
			1.5 & 1.97851 & 1.68565 & 1.89777 & 1.27233 & 1.88569 & 0.73558 & 1.88569 & 0.25593 & 2.82414 & 0 \\ \hline
			1.25 & 1.96851 & 1.65533 & 1.89443 & 1.22690 & 1.88419 & 0.665493 & 1.88419 & 0.21327 & 2.35343 & 0 \\ \hline
			1 & 1.95270 & 1.62502 & 1.89108 & 1.18148 & 1.88269 & 0.57430 & 1.88269 & 0.17062 & 1.88273 & 0 \\ \hline
			0.75 & 1.93980 & 1.59470 & 1.88774 & 1.13606 & 1.88119 & 0.49044 & 1.88119 & 0.12796 & 1.88119 & 0 \\ \hline
			0.5 & 1.92689 & 1.56439 & 1.88440 & 1.09064 & 1.87969 & 0.40658 & 1.87969 & 0.08531 & 1.87969 & 0 \\ \hline
			0.25 & 1.91399 & 1.53407 & 1.88106 & 1.04522 & 1.87820 & 0.32272 & 1.87820 & 0.04265 & 1.87820 & 0 \\ \hline
			0 & 1.90109 & 1.50375 & 1.87772 & 0.99979 & 1.87670 & 0.23886 & 1.87670 & 0 & 1.87670 & 0 \\ \hline
			
		\end{tabular}

\end{table}

\begin{table}[h]
	\centering
	\caption{$e^{-\gamma}uF_{2}(0,u)$ and $e^{-\gamma}uf_{2}(0,u)$}
		\begin{tabular}{||c|c|c||c|c|c||c|c|c||} \hline
			u & $e^{-\gamma}$uF(u) & $e^{-\gamma}$uf(u) & u & $e^{-\gamma}$uF(u) & $e^{-\gamma}$uf(u) & u & $e^{-\gamma}$uF(u) & $e^{-\gamma}$uf(u) \\ \hline
			5.0 & 2.810476 & 2.804123 & 3.9 & 2.227293 & 2.153511 & 2.8 & 1.893647 & 1.306029 \\ \hline
			4.9 & 2.755139 & 2.747114 & 3.8 & 2.179677 & 2.087060 & 2.7 & 1.892139 & 1.226451 \\ \hline
			4.8 & 2.700062 & 2.689884 & 3.7 & 2.133011 & 2.022424 & 2.6 & 1.887881 & 1.112676 \\ \hline
			4.7 & 2.645264 & 2.632382 & 3.6 & 2.088863 & 1.951076 & 2.5 & 1.877724 & 0.999797  \\ \hline
			4.6 & 2.590828 & 2.574554 & 3.5 & 2.046887 & 1.885336 & 2.4 & 1.877724 & 0.870318 \\ \hline
			4.5 & 2.536905 & 2.516300 & 3.4 & 2.008704 & 1.808683 & 2.3 & 1.877175 & 0.731240 \\ \hline
			4.4 & 2.483362 & 2.457531 & 3.3 & 1.974608 & 1.728772 & 2.2 & 1.876697 & 0.581023 \\ \hline
			4.3 & 2.430558 & 2.398088 & 3.2 & 1.945059 & 1.655096 & 2.1 & 1.876697 & 0.417728 \\ \hline
			4.2 & 2.378490 & 2.337796 & 3.1 & 1.921803 & 1.567792 & 2.0 & 1.876697 & 0.238863 \\ \hline
			4.1 & 2.327326 & 2.276432 & 3.0 & 1.901086 & 1.503759 & 1.9 & 1.876697 & 4.1132E-02  \\ \hline
			4.0 & 2.276645 & 2.217810 & 2.9 & 1.893859 & 1.407497 & 1.8 & 1.876697 & 0 \\ \hline
		\end{tabular}
	\label{tab:UF20UAndUf20U}
\end{table}

\section{Part II: Double Sieve, the Goldbach and the twin primes problems} 

\ Corresponding Goldbach's problem, suppose $N$ is an even integer number, define 
\[ \mathcal{A}:=\left\{N-p,p<N\right\} \]
\[ X=\left|\mathcal{A}\right|=\frac{N}{\ln(N)} +O\left(\frac{N}{\ln^{2}(N)} \right)\]
\[ C(N)=\prod_{p|N,P>2}\frac{p-1}{p-2}\prod_{p>2}\left(1-\frac{1}{(p-1)^{2}}\right) \]

\ Or corresponding twin primes problem define
\[ \mathcal{A}:=\left\{p-2,p<N\right\} \]
\[ X=\left|\mathcal{A}\right|=\frac{N}{\ln(N)} +O\left(\frac{N}{\ln^{2}(N)} \right)\]
\[ C(2)=\prod_{p>2}\left(1-\frac{1}{(p-1)^{2}}\right) \]
\ $\Omega(n)$ is the number of all prime factors of the integer $n$.
\[ D_{1,2}(N) := |\left\{n \in \mathcal{A} :\Omega(n)\leq 2 \right\}|.\]
\[ D(N) := |\left\{n \in \mathcal{A} :\Omega(n)= 1 \right\}|.\]

\ Chen \cite{J. R. Chen1973} \cite{J.R. Chen1978} \cite{C.D. Pan and C.B. Pan1992} using his double sieve system proved:
\[ D_{1,2}(N) > 0.67\frac{C(N)}{\ln^{2}(N)} \]
\ And
\[ D(N)\leq 7.8342\frac{C(N)}{\ln^{2}(N)} \]
\ Chen's lower bound of $D_{1,2}(N)$ and upper bound of $D(N)$ has been improved by many authors.
\ In this paper, similar double sieve is used to get better results of the lower and upper bound. This paper will prove:
\subsection{Theorem 3}

\begin{equation}
 D_{1,2}(N) \geq 2.27\frac{C(N)}{\ln^{2}(N)}(1+o(1))
\end{equation}

\ And
\subsection{Theorem 4}

\begin{equation}
 D(N)\leq 6.916\frac{C(N)}{\ln^{2}(N)}(1+o(1))
\end{equation}

\ In this part only suppose $\alpha=2$, and need two parameters $\xi^{2}$ and $\xi_{1}^{2}$, in the beginning of the iteration, we set functions $F^{(0)}(k_{l},u)$ and $f^{(0)}(k_{l},u)$ as the results functions of Part I. \\
\  We need the following two lemmas to estimate the remainder term.
\subsection{Lemma 2.1 (theorem of E. Bombieri\cite{E. Bombieri1965})}
\ Suppose $x >2$,$D=\frac{X^{\frac{1}{2}}}{\ln^{B}(x)}$, For any number $A>0$, $B=A+15$, We have
\[ R(D,x)=\sum_{d<D}\max_{y\leq x}\max_{(l,d)=1}\left|\pi(y,d,l)-\frac{\pi(y)}{\phi(d)}\right| \]
\[ =\sum_{d<D}\max_{y\leq x}\max_{(l,d)=1}\left|\sum_{p<y,p\equiv l(d)}1-\frac{1}{\phi(d)}\sum_{p<y}1\right| \ll x\ln^{-A}(x) \]

\subsection{Corollary 2.1.1}
\ Suppose $x>2$, $D=\frac{X^{\frac{1}{2}}}{\ln^{B}(x)}$, For any number $A>0$, $B=2A+32$, We have
\[ \sum_{d<D}\max_{y\leq x}\max_{(l,d)=1}\mu^{2}(d)3^{v_{1}(d)}\left|\pi(y,d,l)-\frac{\pi(y)}{\phi(d)}\right| \ll x\ln^{-A}(x) \]

\subsection{Lemma 2.2}(A new mean value theorem\cite{C.D. Pan1981}\cite{C.D. Pan and C.B. Pan1992})
\ Suppose $0<\beta \leq 1$, function $E(x)$ and $g_{x}(a)$
 
\begin{equation}
\frac{1}{2}\leq E(x) \ll x^{1-\beta},
\end{equation}

\begin{equation}
g_{x}(a) \ll d^{r}(a),
\end{equation}

\ For any number $A>0$, $B=\frac{3}{2}A+2^{2r+2}+13$,$D=x^{\frac{1}{2}}\ln^{-B}(x)$

\begin{equation}
 R(D,x,\mathcal{E}_{x})=\sum_{d\leq D}\max_{y\leq x}\max_{(l,d)=1}\left|\sum_{a\leq E(x)}g_{x}(a)E(y,a,d,l)\right| 
\end{equation}

\[ =\sum_{d\leq D}\max_{y\leq x}\max_{(l,d)=1}\left|\sum_{a\leq E(x)}g_{x}(a)\left(\sum_{ap<y,ap\equiv l(d)}1-\frac{1}{\phi(d)}\pi(\frac{y}{a})\right)\right| \ll x\ln^{-A}(x) \]
\subsection{Corollary 2.2.1}
\ For any number $A>0$, $B=3A+7\cdot2^{r}+232$, $D=x^{\frac{1}{2}}\ln^{-B}(x)$

\begin{equation}
 R_{1}(D,x,\mathcal{E}_{x})=\sum_{d\leq D}\mu^{2}(d)3^{v_{1}(d)}\max_{y\leq x}\max_{(l,d)=1}\left|\sum_{a\leq E(x)}g_{x}(a)E(y,a,d,l)\right| \ll x\ln^{-A}(x)
\end{equation}

\subsection{Lemma 2.3}
\ Continuum function\cite{W.B. Jurkat and H.-E. Richert1965} $w(u)$ settle for condition
$$
\begin{cases}
\begin{array}{ll}
w(u)=\frac{1}{u} & 1\leq u\leq 2 \\
(uw(u))'=w(u-1) & u>2 \\
\end{array}
\end{cases}
$$
\ We have
 
\begin{equation}
w(u)\leq 1 \ \ ,u\geq 1 
\end{equation}

\subsection{Lemma 2.4}
\ Set 
\[ \mathcal{N}:=\left\{n, n\leq x ,n \ \ is \ \ natural \ \ number\right\} \]
\[ \mathcal{N}_{d}:=\left\{n, n\in \mathcal{N}, d|n ,\right\} \]
\ Suppose $x>1$, $\frac{\ln(x)}{\ln(z)}=u\geq 1$, we have
 
\begin{equation}
S_{0}(\mathcal{N},\mathcal{P},z)=\sum_{1\leq n\leq x,(n,\mathcal{P}_{1}(z))=1}1=w(u)\frac{x}{\ln(z)}+O(\frac{x}{\ln^{2}(z)})+O(\frac{z}{\ln(z)}) 
\end{equation}

\ Where
\[ \mathcal{P}_{1}(z)=\prod_{p<z}p \]
\subsection{Lemma 2.5}
\ Suppose $d\geq 1$, $a(d)>0$
\[ \sum_{w\leq p<z}a(p)S_{k,2}(\mathcal{A}_{pd};\mathcal{P}(d),p,\ln^{2}(\frac{\xi^{2}}{p})) \]
\[ \leq \frac{4F_{2}^{(0)}(0,2)}{\ln(N)}\frac{N}{d}2C(N)e^{-\gamma}\sum_{w\leq p<z,(p,d)=1}a(p)\frac{\ln(\xi^{2}/p)}{p}(t_{p}-1+\frac{k}{2(t_{p}-1)})(1+0(1))\]
\ Where 
\[ t_{p}=\frac{\ln(\xi^{2})}{\ln(p)} \]

\begin{proof}
\ Set $\xi^{2}_{1}=X^{\frac{1}{2}}/\ln^{B}(X)$
\[ \sum_{w\leq p<z,(p,d)=1}a(p)S_{0}(\mathcal{A}_{Pd};\mathcal{P}(d),p)=\sum_{w\leq p<z,(p,d)=1}a(p)\sum_{n \in \mathcal{A},pd|n,(n/d,\mathcal{P}(p))=1}1 \]
\[ =\sum_{w\leq p<z,(p,d)=1}a(p)\sum_{n<N,pd|n,(n/d,\mathcal{P}(p))=1 \atop (N-n,N^{\frac{1}{2}})=1}1 \]
\[ =\sum_{w\leq p<z,(p,d)=1}a(p)\sum_{N-n<N,pd|(N-n),((N-n)/d,\mathcal{P}(p))=1 \atop (n,N^{\frac{1}{2}})=1}1 \]
\[ \leq 2C(N)e^{-\gamma}\frac{1}{\ln(\xi^{2}_{1})}\frac{\ln(\xi^{2}_{1})}{\ln(N^{\frac{1}{2}})}F_{2}(0,\frac{\ln(\xi^{2}_{1})}{\ln(N^{\frac{1}{2}})})\sum_{w\leq p<z}a(p)\sum_{N-n<N,pd|(N-n),((N-n)/d,\mathcal{P}_{1}(p))=1}1\times \]
\[ (1+o(1))+R^{(1)} \]
\[ =X_{d}2C(N)e^{-\gamma}\frac{2}{\ln(N)}2F_{2}(0,2)(1+o(1))+R^{(1)} \]
\ Where
\[ R^{(1)}=\sum_{d|\mathcal{P},d<\xi^{2}_{1}}a(d)3^{v_{1}(d)}\left|r_{d}\right| \]
\ Set $g_{x}(a)=a(d)3^{v_{1}(d)}\ll d^{2}(a)$
\ According to $Lemma 2.2$, $Corollary 2.2.1$, when $\xi^{2}_{1}\leq \frac{N^{\frac{1}{2}}}{\ln^{B}(N)}$, $B=3A+260$ 
\[ R^{(1)} \ll \frac{N}{\ln^{A}(N)} \]
\ On the other hand, by $Lemma 2.4$
 
\begin{equation}
X_{d}=\sum_{w\leq p<z}a(p)\sum_{N-n<N,pd|(N-n),((N-n)/d,\mathcal{P}_{1}(p))=1}1=\sum_{w\leq p<z}a(p)\sum_{n<N,pd|(n),(n/d,\mathcal{P}_{1}(p))=1}1 
\end{equation}

\[ =\frac{N}{d}\sum_{w\leq p<z,(p,d)=1}\frac{a(p)}{p\ln(p)}w(\frac{\ln(N/dp)}{\ln(p)})(1+o(1)) \leq \frac{N}{d}\sum_{w\leq p<z,(p,d)=1}\frac{a(p)}{p\ln(p)}(1+o(1))\]
\ Where used $Lemma 2.3$ $w(u)\leq 1$. So we have
 
\begin{equation}
\sum_{w\leq p<z,(p,d)=1}a(p)S_{0}(\mathcal{A}_{Pd};\mathcal{P}(d),p) \leq \frac{N}{d}\frac{8C(N)e^{-\gamma}}{\ln(N)}F_{2}(0,2)\sum_{w\leq p<z,(p,d)=1}\frac{a(p)}{p\ln(p)}(1+o(1))
\end{equation}

\ Finally we obtain
\[ \sum_{w\leq p<z}a(p)S_{k,2}(\mathcal{A}_{pd};\mathcal{P},p,\ln^{2}(\frac{\xi^{2}}{p}))\]
\[ =\sum_{w\leq p<z}a(p)\ln^{2}(\xi^{2}/p)S_{0}(\mathcal{A}_{Pd};\mathcal{P},p)+\sum_{w\leq p<z}a(p)\sum_{q<p}k\ln^{2}(q)S_{0}(\mathcal{A}_{Pqd};\mathcal{P}(q),p) \]
\[ \leq \frac{N}{d}\frac{8C(N)e^{-\gamma}}{\ln(N)}F_{2}(0,2)\sum_{w\leq p<z}\frac{a(p)\ln(\xi^{2}/p)}{p\ln(p)}(1+o(1)) \]
\[ +\frac{N}{d}\frac{8C(N)e^{-\gamma}}{\ln(N)}F_{2}(0,2)\sum_{w\leq p<z}\frac{a(p)}{p\ln(p)}\sum_{q<p}k\ln^{2}(q)(1+o(1)) \]
\[ =\frac{N}{d}\frac{8C(N)e^{-\gamma}}{\ln(N)}F_{2}(0,2)\sum_{w\leq p<z}\frac{a(p)}{p\ln(p)}\left(\ln^{2}(\xi^{2}/p)+\frac{k}{2}\ln^{2}(p)\right)(1+o(1)) \]
\[ =\frac{N}{d}\frac{2C(N)e^{-\gamma}}{\ln(N)}4F_{2}(0,2)\sum_{w\leq p<z}\frac{a(p)\ln(\xi^{2}/p)}{p}\left(t_{p}-1+\frac{k}{2(t_{p}-1)}\right)(1+o(1)) \]

\end{proof}

\ The following Lemmas (2.6 to 2.9) need to be added to the iteration.
\subsection{Lemma 2.6}
\ Suppose $t^{2}>k_{l}$ 
\[ k_{l}(t)=k_{l}\frac{(1-\frac{1}{t})^{2}}{1-\frac{k_{l}}{t^2}}  \]
If $0<k_{l}(t)\leq k_{n}$ define the function $\breve{F}_{2}^{(i-1)}(k_{l}(t),t-1)$ as\\

\[ \breve{F}_{2}^{(i-1)}(k_{l}(t),t-1)=\beta (t-1+\frac{k_{l_{1}}}{2(t-1)})F_{2}^{(i-1)}(k_{l_{1}},t-1)+(1-\beta)(t-1+\frac{k_{l_{2}}}{2(t-1)})F_{2}^{(i-1)}(k_{l_{2}},t-1) \]
\ Where
\[ k_{l_{1}}=\min_{0<m\leq n \atop k_{m}\geq k_{l}(t)}(k_{m}),\ \ k_{l_{2}}=\max_{0\leq m<n \atop k_{m}\leq k_{l}(t)}(k_{m}) \]

\[ \beta= \frac{k_{l}(t)-k_{l_{2}}}{k_{l_{1}}-k_{l_{2}}}\]
\ Is a solution of equation
\[ \beta k_{l_{1}}+(1-\beta)k_{l_{2}}=k_{l}(t) \]
If $k_{t}> k_{n}$ define the functions $\breve{F}_{2,h}^{(i-1)}(k_{l}(t),t-1)$ as
\[ \breve{F}_{2,1}^{(i-1)}(k_{l}(t),t-1)=(t-1+\frac{k_{n}}{2(t-1)})F_{2}^{(i-1)}(k_{n},t-1) \]
\[ \breve{F}_{2,2}^{(i-1)}(k_{l}(t),t-1)=(t-1+\frac{k_{l}(t)}{2(t-1)})4F_{2}^{(0)}(0,2) \]
\[ \breve{F}_{2}^{(i-1)}(k_{l}(t),t-1)=\min\left(\breve{F}_{2,1}^{(i-1)}(k_{l}(t),t-1),\breve{F}_{2,2}^{(i-1)}(k_{l}(t),t-1)\right) \]
\ And
\[ \breve{F}_{2,3}^{(i-1)}(k_{l}(t),t-1)=\frac{(k_{l}(t)-k_{n})}{2t^{2}}(t-1)4F_{2}^{(0)}(2), \ \ k_{l}(t)>k_{n}\]
\[ \breve{F}_{2,3}^{(i-1)}(k_{l}(t),t-1)=0, \ \ k_{l}(t)<k_{n} \]

\[ \hat{f}_{2,5}^{(i)}(k_{l},u,v)=\frac{1}{(u+\frac{k_{l}}{2u})}(v+\frac{k_{l}}{2v})f_{2}^{i-1}(k_{l},v) \]
\[ -\frac{1}{(u+\frac{k_{l}}{2u})}\int_{u}^{v}\frac{1-\frac{k_{l}}{t^2}}{t(1-\frac{1}{t})}\left((\breve{F}_{2}^{(i-1)}(k_{l}(t),t-1)+\breve{F}_{2,3}^{(i-1)}(k_{l}(t),t-1)\right)dt \]

\ We have \\

\begin{equation}
 S_{k_{l},2}(\mathcal{A};\mathcal{P},z,\ln^{2}(\xi^{2})) \geq Xe^{-\gamma}\ln(\xi^{2})\max_{v>u}(u+\frac{k_{l}}{2u})\hat{f}_{2,5}^{(i)}(k_{l},u,v)(1+o(1)) -\ln^{2}(\xi^{2})R
\end{equation}

\begin{proof}
\ Suppose $\ln^{4+c_{1}}(x)\leq w,z$, $v=\frac{\ln(\xi^{2})}{\ln(w)}$, $u=\frac{\ln(\xi^{2})}{\ln(z)}$, $1\leq u<u_{1}\leq u$, $k_{l}\leq u^{2}$
\[ S_{k_{l},2}(\mathcal{A};\mathcal{P},z,\ln^{2}(\xi^{2}))=S_{k_{l},2}(\mathcal{A};\mathcal{P},w,\ln^{2}(\xi^{2}))
-\sum_{w\leq p<z}S_{k_{l},2}(\mathcal{A}_{p};\mathcal{P},p,\ln^{2}(\xi^{2})-k_{l}\ln^{2}(p))+O\left(\frac{X}{\ln^{2}(X)}\right)\]
\[ =S_{k_{l},2}(\mathcal{A};\mathcal{P},w,\ln^{2}(\xi^{2}))-\Sigma_{1}+O\left(\frac{X}{\ln^{2}(X)}\right)\]

\ Where
\[ \Sigma_{1}=\sum_{w\leq p<z}\frac{\ln^{2}(\xi^{2})-k_{l}\ln^{2}(p)}{\ln^{2}(\xi^{2}/p)}S_{k_{l}(t_{p}),2}(\mathcal{A}_{p};\mathcal{P},p,\ln^{2}(\xi^{2}/p)) \]
\[ =\sum_{w\leq p<z, k_{l}(t_{p})\leq k_{n}}\frac{\ln^{2}(\xi^{2})-k_{l}\ln^{2}(p)}{\ln^{2}(\xi^{2}/p)}S_{k_{l}(t_{p}),2}(\mathcal{A}_{p};\mathcal{P},p,\ln^{2}(\xi^{2}/p))\]
\[ +\sum_{w\leq p<z, k_{l}(t_{p})>k_{n}}\frac{\ln^{2}(\xi^{2})-k_{l}\ln^{2}(p)}{\ln^{2}(\xi^{2}/p)}S_{k_{n},2}(\mathcal{A}_{p};\mathcal{P},p,\ln^{2}(\xi^{2}/p))\]
\[ +\sum_{w\leq p<z, k_{l}(t_{p})>k_{n}}\frac{\ln^{2}(\xi^{2})-k_{l}\ln^{2}(p)}{\ln^{2}(\xi^{2}/p)}(k_{l}(t_{p})-k_{n})\sum_{q<p}\ln^{2}(q)S_{0,2}(\mathcal{A}_{pq};\mathcal{P}(q),p) \]
\[ =\Sigma_{11}+\Sigma_{12}+\Sigma_{13} \]

\ Where $t_{p}=\frac{\ln(\xi^{2})}{\ln(p)}$
\ Since that
\[ S_{k_{l}(t_{p}),2}(\mathcal{A}_{p};\mathcal{P},p,\ln^{2}(\xi^{2}/p))\leq \frac{Xe^{-\gamma}2C(N)}{p}\ln(\xi^{2}/p)(t_{p}-1+\frac{k_{l}(t_{p})}{2(t_{p}-1)})F_{2}^{(i-1)}(k_{l}(t_{p}),t_{p}-1)\left(1+o(1)\right)\]
\[ +\ln^{2}(p)R_{p} \]
\[ =\frac{Xe^{-\gamma}2C(N)}{p}\ln(\xi^{2}/p)\breve{F}_{2}^{(i-1)}(k_{l}(t_{p}),t_{p}-1)\left(1+o(1)\right)+\ln^{2}(p)R_{p} \]
\ So that
\[ \Sigma_{11}\leq \frac{N2C(N)e^{-\gamma}}{\ln(N)}\sum_{w\leq p<z,k_{l}(t_{p})\leq k_{n}}\frac{\ln^{2}(\xi^{2})-k_{l}\ln^{2}(p)}{p\ln^{2}(\xi^{2}/p)}\ln(\xi^{2}/p)(t_{p}-1+\frac{k_{l}(t_{p})}{2(t_{p}-1)})\breve{F}_{2}^{(i-1)}(k_{l}(t_{p}),t_{p}-1)\times \]
\[ (1+o(1)) +\sum_{w\leq p<z,k_{l}(t_{p})\leq k_{n}}\ln^{2}(p)R_{p} \]
\[ =\frac{N2C(N)e^{-\gamma}}{\ln(N)}\ln(\xi^{2})\int_{u_{1}}^{v}\frac{1-\frac{k_{l}}{t^{2}}}{t(1-\frac{1}{t})}(t-1+\frac{k_{l}(t)}{2(t-1)})\breve{F}_{2}^{(i-1)}(k_{l}(t),t-1)dt(1+o(1))+\ln^{2}(\xi^{2})R \]
\[ =\frac{N2C(N)e^{-\gamma}}{\ln(N)}\ln(\xi^{2})\int_{u_{1}}^{v}\frac{1-\frac{1}{t^{2}}}{t(1-\frac{1}{t})}\breve{F}_{2,1}^{(i-1)}(k_{l},t-1)dt(1+o(1))+\ln^{2}(\xi^{2})R \]

\ and
\[ S_{k_{n},2}(\mathcal{A}_{p};\mathcal{P},p,\ln^{2}(\xi^{2}/p))\leq \frac{Xe^{-\gamma}2C(N)}{p}\ln(\xi^{2}/p)(t_{p}-1+\frac{k_{n}}{2(t_{p}-1)})F_{2}^{(i-1)}(k_{n},t_{p}-1)\left(1+o(1)\right)+\ln^{2}(p)R_{p} \]
\[ =\frac{Xe^{-\gamma}2C(N)}{p}\ln(\xi^{2}/p)\breve{F}_{2,1}^{(i-1)}(k_{n},t_{p}-1)\left(1+o(1)\right)+\ln^{2}(p)R_{p} \]
\[ \sum_{w\leq p<z,k_{l}(t_{p})> k_{n}}\frac{\ln^{2}(\xi^{2})-k_{l}\ln^{2}(p)}{\ln^{2}(\xi^{2}/p)}S_{k_{n},2}(\mathcal{A}_{p};\mathcal{P},p,\ln^{2}(\xi^{2}/p)) \]
\[ \leq \frac{N2C(N)e^{-\gamma}}{\ln(N)}\sum_{w\leq p<z,k_{l}(t_{p})> k_{n}}\frac{\ln^{2}(\xi^{2})-k_{l}\ln^{2}(p)}{p\ln^{2}(\xi^{2}/p)}\ln(\xi^{2}/p)(t_{p}-1+\frac{k_{n}}{2(t_{p}-1)})\breve{F}_{2}^{(i-1)}(k_{n},t_{p}-1)(1+o(1))\]
\[ +\sum_{w\leq p<z,k_{l}(t_{p})> k_{n}}\ln^{2}(p)R_{p} \]
\[ =\frac{N2C(N)e^{-\gamma}}{\ln(N)}\ln(\xi^{2})\int_{u}^{u_{1}}\frac{1-\frac{1}{t^{2}}}{t(1-\frac{1}{t})}(t-1+\frac{k_{n}}{2(t-1)})\breve{F}_{2}^{(i-1)}(k_{n},t-1)dt(1+o(1))+\ln^{2}(\xi^{2})R \]
\[ =\frac{N2C(N)e^{-\gamma}}{\ln(N)}\ln(\xi^{2})\int_{u}^{u_{1}}\frac{1-\frac{1}{t^{2}}}{t(1-\frac{1}{t})}\breve{F}_{2,1}^{(i-1)}(k_{n},t-1)dt(1+o(1))+\ln^{2}(\xi^{2})R \]

\ On the other hand, similar of Chen's double sieve according to $Lemma 2.5$ set $u_{1}$ as the solution of
\[ k_{l}(u_{1})=k_{n} \]
\ We have
\[ \Sigma_{11}=\sum_{w\leq p<z,k_{l}(t_{p})\leq k_{n}}\frac{\ln^{2}(\xi^{2})-(k_{l}\ln^{2}(p))}{\ln^{2}(\xi^{2}/p)}S_{k_{l}(t_{p}),2}(\mathcal{A}_{p};\mathcal{P},p,\ln^{2}(\xi^{2}/p))\]
\[ \leq N2C(N)e^{-\gamma}\frac{4F_{2}(0,2)}{\ln(N)}\sum_{w\leq p<z,k_{l}(t_{p})\leq k_{n}}\frac{\ln^{2}(\xi^{2})-(k_{l}\ln^{2}(p))}{p\ln(\xi^{2}/p)}\left(t_{p}-1+\frac{k_{p}}{2(t_{p}-1)}\right)\left(1+o(1)\right) \]
\[ =X2C(N)e^{-\gamma}\frac{4F_{2}(0,2)}{\ln(N)}\ln(\xi^{2})\int_{u_{1}}^{v}\frac{1-\frac{k_{l}(t)}{t^{2}}}{t(1-\frac{1}{t})}(t-1+\frac{k_{l}(t)}{2(t-1)})dt\left(1+o(1)\right) \]
\[ =\frac{N2C(N)e^{-\gamma}}{\ln(N)}\ln(\xi^{2})\int_{u_{1}}^{v}\frac{1-\frac{k_{l}(t)}{t^{2}}}{t(1-\frac{1}{t})}\breve{F}_{2,2}^{(i-1)}(k_{l}(t),t-1)dt\left(1+o(1)\right) \]
\ and
\[ \Sigma_{12}=\sum_{w\leq p<z,k_{l}(t_{p})> k_{n}}\frac{\ln^{2}(\xi^{2})-(k_{l}\ln^{2}(p))}{\ln^{2}(\xi^{2}/p)}S_{k_{n},2}(\mathcal{A}_{p};\mathcal{P},p,\ln^{2}(\xi^{2}/p))\]
\[ \leq N2C(N)e^{-\gamma}\frac{4F_{2}(0,2)}{\ln(N)}\sum_{w\leq p<z,k_{l}(t_{p})> k_{n}}\frac{\ln^{2}(\xi^{2})-(k_{n}\ln^{2}(p))}{p\ln(\xi^{2}/p)}\left(t_{p}-1+\frac{k_{n}}{2(t_{p}-1)}\right)\left(1+o(1)\right) \]

\[ =\frac{X2C(N)e^{-\gamma}}{\ln(N)}4F_{2}(0,2)\ln(\xi^{2})\int_{u}^{u_{1}}\frac{1-\frac{k_{n}}{t^{2}}}{t(1-\frac{1}{t})}(t-1+\frac{k_{n}(t)}{2(t-1)})dt\left(1+o(1)\right) \]

\[ =\frac{N2C(N)e^{-\gamma}}{\ln(N)}\ln(\xi^{2})\int_{u}^{u_{1}}\frac{1-\frac{k_{l}(t)}{t^{2}}}{t(1-\frac{1}{t})}\breve{F}_{2,2}^{(i-1)}(k_{n},t-1)dt\left(1+o(1)\right) \]

\ and 
\[ \Sigma_{13}=\sum_{w\leq p<z, k_{l}(t_{p})>k_{n}}\frac{\ln^{2}(\xi^{2})-k_{l}\ln^{2}(p)}{\ln^{2}(\xi^{2}/p)}(k_{l}(t_{p})-k_{n})\sum_{q<p}\ln^{2}(q)S_{0,2}(\mathcal{A}_{pq};\mathcal{P}(q),p) \]
\[ \leq N2C(N)e^{-\gamma}\frac{4F_{2}(0,2)}{\ln(N)}\sum_{q<p}\frac{\ln^{2}(q)}{q\ln^{2}(\xi^{2})}\times \]
\[ \sum_{w\leq p<z,k_{l}(t_{p})> k_{n}}\frac{\ln^{2}(\xi^{2})-(k_{n}\ln^{2}(p))}{p\ln(\xi^{2}/p)}\left(k_{l}(t_{p})-k_{n}\right))\left(t_{p}-1\right)\left(1+o(1)\right) \]
\[ =N2C(N)e^{-\gamma}\frac{4F_{2}(0,2)}{\ln(N)}\sum_{w\leq p<z,k_{l}(t_{p})> k_{n}}\frac{\ln^{2}(\xi^{2})-(k_{n}\ln^{2}(p))}{p\ln(\xi^{2}/p)}\left(k_{l}(t_{p})-k_{n}\right)\frac{\ln^{2}(p)}{2\ln^{2}(\xi^{2})}\left(t_{p}-1\right)\left(1+o(1)\right) \]
\[ =X2C(N)e^{-\gamma}\frac{4F_{2}(0,2)}{\ln(N)}\ln(\xi^{2})\int_{u}^{u_{1}}\frac{1-\frac{k_{n}}{t^{2}}}{t(1-\frac{1}{t})}\frac{k_{l}(t)-k_{n}}{2t^{2}}(t-1)dt\left(1+o(1)\right) \]
\[ =\frac{X2C(N)e^{-\gamma}}{\ln(N)}\ln(\xi^{2})\int_{u}^{u_{1}}\frac{1-\frac{k_{n}}{t^{2}}}{t(1-\frac{1}{t})}\breve{F}_{2,3}^{(i-1)}(k_{l},t-1)dt\left(1+o(1)\right) \]

\ Finally we obtain
\[ \Sigma_{11}+\Sigma_{12}+\Sigma_{12}\leq \frac{N2C(N)e^{-\gamma}}{\ln(N)}\ln(\xi^{2})\int_{u}^{v}\frac{1-\frac{k_{n}}{t^{2}}}{t(1-\frac{1}{t})}\left(\breve{F}_{2}^{(i-1)}(k_{l},t-1)+\breve{F}_{2,3}^{(i-1)}(k_{l},t-1)\right)dt\left(1+o(1)\right) \]
\[ S_{k_{l},2}(\mathcal{A};\mathcal{P},z,\ln^{2}(\xi^{2}))\geq \frac{N2C(N)e^{-\gamma}}{\ln(N)}\ln(\xi^{2})(v+\frac{k_{l}}{2v})f_{2}^{(i-1)}(k_{l},v) \]
\[ -\frac{N2C(N)e^{-\gamma}}{\ln(N)}\ln(\xi^{2})\int_{u}^{v}\frac{1-\frac{k_{n}}{t^{2}}}{t(1-\frac{1}{t})}\left(\breve{F}_{2}^{(i-1)}(k_{l},t-1)+\breve{F}_{2,3}^{(i-1)}(k_{l},t-1)\right)dt\left(1+o(1)\right) -\ln^{2}(\xi^{2})R\]
\[ =\frac{N2C(N)e^{-\gamma}}{\ln(N)}\ln(\xi^{2})(u+\frac{k_{l}}{2u})\hat{f}_{2,5}^{(i)}(k_{l},u)-\ln^{2}(\xi^{2})R\]
 
\end{proof}

\subsection{Lemma 2.7}
\ Suppose $l>n$, $k_{l}=\frac{v^{2}}{2}$, $\frac{v}{2} \geq u_{1} \geq u$, $v\geq 3$ \\

\ When $t \geq u_{1}$
\[ \breve{F}_{2,1}^{(i-1)}(0,t)=\frac{1}{2t}vF_{2}^{(i-1)}(0,v-\frac{v}{t}) \]
\ When $u_{1}> t \geq u$
\[ \breve{F}_{2,1}^{(i-1)}(0,t)=\left(\frac{1}{2}\ln(\frac{v}{u_{1}})+\ln(\frac{u_{1}}{t})\right)4F_{2}^{(0)}(0,2) \]

\[ \hat{f}_{2,6}^{(i)}(k_{l},u,v)=\frac{1}{(u+\frac{k_{l}}{2u})}vf_{2}^{(i-1)}(0,v) \]
\[ -\frac{1}{(u+\frac{k_{l}}{2u})}\int_{u}^{v}\breve{F}_{2,1}^{(i-1)}(0,t)dt \]
\ We have
 
\begin{equation}
 S_{k_{l},2}(\mathcal{A};\mathcal{P},z,\ln^{2}(\xi^{2})) \geq Xe^{-\gamma}\ln(\xi^{2})(u+\frac{k_{l}}{2u})\max_{v>u,v\geq 3}\hat{f}_{2,6}^{(i)}(k_{l},u,v)
\end{equation}

\begin{proof}
\ Suppose $v=3,4,4.5,5,5.5$, $v=\frac{\ln(\xi^{2})}{\ln(w)}$, $\ln^{4+c_{1}}(\xi^{2})<w^{\sqrt{2}}\leq z_{1}<z$, According to the definition of $k_{l}$ in $Lemma 1.12$, when $k_{l}=\frac{v^{2}}{2}$ we have \\

\[ \frac{1}{2}\ln^{2}(\xi^{2})\leq k_{l}\ln^{2}(p), \ \ p>w \]
\[ \ln^{2}(\xi^{2})\leq k_{l}\ln^{2}(p), \ \ p>w^{\sqrt{2}} \]
\ Thus
\[ S_{k_{l},2}(\mathcal{A};\mathcal{P},z,\ln^{2}(\xi^{2}))\geq \ln^{2}(\xi^{2})S_{0}(\mathcal{A};\mathcal{P},z) \]
\[ +\frac{1}{2}\ln^{2}(\xi^{2})\sum_{w\leq q<z_{1}}S_{0}(\mathcal{A}_{p};\mathcal{P}(q),z) +\ln^{2}(\xi^{2})\sum_{z_{1}\leq q<z}S_{0}(\mathcal{A}_{p};\mathcal{P}(q),z)\]
\[ =\ln^{2}(\xi^{2})S_{0}(\mathcal{A};\mathcal{P},w)-\frac{\ln^{2}(\xi^{2})}{2}\sum_{w\leq p<z_{1}}\left(S_{0}(\mathcal{A}_{p};\mathcal{P},p)+\sum_{w\leq q<p}S_{0}(\mathcal{A}_{pq};\mathcal{P}(q),p)\right) \]
\[ -\frac{\ln^{2}(\xi^{2})}{2}\sum_{w\leq q<p<z, q<z_{1}<p}S_{0}(\mathcal{A}_{pq};\mathcal{P}(q),p)-\ln^{2}(\xi^{2})\sum_{z_{1}\leq q<p<z}S_{0}(\mathcal{A}_{pq};\mathcal{P}(q),p) _O\left(\frac{X}{\ln^{2}(X)}\right)\]
\[ \geq \ln^{2}(\xi^{2})S_{0}(\mathcal{A};\mathcal{P},w)-\frac{\ln^{2}(\xi^{2})}{2}\sum_{w\leq p<z_{1}}S_{0}(\mathcal{A}_{p};\mathcal{P}(p),w) \]
\[ -\frac{\ln^{2}(\xi^{2})}{2}\sum_{w\leq q<p<z, q<z_{1}<p}S_{0}(\mathcal{A}_{pq};\mathcal{P}(q),p)-\ln^{2}(\xi^{2})\sum_{z_{1}\leq q<p<z}S_{0}(\mathcal{A}_{pq};\mathcal{P}(q),p) +O\left(\frac{X}{\ln^{2}(X)}\right)\]
\[ =\Sigma_{1}-\frac{1}{2}\Sigma_{2}-\frac{1}{2}\Sigma_{3}-\Sigma_{4} +O\left(\frac{X}{\ln^{2}(X)}\right)\]
\ It is easy to see that
\[ \Sigma_{1} \geq X2C(N)e^{-\gamma}\ln(\xi^{2})vf_{2}^{(i-1)}(0,v)(1+o(1))-\ln^{2}(\xi^{2})R \]
\ and
\[ \Sigma_{2} \leq X2C(N)e^{-\gamma}\ln^{2}(\xi^{2})\sum_{w\leq p<z_{1}}\frac{\omega(p)}{p}\frac{F_{2}^{(i-1)}(0,v-v\frac{\ln(p)}{\ln(\xi^{2}/p)})}{\ln(w)}(1+o(1))+\ln^{2}(\xi^{2})R \]
\[ =X2C(N)e^{-\gamma}\ln(\xi^{2})v\int_{u_{1}}^{v}\frac{F_{2}^{(i-1)}(0,v-\frac{v}{t})}{t}dt(1+o(1))+\ln^{2}(\xi^{2})R \]
\[ =X2C(N)e^{-\gamma}\ln(\xi^{2})2\int_{u_{1}}^{v}\breve{F}_{2,1}^{(i-1)}(0,t)dt (1+0(1))+\ln^{2}(\xi^{2})R\]
\ Accordinh to equation (3.10) of $Lemma 2.5$
\[ \sum_{z_{1}\leq p<z}S_{0}(\mathcal{A}_{pq};\mathcal{P}(q),p)\leq \frac{N2C(N)}{q\ln(N)}4e^{-\gamma}F_{2}^{(0)}(0,2)\sum_{z_{1}\leq p<z}\frac{1}{p\ln(p)}(1+o(1)) \]

\ Thus

\[ \Sigma_{3}\leq \frac{N2C(N)}{\ln(N)}4e^{-\gamma}F_{2}^{(0)}(0,2)\ln^{2}(\xi^{2})\sum_{z_{1}\leq p<z}\frac{1}{p\ln(p)}\sum_{w\leq q<z_{1}}\frac{1}{q}(1+o(1))\]
\[ =\frac{N2C(N)}{\ln(N)}4e^{-\gamma}F_{2}^{(0)}(0,2)\ln^{2}(\xi^{2})\sum_{z_{1}\leq p<z}\frac{1}{p\ln(p)}\ln(\frac{\ln(w)}{\ln(z_{1})})(1+o(1))\]
\[ =\frac{N2C(N)}{\ln(N)}4e^{-\gamma}F_{2}^{(0)}(0,2)\ln(\frac{v}{u_{1}})\ln(\xi^{2})\int_{u_{1}}^{v}dt(1+o(1))\]
\ Where
\[ \int_{u_{1}}^{v}dt=v-u_{1}\leq u_{1}-u=\int_{u}^{u_{1}}dt \]
\ and
\[ \Sigma_{4}\leq \frac{N2C(N)}{\ln(N)}4e^{-\gamma}F_{2}^{(0)}(0,2)\ln^{2}(\xi^{2})\sum_{z_{1}\leq p<z}\frac{1}{p\ln(p)}\sum_{z_{1}\leq q<p}\frac{1}{q}(1+o(1))+\ln^{2}(\xi^{2})R\]
\[ =\frac{N2C(N)}{\ln(N)}4e^{-\gamma}F_{2}^{(0)}(0,2)\ln^{2}(\xi^{2})\sum_{z_{1}\leq p<z}\frac{1}{p\ln(p)}\ln(\frac{\ln(z_{1})}{\ln(p)})(1+o(1))+\ln^{2}(\xi^{2})R\]
\[ =\frac{N2C(N)}{\ln(N)}4e^{-\gamma}F_{2}^{(0)}(0,2)\ln^{2}(\xi^{2})\int_{z_{1}}^{z}\frac{1}{t\ln^{2}(t)}\ln(\frac{\ln(z_{1})}{\ln(t)})dt(1+o(1))+\ln^{2}(\xi^{2})R\]
\[ =\frac{N2C(N)}{\ln(N)}4e^{-\gamma}F_{2}^{(0)}(0,2)\ln(\xi^{2})\int_{u}^{u_{1}}\ln(\frac{u_{1}}{t})dt(1+o(1))+\ln^{2}(\xi^{2})R\]

\ Thus
\[ \frac{1}{2}\Sigma_{3}+\Sigma_{4} \leq \frac{N2C(N)}{\ln(N)}4e^{-\gamma}F_{2}^{(0)}(0,2)\ln(\xi^{2})\int_{u}^{u_{1}}\left(\frac{1}{2}\ln(\frac{v}{u_{1}})+\ln(\frac{u_{1}}{t})\right)dt(1+o(1))+\ln^{2}(\xi^{2})R\]
\[ =\frac{N2C(N)}{\ln(N)}4e^{-\gamma}\ln(\xi^{2})\int_{u}^{u_{1}}\breve{F}_{2,1}^{(i-1)}(0,t)dt(1+o(1))+\ln^{2}(\xi^{2})R\]

\[ S_{k_{l},2}(\mathcal{A};\mathcal{P},z,\ln^{2}(\xi^{2}))\geq \frac{N2C(N)e^{-\gamma}}{\ln(N)}\ln(\xi^{2})\left(vf_{2}^{(i-1)}(0,v)-\int_{u}^{v}\breve{F}_{2,1}^{(i-1)}(0,t)dt\right)(1+0(1)) -\ln^{2}(\xi^{2})R\]
\[ =\frac{N2C(N)e^{-\gamma}}{\ln(N)}\ln(\xi^{2})(u+\frac{k_{l}}{2u})\hat{f}_{2,6}^{(i)}(k_{l},u,v)(1+0(1))-\ln^{2}(\xi^{2})R \]

\ Combining these to get lemma 2.7 .
\end{proof}
\subsection{Lemma 2.8}
\ Suppose $v>u$
\[ \breve{F}_{2,1}^{(i-1)}(0,t)=\frac{1}{2t}vF_{2}^{(i-1)}(0,v-\frac{v}{t}) \]
\[ \breve{F}_{2,2}^{(i-1)}(0,t)=\frac{1}{2}\ln(\frac{v}{t})4F_{2}^{(0)}(0,2) \]
\[ \hat{f}_{2,7}^{(i)}(0,u,v)=\frac{1}{u}vf_{2}^{i-1}(0,v) \]
\[ -\frac{1}{u}\int_{u}^{v}\breve{F}_{2,1}^{(i-1)}(0,t)dt-\frac{1}{u}\breve{F}_{2,2}^{(i-1)}(0,u) \]
\ We have
 
\begin{equation}
 S_{0,2}(\mathcal{A};\mathcal{P},z,\ln^{2}(\xi^{2})) \geq Xu\ln(\xi^{2})\max_{v>u,v\geq 3}\hat{f}_{2,7}^{(i)}(0,u,v)
\end{equation}

\begin{proof}In $Lemma1.11$ we take k as $0$ 
\[ S_{0,2}(\mathcal{A};\mathcal{P},z,\ln^{2}(\xi^{2}))=\ln^{2}(\xi^{2})S_{0}(\mathcal{A};\mathcal{P},z) \]
\[ \geq \ln^{2}(\xi^{2})S_{0}(\mathcal{A};\mathcal{P},w)-\frac{\ln^{2}(\xi^{2})}{2}\sum_{w\leq p<z}S_{0}(\mathcal{A}_{p};\mathcal{P},w)-\frac{\ln^{2}(\xi^{2})}{2}\sum_{w\leq p<z}S_{0}(\mathcal{A}_{p};\mathcal{P}(p),z) +O\left(\frac{X}{\ln^{2}(X)}\right)\]
\[ =\Sigma_{1}-\frac{1}{2} \Sigma_{2}-\frac{1}{2} \Sigma_{3} +O\left(\frac{X}{\ln^{2}(X)}\right)\]
\ Hence
\[ \Sigma_{1} \geq X2C(N)e^{-\gamma}\ln(\xi^{2})vf_{2}^{(i-1)}(0,v)(1+0(1))-\ln^{2}(\xi^{2})R \]
\ and
\[ \Sigma_{2} \leq X2C(N)e^{-\gamma}\ln(\xi^{2}) \int_{u}^{v}v\frac{F_{2}^{(i-1)}(0,v-v/t)}{t}dt(1+0(1))+\ln^{2}(\xi^{2})R \]
\[ =X2C(N)e^{-\gamma}\ln(\xi^{2})\int_{u}^{v}2\breve{F}_{2,1}^{(i-1)}(0,t)dt(1+0(1))+\ln^{2}(\xi^{2})R \]
\ Similar of $Lemma 2.5$ we have
\[ \Sigma_{3}=\sum_{w\leq p<z}\ln^{2}(\xi^{2})S_{0}(\mathcal{A}_{p};\mathcal{P}(p),z) \leq X2C(N)4e^{-\gamma}F_{2}^{(0)}(0,2)\frac{\ln^{2}(\xi^{2})}{\ln(z)}\sum_{w\leq p<z}\frac{1}{p}(1+o(1)) \]

\[ =X2C(N)e^{-\gamma}4F_{2}^{(0)}(0,2)\ln(\xi^{2})u\ln(\frac{\ln(z)}{\ln(w)})(1+o(1)) \]

\[ =X2C(N)e^{-\gamma}4F_{2}^{(0)}(0,2)\ln(\xi^{2})u\ln(\frac{v}{u})(1+o(1)) \]
\[ =X2C(N)e^{-\gamma}\ln(\xi^{2})2u\breve{F}_{2,2}^{(i-1)}(0,u)(1+o(1)) \]
\ So we obtain
\[ S_{0,2}(\mathcal{A};\mathcal{P},z,\ln^{2}(\xi^{2}))\geq  X2C(N)e^{-\gamma}\ln(\xi^{2})\left(vf_{2}^{(i-1)}(0,v)-\frac{1}{2} \int_{u}^{v}\breve{F}_{2,1}^{(i-1)}(0,t)dt-u\breve{F}_{2,2}^{(i-1)}(0,u)\right)\times \]
\[ (1+o(1))-\ln^{2}(\xi^{2})R \]
\[ =X2C(N)e^{-\gamma}\ln(\xi^{2})u\hat{f}_{2,7}^{(i)}(0,u,v)(1+o(1))-\ln^{2}(\xi^{2})R \]

\end{proof} 

\subsection{Lemma 2.9}
\ Suppose 
\[ u_{1}=\frac{1}{1-\frac{F_{2}^{(i-1)}(0,u)}{4F_{2}^{(0)}(0,2)}} \]
\ Where $1/u_{1}$ is the solution of
\[ F_{2}^{(i-1)}(0,u)\frac{k_{l}t}{1-t}=4F_{2}^{(0)}(0,2)k_{l}t \]
\ when $u>0$, $u_{1}>u$, $0<k_{l}\leq k_{n}$
\[ \breve{F}_{2,1}^{(i-1)}(0,u)=F_{2}^{(i-1)}(0,u)\int_{0}^{1/u_{1}}\frac{k_{l}t}{1-t}dt, \ \ t<u_{1}\]
\[ \breve{F}_{2,2}^{(i-1)}(0,u)=4F_{2}^{(0)}(0,2)\int_{1/u_{1}}^{1/u}k_{l}tdt, \ \ t\geq u_{1}\]

\[ \hat{f}_{2,8}^{(i)}(0,u)=\frac{1}{u}\left((u+\frac{k_{l}}{2u})f_{2}^{(i-1)}(k_{l},u)-u\breve{F}_{2,1}^{(i-1)}(0,u)-u\breve{F}_{2,2}^{(i-1)}(0,u)\right) \]
\ We have
\[ S_{0,2}(\mathcal{A},\mathcal{P},z,\ln^{2}(\xi^{2}))\geq X2C(N)e^{-\gamma}\ln(\xi^{2})u\hat{f}_{2,8}^{(i)}(0,u)(1+o(1))-\ln^{2}(\xi^{2})R\]

\begin{proof} By equation (2.6)
\[ S_{0,2}(\mathcal{A},\mathcal{P},z,\ln^{2}(\xi^{2}))=S_{k_{l},2}(\mathcal{A},\mathcal{P},z,\ln^{2}(\xi^{2}))\]
\[ -\sum_{q<z}k_{l}\ln^{2}(q)S_{0}(\mathcal{A}_{q},\mathcal{P}(q),z) \]
\ The sum on the right
\[ \sum_{q<z}k_{l}\ln^{2}(q)S_{0}(\mathcal{A}_{q},\mathcal{P}(q),z)=\sum_{q<z_{1}}k_{l}\ln^{2}(q)S_{0}(\mathcal{A}_{q},\mathcal{P}(q),z)\]
\[ +\sum_{z_{1}\leq q<z}k_{l}\ln^{2}(q)S_{0}(\mathcal{A}_{q},\mathcal{P}(q),z)=\Sigma_{1}+\Sigma_{2} \]
\ Where
\[ \Sigma_{1}\leq X2C(N)e^{-\gamma}\frac{\ln^{2}(\xi^{2})}{\ln(z)}\sum_{q<z_{1}}\frac{k_{l}\omega(q)\ln^{2}(q)}{q}F_{2}^{(i-1)}(0,\frac{\ln{\xi^{2}/q}}{\ln(z)})(1+o(1)) +\ln^{2}(\xi^{2})R\]
\[ = X2C(N)e^{-\gamma}\ln(\xi^{2})u\int_{0}^{1/u_{1}}F_{2}^{(i-1)}(0,u-ut)k_{l}tdt(1+o(1))+\ln^{2}(\xi^{2})R \]
\[ \leq X2C(N)e^{-\gamma}\ln(\xi^{2})uF_{2}^{(i-1)}(0,u)\int_{0}^{1/u_{1}}\frac{k_{l}t}{1-t}dt(1+o(1))+\ln^{2}(\xi^{2})R \]
\[ =X2C(N)e^{-\gamma}\ln(\xi^{2})u\breve{F}_{2,1}^{(i-1)}(0,u)(1+o(1))+\ln^{2}(\xi^{2})R \]
\ Similar as $Lemma 2.5$
\[ \Sigma_{2} \leq X2C(N)e^{-\gamma}4F_{2}^{(0)}(0,2)\frac{\ln^{2}(\xi^{2})}{\ln(z)}\sum_{z_{1}\leq q<z}\frac{k_{l}\ln^{2}(q)}{q}(1+o(1)) \]
\[ =X2C(N)e^{-\gamma}4F_{2}^{(0)}(0,2)\ln(\xi^{2})u\int_{1/u_{1}}^{1/u}k_{l}tdt(1+o(1)) \]
\[ =X2C(N)e^{-\gamma}\ln(\xi^{2})u\breve{F}_{2,2}^{(i-1)}(0,u)(1+o(1)) \]
\ So we obtain
\[ S_{0,2}(\mathcal{A},\mathcal{P},z,\ln^{2}(\xi^{2}))\geq X2C(N)e^{-\gamma}\ln{\xi^{2}}(u+\frac{k_{l}}{2u})f_{2}^{(i-1)}(k_{l},u) \]
\[ -X2C(N)e^{-\gamma}\ln(\xi^{2})u\left(\breve{F}_{2,1}^{(i-1)}(0,u)-\breve{F}_{2,2}^{(i-1)}(0,u)\right)(1+o(1))\]
\[ =X2C(N)e^{-\gamma}\ln(\xi^{2})u\hat{f}_{2,8}^{(i)}(0,u,v)(1+o(1))-\ln^{2}(\xi^{2})R \]
\end{proof}
\ Same as $Theorem 1$,do interation with $4$ cycles of each $k_{l}$; 4 cycles from $k_{l}=0$ to $k_{l}=k_{n}$. We obtain $Theorem 2$

\subsection{Theorem 2}
\ If for any $i \geq 1$, and $\alpha=2$ define the functions 
\[ F_{\alpha}^{(i)}(k_{l},u)=\min_{k_{l}^{\frac{1}{\alpha}}\leq u<v}\left(F_{\alpha}^{(i-1)}(k_{l},u),F_{\alpha,1}^{(i)}(k_{l},u,v),\hat{F}_{\alpha,2}^{(i)}(k_{l},u),\hat{F}_{\alpha,3}^{(i)}(k_{l},u.v),\hat{F}_{\alpha,4}^{(i)}(k_{l},u)\right), \ \ u\geq k_{l}^{\frac{1}{\alpha}} \]
\[ F_{\alpha}^{(i)}(k_{l},u)=\hat{F}_{\alpha,5}^{(i)}(k_{l},u), \ \ 0<u<k_{l}^{\frac{1}{\alpha}} \]
\ And
\[ f_{\alpha}^{(i)}(k_{l},u)=\max\left(\max_{j=1,2,v>u}\hat{f}_{\alpha,j}^{(i)}(k_{l},u,v),\max_{j=3,4,5,6}\hat{f}_{\alpha,j}^{(i)}(k_{l},u)\right), \  \ u>0 \]


\[ f_{\alpha}^{(i)}(0,u)=\max\left(\max_{j=1,2,v>u}\hat{f}_{\alpha,j}^{(i)}(0,u,v),\max_{3\leq j \leq 8}\hat{f}_{\alpha,j}^{(i)}(0,u)\right), \  \ u>0 \]

\ We have:
\ When , $0\leq k_{l}\leq k_{n}$

\begin{equation}
 S_{k_{l},\alpha}(\mathcal{A};\mathcal{P},z,\ln^{\alpha}(\xi^{2}))\leq Xe^{-\gamma}C(\omega)\ln^{\alpha-1}(\xi^{2})(u+\frac{k_{l}}{\alpha u^{\alpha-1}})F_{\alpha}^{(i)}(k_{l},u)(1+O(\frac{1}{\ln^{\frac{1}{14}}(\xi^{2})}))
\end{equation} 

\[ +\ln^{\alpha}(\xi^{2})\sum_{d|\mathcal{P}(z),d<\xi^{2}}3^{v_{1}(d)}|r_{d}| \]
\ And
\ When , $0\leq k_{l}\leq k_{n+4}$

\begin{equation}
 S_{k_{l},\alpha}(\mathcal{A};\mathcal{P},z,\ln^{\alpha}(\xi^{2}))\geq Xe^{-\gamma}C(\omega)\ln^{\alpha-1}(\xi^{2})(u+\frac{k_{l}}{\alpha u^{\alpha-1}})f_{\alpha}^{(i)}(k_{l},u)(1+O(\frac{1}{\ln^{\frac{1}{14}}(\xi^{2})}))
\end{equation} 

\[ -\ln^{\alpha}(\xi^{2})\sum_{d|\mathcal{P}(z),d<\xi^{2}}3^{v_{1}(d)}|r_{d}| \]

\ $Table 4$ and $Table 5$ are results of the Double Sieve \\

\begin{table}[h]
	\centering
	\caption{$e^{-\gamma}(u+\frac{k_{l}}{2u})F_{2}(k_{l},u)$ and $e^{-\gamma}(u+\frac{k_{l}}{2u})f_{2}(k_{l},u)$ of Double Sieve}
		\begin{tabular}{|c|c|c|c|c|c|c|c|c|c|c|} \hline
			 u & \multicolumn{2}{|c|}{\textbf{3}} & \multicolumn{2}{|c|}{\textbf{2.5}} & \multicolumn{2}{|c|}{\textbf{2}} & \multicolumn{2}{|c|}{\textbf{1.5}} & \multicolumn{2}{|c|}{\textbf{1}} \\ \hline
			 $k_{l}$ & F & f & F & f & F & f & F & f & F & f \\ \hline
			4.5 &   & 1.98678 &  & 1.72368 &  & 1.45786 &  & 0.87942 &  & 0.0 \\ \hline
			4 & 2.12780 & 1.98678 & 1.97904 & 1.72368 & 1.92609 & 1.45786 & 3.45385 & 0.87942 & 7.62429 & 0.0 \\ \hline
			3.75 & 2.08345 & 1.98678 & 1.91941 & 1.72368 & 1.88222 & 1.45787 & 2.66467 & 0.87942 & 3.883690  & 0.0 \\ \hline
			3.5 & 2.06380 & 1.95927 & 1.89666 & 1.68620 & 1.83588 & 1.43309 & 2.50990 & 0.87942 & 3.69096 & 0.0 \\ \hline
			3.25 & 2.04590 & 1.93176 & 1.87680 & 1.64872 & 1.80742 & 1.39642 & 2.35513 & 0.87942 & 3.49502 & 0.0 \\ \hline
			3 & 2.02717 & 1.90425 & 1.85984 & 1.61124 & 1.78239 & 1.33625 & 2.20035 & 0.87942 & 3.29908 & 0.0 \\ \hline
			2.75 & 2.00935 & 1.87674 & 1.84315 & 1.57376 & 1.76010 & 1.27608 & 2.04558 & 0.87942 & 3.10314 & 0.0 \\ \hline
			2.5 & 1.99910 & 1.84923 & 1.82406 & 1.53628 & 1.74140 & 1.21591 & 1.89080 & 0.87942 & 2.90720 & 0.0  \\ \hline
			2.25 & 1.97322 & 1.82172 & 1.81642 & 1.49881 & 1.73603 & 1.15574 & 1.73603 & 0.87942 & 2.71126 & 0.0 \\ \hline
			2 & 1.95545 & 1.79421 & 1.80879 & 1.46133 & 1.73420 & 1.09557 & 1.73420 & 0.78174 & 2.51532 & 0.0 \\ \hline
			1.75 & 1.93768 & 1.76670 & 1.80115 & 1.42385 & 1.73354 & 1.03540 & 1.73354 & 0.68402 & 2.31938 & 0.0 \\ \hline
			1.5 & 1.91990 & 1.73919 & 1.79351 & 1.38638 & 1.73288 & 0.97523 & 1.73288 & 0.58630 & 2.12343 & 0.0 \\ \hline
			1.25 & 1.90213 & 1.71168 & 1.78587 & 1.34889 & 1.73221 & 0.91506 & 1.73221 & 0.48859 & 1.92749 & 0.0 \\ \hline
			1 & 1.88436 & 1.68417 & 1.77823 & 1.31142 & 1.73155 & 0.85489 & 1.73155 & 0.39087 & 1.73155 & 0.0 \\ \hline
			0.75 & 1.86659 & 1.65666 & 1.77059 & 1.27394 & 1.73089 & 0.79471 & 1.73089 & 0.29315 & 1.73089 & 0.0 \\ \hline
			0.5 & 1.84882 & 1.62915 & 1.76295 & 1.23646 & 1.73023 & 0.73454 & 1.73023 & 0.19543 & 1.73023 & 0.0 \\ \hline
			0.25 & 1.83104 & 1.60164 & 1.75531 & 1.19898 & 1.72957 & 0.67437 & 1.72957 & 0.09771 & 1.72957 & 0.0 \\ \hline
			0 & 1.81327 & 1.57413 & 1.74767 & 1.16150 & 1.72891 & 0.61420 & 1.72891 & 0 & 1.72891 & 0 \\ \hline
			
		\end{tabular}

\end{table}

\begin{table}[h]
	\centering
	\caption{Double Sieve $e^{-\gamma}uF_{2}(0,u)$ and $e^{-\gamma}uf_{2}(0,u)$}
		\begin{tabular}{||c|c|c||c|c|c||c|c|c||} \hline
			u & $e^{-\gamma}$uF(u) & $e^{-\gamma}$uf(u) & u & $e^{-\gamma}$uF(u) & $e^{-\gamma}$uf(u) & u & $e^{-\gamma}$uF(u) & $e^{-\gamma}$uf(u) \\ \hline
			5.0 & 2.808880 & 2.805636 & 3.7 & 2.108948 & 2.049196 & 2.4 & 1.747668 & 1.081561 \\ \hline
			4.9 & 2.753155 & 2.749037 & 3.6 & 2.060095 & 1.986305 & 2.3 & 1.747668 & 0.981977 \\ \hline
			4.8 & 2.697544 & 2.692357 & 3.5 & 2.012771 & 1.921903 & 2.2 & 1.746616 & 0.875651 \\ \hline
			4.7 & 2.642076 & 2.635561 & 3.4 & 1.967290 & 1.856157 & 2.1 & 1.740615 & 0.761280 \\ \hline
			4.6 & 2.586792 & 2.578614 & 3.3 & 1.924047 & 1.788645 & 2.0 & 1.728908 & 0.637005 \\ \hline
			4.5 & 2.531744 & 2.521477 & 3.2 & 1.883539 & 1.719210 & 1.9 & 1.728908 & 0.459369 \\ \hline
			4.4 & 2.476986 & 2.464103 & 3.1 & 1.846355 & 1.647753 & 1.8 & 1.728908 & 0.260835 \\ \hline
			4.3 & 2.422583 & 2.406434 & 3.0 & 1.813272 & 1.574131 & 1.702 & 1.728908 & 2.4275E-03  \\ \hline
			4.2 & 2.368612 & 2.348405 & 2.9 & 1.786400 & 1.498156 & 1.6 & 1.728908 & 0.0   \\ \hline
			4.1 & 2.315164 & 2.289931 & 2.8 & 1.768490 & 1.419319 & 1.5 & 1.728908 & 0.0 \\ \hline
			4.0 & 2.262342 & 2.230915 & 2.7 & 1.761952 & 1.337257 & 1.4 & 1.728908 & 0.0 \\ \hline
			3.9 & 2.210264 & 2.171225 & 2.6 & 1.761952 & 1.251518 & 1.3 & 1.728908 & 0.0 \\ \hline
			3.8 & 2.159074 & 2.110724 & 2.5 & 1.747668 & 1.161508 & 1.2 & 1.728908 & 0.0 \\ \hline

		\end{tabular}
	\label{tab:DoubleSieve1}
\end{table}

\ We are now in a position to prove $Theorem 3$ and $Theorem 4$

\begin{proof} of $Teorem 4$ \\

\ Suppose $\xi^{2}=\frac{N^{0.5}}{\ln^{B}(N)}$, $u=\frac{\ln(\xi^{2})}{\ln(N^{0.5})}=1+O(\frac{1}{\ln(N)})$ \\
\ From $Table 4$ or $Table 5$ we have
\[ D(N)\leq S_{0}(\mathcal{A};\mathcal{P},N^{0.5})=\frac{1}{\ln^{2}(\xi^{2})}S_{0,2}(\mathcal{A};\mathcal{P},N^{0.5},\ln^{2}(\xi^{2})) \]
\[ \leq \frac{X2C(N)e^{-\gamma}}{\ln(\xi^{2})}uF_{2}^{(i)}(0,u)(1+o(1))+\ln^{2}(\xi^{2})R=\frac{4NC(N)1.728908}{\ln^{2}(N)}(1+o(1))+\ln^{2}(\xi^{2})R\]
\[ \leq \frac{6.916NC(N)}{\ln^{2}(N)}(1+o(1)) +\ln^{2}(\xi^{2})R\]
\ Where
\[ R=\sum_{d|\mathcal(P)(N^{0.5}),d<\xi^{2}}3^{v_{1}(d)}\left|r_{d}\right|\ll \frac{N}{\ln^{A}(N)} \]
\ We obtain $Theorem 4$. \\
\end{proof}

\begin{proof} of $Theorem 3$
\ Using similar double sieve method as the Jing Run Chen\cite{J. R. Chen1973} proof of $D_{1,2}(N)>0$. Suppose $\xi^{2}=\frac{\sqrt{N}}{\ln^{B}(N)}$, we have: \\
\[ \ln^{2}(\xi^{2})D_{1,2}(N)\geq S_{2.25,2}(\mathcal{A};\mathcal{P},N^{\frac{1}{3}},\ln^{2}(\xi^{2}))-2.25\sum_{q<N^{\frac{1}{3}}}\ln^{2}(q)\sum_{N^{\frac{1}{3}}\leq p<\sqrt{\frac{N}{q}}}S_{0}(\mathcal{A}_{qp},\mathcal{P}(pq),\sqrt{\frac{N}{q}})\]
\[ \geq \ln{\xi^{2}}\frac{N}{\ln(N)}2C(N)e^{-\gamma}(1.5+\frac{2.25}{2\times1.5})f_{2}^{(i)}(2.25,1.5)(1+o(1))-2.25\Omega_{1} \]
\ By $Table 4$ we have 
\[ e^{-\gamma}(1.5+\frac{2.25}{2\times1.5})f_{2}^{(i)}(2.25,1.5)\geq 0.8794 \]
\[ \Omega_{1}\leq \frac{X4C(N)}{\ln(N)}e^{-\gamma}2F_{2}^{(0)}(0,2)(1+o(1)) \leq \frac{X4C(N)}{\ln(N)}\times 1.876677(1+o(1))\]
\ Where 
\[ X=\sum_{q<N^{\frac{1}{3}}}\ln^{2}(q)\sum_{N^{\frac{1}{3}}<p<\sqrt{\frac{N}{q}}\pi(\frac{N}{pq})}(1+o(1)) \]
\[ =\sum_{q<N^{\frac{1}{3}}}\ln^{2}(q)\sum_{N^{\frac{1}{3}}<p<\sqrt{\frac{N}{q}}}\frac{N}{pq\ln(\frac{N}{pq})}(1+o(1)) \]
\[ =\frac{N}{\ln(N)}\int_{2}^{N^{\frac{1}{3}}}\frac{\ln(t)}{t}dt\int_{N^{\frac{1}{3}}}^{(\frac{N}{t})^{\frac{1}{2}}}\frac{ds}{s\ln(s)(1-\frac{\ln(st)}{\ln(N)})}(1+o(1)) \]
\[ =N\ln(N)\int_{0}^{\frac{1}{3}}tdt\int_{\frac{1}{3}}^{\frac{1-t}{2}}\frac{ds}{s(1-s-t)} (1+o(1))\]

\[ =N\ln(N)\int_{0}^{\frac{1}{3}}\frac{t\ln(2-3t)}{(1-t)}dt \leq 0.01846N\ln(N)(1+o(1))\]
\ we obtain
\[ 2.25\Omega_{1}\leq NC(N)(2.25\times 4\times 1.876677 \times 0.01846)(1+o(1)) \leq 0.3118NC(N)(1+o(1)) \]
\[ D_{1,2}(N)\geq \frac{NC(N)}{\ln^{2}(\xi^{2})}\left(\frac{2\times 0.8794\ln(\xi^{2})}{\ln(N)}-0.3118\right)(1+o(1)) \]
\[ =\frac{NC(N)}{\ln^{2}(N)}4(0.8794-0.3118)(1+o(1))\geq 2.27\frac{NC(N)}{\ln^{2}(N)}(1+o(1))\]

\ This completes the proof.
\end{proof}

\section{Part III: Application in estimate of the exception set in Goldbach's number}

\ Define set $\mathcal{B}$ as:
\[ \mathcal{B}:=\left\{n,2|n,\frac{X}{2}\leq n \leq X,n \ \ not \ \ a \ \ Goldbach's \ \ number\right\} \]
\[ E(X):=\left\{n,2|n, 2 \leq n\leq X,n\in \mathcal{B}  \right\} \]
\ This Part will prove:\\
\subsection{Theorem 5}
\ For any small positive number $\epsilon$ we have
\[ \left|E(X)\right| \ll O(X^{0.702+\epsilon})  \]

\ In order to prove $Theorem 5$ we need the follwing two $Lemmas$ \\

\ Suppose $\mathcal{N}$ is any natural number set that fulfils the following expression: \\

\[ n \in \mathcal{N} \Rightarrow 2|n, \frac{X}{2} \leq n \leq X \]
\ And 
\[ \left|\mathcal{N}\right| = X^{1-2\Delta} \]

\ Suppose $\xi^{2}=\frac{X^{1-\Delta}}{\ln^{B}(X)}$, $\mathcal{A}$ is a number set.
\[ \mathcal{A}:=\left\{n-p;\xi^{2}\leq p<\frac{1}{2}X,(p,n)=1,n \in \mathcal{N} \right\} \]

\[ \left|\mathcal{A}\right|=\frac{X^{2-2\Delta}}{2\ln(X)}(1+o(1)) \]
\ When $d\leq \xi^{2}$
 
\begin{equation}
\left|\mathcal{A}_{d}\right|=\sum_{n \in \mathcal{N}}\sum_{\xi^{2}\leq p<\frac{1}{2}X,d|n-p,(d,n)=1}1=\frac{1}{\phi(d)}\sum_{\chi_{d}}\sum_{n \in \mathcal{N}}\bar{\chi}(n)\sum_{\xi^{2}\leq p<\frac{1}{2}X}\chi(p)
\end{equation}

\[ =\frac{1}{\phi(d)}\sum_{\xi^{2}\leq p<\frac{1}{2}X}1\sum_{n \in \mathcal{N},(n,d)=1}1 \]
\[ +\frac{1}{\phi(d)}\sum_{\chi_{d}\neq \chi_{0}}\sum_{n \in \mathcal{N}}\bar{\chi}(n)\sum_{\xi^{2}\leq p<\frac{1}{2}X}\chi(p)\]
\[ =\frac{1}{d}\sum_{\xi^{2}\leq p<\frac{1}{2}X}1\sum_{n \in \mathcal{N}}\omega_{n}(d)+r_{d} \]
\ Where 
\[ \phi(d) = d\prod_{p|d}(1-\frac{1}{p}) \ \ is \ \ the \ \ Euler \ \ function. \]
\ and function $\chi$ is the Dirichlet character, function $\chi_{0}$ is the principal character.

\begin{equation}
 r_{d}=\frac{1}{\phi(d)}\sum_{\chi_{d}\neq \chi_{0}}\sum_{n \in \mathcal{N}}\bar{\chi}(n)\sum_{\xi^{2}\leq p<\frac{1}{2}X}\chi(p)
\label{RD2}
\end{equation}

\[ \omega_{n}(d)=\frac{d}{\phi(d)} \  \ if \ \ (n,d)=1 \ \; or \ \ \omega_{n}(d)=0 \ \ if \ \ (n,d)>1 \]

\subsection{Lemma 3.1}(Estimate the sum of the character with Large Sieve \cite{C.D. Pan and C.B. Pan1992})\\
\ Suppose $Q\geq 2$, $1<D<Q$, we have \\

\[ \sum_{D<q\leq Q}\frac{1}{\phi(q)}\sideset{}{^*}\sum_{\chi_{q}} \left|\sum_{n=M+1}^{M+N}a_{n}\chi(n)\right|^{2} \ll (Q+\frac{N}{D})\sum_{n=M+1}^{M+N}\left|a_{n}\right|^{2} \]

\ Where $\sideset{}{^*}\sum_{\chi_{q}}$ denotes (here and later) summation over all primitive characters mod $q$.
\subsection{Corollary 3.1.1}
\ Suppose $Q\geq 2$, $1<D<Q$, we have \\

\[ \sum_{D<q\leq Q}\frac{3^{v_{1}(q)}}{\phi(q)}\sideset{}{^*}\sum_{\chi_{q}} \left|\sum_{n=M+1}^{M+N}a_{n}\chi(n)\right|^{2} \ll \ln^{A+17}(N)(Q+\frac{N}{D})\sum_{n=M+1}^{M+N}\left|a_{n}\right|^{2} \]


\subsection{Lemma 3.2}(Estimate the sum of the primer number character )\cite{E. Bombieri1965}\cite{C.D. Pan and C.B. Pan1992}\\
\ when $D=\frac{X^{0.5}}{\ln^{B}(X)}$, $B=2A+32$
\[ \sum_{1<d|\mathcal{P}(z),d\leq D}3^{v_{1}(d)}\sideset{}{^*}\sum_{\chi_{d} }\frac{1}{\phi(d)}\left|\sum_{p<n}\chi(p)\right| \ll \frac{X}{ln^{A}(x)}\]

\subsection{Lemma 3.3}
\ Suppose $\xi^{2}=\frac{X^{1-\Delta}}{\ln^{B}(X)}$, $B=2A+32$
\[ E_{\mathcal{A}}=\sum_{d|\mathcal{P}(z),d<\xi^{2}}3^{v_{1}(d)}|r_{d}| \ll \frac{X}{ln^{A}(x)}\]

\begin{proof}
\ Suppose $D=\frac{X^{0.5}}{\ln^{B}(X)}$, $B=2A+32$ we have
\[ E_{\mathcal{A}}=\sum_{d|\mathcal{P}(z),d<\xi^{2}}3^{v_{1}(d)}\sum_{\chi_{d} \neq \chi_{0}}\frac{1}{\phi(d)}\sum_{n \in \mathcal{N}}\overline{\chi}(n)\sum_{p<n}\chi(p) \]
\[ \ll \ln(\xi^{2})\sum_{1<d|\mathcal{P}(z),d<\xi^{2}}3^{v_{1}(d)}\sideset{}{^*}\sum_{\chi_{d} }\frac{1}{\phi(d)}\left|\sum_{n \in \mathcal{N}}\overline{\chi}(n)\right|\left|\sum_{p<n}\chi(p)\right| \]
\[ =\ln(\xi^{2})\sum_{1<d|\mathcal{P}(z),d\leq D}3^{v_{1}(d)}\sideset{}{^*}\sum_{\chi_{d} }\frac{1}{\phi(d)}\left|\sum_{n \in \mathcal{N}}\overline{\chi}(n)\right|\left|\sum_{p<n}\chi(p)\right| \]
\[ +\ln(\xi^{2})\sum_{d|\mathcal{P}(z),D<d\leq \xi^{2}}3^{v_{1}(d)}\sideset{}{^*}\sum_{\chi_{d} }\frac{1}{\phi(d)}\left|\sum_{n \in \mathcal{N}}\overline{\chi}(n)\right|\left|\sum_{p<n}\chi(p)\right|=\Sigma_{1}+\Sigma_{2} \]

\ By $Lemma 3.2$, we obtain
\[ \Sigma_{1} \ll \ln(\xi^{2})\left|\mathcal{N}\right|\frac{X}{\ln^{A}(X)}=\frac{X^{2-2\Delta}}{\ln^{A-1}(X)} \]
\[ (\Sigma_{2})^{2} \leq \ln^{2}(\xi^{2})\sum_{d|\mathcal{P}(z),D<d\leq \xi^{2}}3^{v_{1}(d)}\sideset{}{^*}\sum_{\chi_{d} }\frac{1}{\phi(d)}\left|\sum_{n \in \mathcal{N}}\overline{\chi}(n)\right|^{2} \sum_{d|\mathcal{P}(z),D<d\leq \xi^{2}}3^{v_{1}(d)}\sideset{}{^*}\sum_{\chi_{d} }\frac{1}{\phi(d)}\left|\sum_{p <X}\chi(p)\right|^{2} \]
\[ \leq \ln^{36+2A}(X)(\xi^{2}+\frac{X}{D})^{2}\sum_{n \in \mathcal{N}}1\sum_{p<X}1 \]
\[ \leq \ln^{36+2A}(X) (\xi^{2}+\frac{X}{D})^{2}X^{2-2\Delta} \]
\[ \Sigma_{2}\leq \ln^{18+A}(x)(\xi^{2}+\frac{X}{D})X^{1-\Delta} \]

\ Where $\xi^{2}=X^{1-\Delta}/\ln^{B}(X)>\frac{X}{D}=X^{0.5}\ln^{B}(X)$ 
\ By $Lemma3.1$
\[ \Sigma_{2} \leq \ln^{A+18}(X)\xi^{2}X^{1-\Delta} \leq X^{2-2\Delta}/\ln^{A}(X)\]
\ So we obtain
\[ \left|E_{\mathcal{A}}\right| \leq O\left(\frac{X^{2-2\Delta}}{\ln^{A-1}(X)}\right) \]

\ Where $A$ is a any big integer number. \\
\end{proof}
 
\begin{proof} of $Theorem 5$\\

\ Suppose $\frac{X}{2}\leq n<X$ , $u=\frac{\ln(\xi^{2})}{\ln(z)}$
\[ \mathcal{A}^{(n)}:=\left\{m:m=n-p,\xi^{2}\leq p<X \right\}\]

\[ S_{2,0}(\mathcal{A};\mathcal{P},z,\ln^{2}(\xi^{2}))=\ln(\xi^{2})\sum_{n \in \mathcal{N}}S_{2,0}(\mathcal{A}^{(n)};\mathcal{P},z) \]
\ It's main term is the sum of different $n$, but we can combine it's remainder term as the sum of equation (4.2).
\[ S_{2,0}(\mathcal{A};\mathcal{P},z,\ln^{2}(\xi^{2}))\geq Xe^{-\gamma}\ln(\xi^{2})uf_{2}(0,u)\sum_{n \in \mathcal{N}}2C(n) +\sum_{d|\mathcal{P},d<\xi^{2}}3^{v_{1}(d)}\left|r_{d}\right|\]
\ According to $Lemma3.1$ and $Lemma3.3$, for the remainder term we have
\[ \sum_{d|\mathcal{P},d<\xi^{2}}3^{v_{1}(d)}\left|r_{d}\right|\leq \frac{X^{2-2\Delta}}{\ln^{A-1}(X)} \]
\ In $Theorem 2$ we take:
\[ z=X^{\frac{1}{2}}, \ \ u=\frac{\ln(\xi^{2})}{\ln(z)}=2-2\Delta+O(\frac{\ln(\ln(X))}{\ln(X)}) \]
\ Since 
\[ n \in \mathcal{N}\Rightarrow n-p \in \mathcal{A}\Rightarrow n-p\leq X\]

\[ (n-p,\mathcal{P}(z))=(n-p,\mathcal{P}(X^{\frac{1}{2}}))=1\Rightarrow n-p \ \ is \ \ a \ \ primer \ \ number. \]
\ and for any even number $n$
\[ C(n)\geq C(2)>0 \]
\[ \sum_{n \in \mathcal{N}}C(n) \geq C(2)\left|\mathcal{N}\right|= C(2)X^{1-2\Delta} \]
\ So we have:
\[ \ln^{2}(\xi^{2})\sum_{n \in \mathcal{A},(n,\mathcal{P}(z))=1}1=S_{2,0}(\mathcal{A},\mathcal{P},z,\ln^{2}(\xi^{2}))\geq \ln(\xi^{2})uf_{2}(0,u)\frac{X^{2-2\Delta}}{\ln(n)}2e^{-\gamma}C(2)(1+0(1)) \]

\ Set $\epsilon$ is any smal positive number, when 
\[ \Delta \leq 0.149-\frac{\epsilon}{2} \]
\ For sufficiently large number $X$
\[ u=2-2\times 0.149+\epsilon-O(\frac{\ln(\ln(X))}{\ln(X)}) \geq 1.702 \]
\[ f_{2}(0,u)\geq f_{2}(0,1.702)>0 \]

\ So for any small positive number $\epsilon$ we obtain
\[ \left|\mathcal{B}\right|\ll X^{0.702+\epsilon} \]

\ From this formula we know: in region $(\frac{X}{2},X)$ non Goldbach number less then $X^{0.702+\epsilon}$.
\ Finally we obtain
\[ \left|E(X)\right| \leq \sum_{2^{i}\leq X}O(\frac{X}{2^{i}})^{0.702+\epsilon} \ll O(X^{0.702+\epsilon}) \]
\ This completes the proof
\end{proof}

\section{Discuss}

1.	The process of proving the Jurkat-Richert Theorem is iterative and determined by the structure of the decomposed Sieve function.  The expansion of the Sieve function in Part I, changed its decomposed structure such that weighted Sieve functions can be used in the iterative operation.  The relationships between Sieve functions with different parameters make it possible to use the expanded Sieve functions to improve the traditional Sieve functions.  Lemma 1.21 is introduced for the estimation of the upper limit.  There is great improvement in determining the upper and lower limit of the Sieve function by using this Lemma 1.21.  These advantages are more greatly seen in the Double Sieve in Part II. \\

\ 2.	Part II Lemma 2.3 uses the weak condition, $w(u)\leq 1$  , in the analysis of the Double Sieve.  The results have a large margin but could be improved through further detailed analysis of $w(u)$  with different parameters $\xi^{2}$.  \\

\ 3.	In Part III, Theorem 5 uses the approved Sieve methods to develop a new estimation of the exception set in Goldbach's number that is better than using the Circle Method.  If we change the number set $\mathcal{N}$ to 
\[ \mathcal{N}:=\left\{n;n=N-p,p<\frac{N}{2}\right\}, \]
\ Where N is a Sufficiently large odd number. \ 
\[ \left|\mathcal{N}\right|=\frac{N}{2\ln(N)}(1+o(\frac{1}{\ln(N)}) \geq N^{0.702} \]
\ and now
\[ \mathcal{A}:=\left\{n-p_{2};\xi^{2}\leq p_{2}<\frac{1}{2}X,(p_{2},n)=1,n=N-p_{1},p_{1}<\frac{N}{2} \right\}, \]

\ that actually gives a new proof of the Goldbach Conjecture about odd numbers that is different than the circle method.    \\

\ 4.	Theorem 2 also inferred that the Goldbach Conjecture and the Twin primes Conjecture will hold true if the parameters $D$ in Lemma 2.1 are increased to more than $N^{0.851}$ .  There will be less restrictions on the parameter, if the parameter   on Lemma 2.2 is improved as well.


\begin{thebibliography}{9}                                                                                                
\bibitem {V. Brun 1915}V. Brun. \textit{Uber das Goldbachsche Gesetz und die Anzahl der Primzahlpaare.} Arch. Mat. Natur. B, 34, no. 8, 1915.
\bibitem {V. Brun 1920}V. Brun. \textit{Le crible d'Eratosth`ene et le th'eor`eme de Goldbach.} Videnskaps. Skr. Mat. Natur. Kl. Kristiana, no. 3, 1920.
\bibitem {Yu. V. Linnik 1941}Yu. V. Linnik. \textit{The large sieve. C.R. Acad.} Sci. URSS (N.S.), 30 (1941), 292-294.
\bibitem {A. Selberg 1947}A. Selberg. \textit{On an elementary method in the theory of primes. Norske Vid. Selsk. Forh.}, Trondhjem, 19 (1947), 64-67.
\bibitem {A. Selberg 1949}A. Selberg. \textit{On elementary methods in prime number theory and their limitations.} 11th. Skand. Math. Kongr., Trondhjem, (1949), 13-22.
\bibitem {A. Selberg1950}A. Selberg. \textit{The general sieve-method and its place in prime number theory.} Proc. Intern. Cong. Math., Cambridge, Mass., 1 (1950), 286-292.
\bibitem {H. Iwaniec1980}H. Iwaniec. \textit{Rosser's sieve.} Acta Arith., 36 (1980), 171-202.
\bibitem {W.B. Jurkat and H.-E. Richert1965}W.B. Jurkat and H.-E. Richert. \textit{An improvement of Selberg's sieve method. I.} Acta Arith., 11 (1965), 217-240.
\bibitem {E. Bombieri1965}E. Bombieri. \textit{On the large sieve} Mathematika, 12 (1965), 201-225.

\bibitem {J. R. Chen1973}\textsc{J. R. Chen}:\ \textit{On the Representation of a Large Even Integer as the Sum of a Prime and the Product of at Most Two Primes.}, Sci,Sin. 17 (1973) 157-176.

\bibitem {D. R. Heath-Brown}D. R. Heath-Brown, J. C. Puchta, \textit{Integers represented as a sum of primes and powers of two.}, The Asian Journal of Mathematics, 6 (2002), no. 3, pages 535-565.

\bibitem {H.L. Montgomery1975}H.L. Montgomery, Vaughan, R. C., \textit{The exceptional set in Goldbach's problem.}, Acta Arith. 27 (1975), 353-370.

\bibitem {Halberstam1974} \textsc{Halberstam H, Richert H E.}:\ \textit{Sieve Methods.}, Academic Press \textbf{1974}.

\bibitem {C.D. Pan1981} \textsc{C. D. Pan}:\ \textit{A new mean value theorem and its applications.}, Recent Progress in Analytic Number Theory I, Academic Press, \textbf(1981:) 275-287.

\bibitem {Estermann, T.1938} \textsc{Estermann, T.}:\ \textit{On Goldbach's Problem: Proof that Almost All Even Positive Integers are Sums of Two Primes.}, Proc. London Math. Soc. Ser. 2 44, \textbf(1938:) 307-314.

\bibitem {Vinogradov, I. M.1937} \textsc{Vinogradov, I. M.}:\ \textit{Representation of an Odd Number as a Sum of Three Primes.}, Comptes rendus (Doklady) de l'Academie des Sciences de l'U.R.S.S. 15, \textbf(1937a:) 169-172.

\bibitem {J.-M. Deshouillers1997} \textsc{J.-M. Deshouillers; G. Effinger; H. te Riele; D. Zinoviev}:\ \textit{A complete Vinogradov 3-primes theorem under the Riemann hypothesis}, Electron. Res. Announc. Amer. Math. Soc. 3 \textbf(1997:) 99-104.

\bibitem {J.R. Chen1978} \textsc{J.R. Chen 1978}:\ \textit{On the Goldbach's problem and the sieve methods}, Sci. Sin, 21 \textbf(1978),701-739.

\bibitem {C.D. Pan and C.B. Pan1992} \textsc{C.D. Pan and C.B. Pan}:\ \textit{Goldbach Conjecture}, Science Press, Beijing, China, \textbf(1992).

\bibitem {Y. Motohashi2005} \textsc{Y. Motohashi}:\ \textit{An Overview of the Sieve Method and its History}, math.NT/0505521, \textbf(2005).

\end{thebibliography}
\end{document}